\documentclass[12pt]{article}

\usepackage{enumitem}

\usepackage{amsmath,amsthm,amssymb,amscd}
\usepackage[pdftex]{graphicx}
\usepackage{graphicx}
\usepackage{epsfig}
\usepackage{verbatim}
\usepackage{footmisc}
\usepackage{bbm,bm}
\usepackage{lscape}
\usepackage[longnamesfirst,round]{natbib}
\usepackage{color}
\usepackage{float}
\usepackage{dsfont}
\usepackage{longtable}
\usepackage{multirow}
\usepackage{xr}

\allowdisplaybreaks 

\def\eps{\varepsilon}

\newcommand{\Y}{\mathcal{Y}}

\def\a{\alpha}
\def\N{\mathbb{N}}

\def\R{\mathbb{R}}

\def\cx{\mathcal{X}}
\def\cy{\mathcal{Y}}
\def\cyw{\mathcal{Y}_{w}}

\def\T{\Theta}

\def\sjn{\sum_{j=1}^n}
\def\skn{\sum_{k=1}^n}

\def\beq{\begin{eqnarray*}}
\def\eeq{\end{eqnarray*}}

\def\Lt{\Lambda_\theta}

\def\t{\theta}

\def\H{\mathcal{H}}

\newtheorem{theo}{Theorem}[section]
\newtheorem{lemma}[theo]{Lemma}

\newtheorem{rem}[theo]{Remark}

\newtheorem{alg}[theo]{Algorithm}

\makeatletter

\@addtoreset{equation}{section}
\makeatother

\renewcommand{\baselinestretch}{1.2}

\usepackage{color}

\begin{document}

\title{\bf Specification testing in semi-parametric transformation models}


\author{{\sc Nick Kloodt}{\small$^1$}, {\sc Natalie Neumeyer}{\small$^1$} and {\sc Ingrid Van Keilegom}{\small$^2$}\\\small $^1$Department of Mathematics, University of Hamburg\\\small $^2$Research Centre for Operations Research and Statistics
(ORSTAT), KU Leuven}

\maketitle

\renewcommand{\baselinestretch}{1.1}

\begin{abstract}
In transformation regression models the response is transformed before fitting a regression model to covariates and transformed response. We assume such a model where the errors are independent from the covariates and the regression function is modeled nonparametrically. We suggest a test for goodness-of-fit of a parametric transformation class based on a distance between a nonparametric transformation estimator and the parametric class. We present asymptotic theory under the null hypothesis of validity of the semi-parametric model and under local alternatives. A bootstrap algorithm is suggested in order to apply the test.  We also consider relevant hypotheses to distinguish between large and small distances of the parametric transformation class to the `true' transformation. 
\end{abstract}

\vspace*{.5cm}



\section{Introduction}  \label{introduction}

It is very common in applications to transform data before investigation of functional dependence of variables by regression models. The aim of the transformation is to obtain a simpler model, e.g.\ with a specific structure of the regression function, or a homoscedastic instead of a heteroscedastic model.  Typically flexible parametric classes of transformations are considered from which a suitable one is selected data-dependently.  A classical example is the class of Box-Cox power transformations (see \citet{BC1964}).
For purely  parametric transformation models see \citet{CR1988} and references therein.
\citet{Pow:1991:Proc} and \citet{MH2007} consider transformation quantile regression models. 
Nonparametric estimation of the transformation in the context of parametric regression models has been considered by \citet{Hor1996} and \citet{Chen2002},  among others. 
\citet{Hor2009} reviews estimation in transformation models with parametric regression in the cases where either the transformation or the error distribution or both are modeled nonparametrically. 
\citet{LSvK2008} suggest a profile likelihood estimator for a parametric class of transformations, while the error distribution is estimated nonparametrically and the regression function semi-parametrically.
\citet{HSvK2015} suggest an estimator for the error distribution in the same model.
\citet{NNvK2016} consider profile likelihood estimation in heteroscedastic semi-parametric transformation regression models, i.e.\ the mean and variance function are modeled nonparametrically, while the transformation function is chosen from a parametric class.  A completely nonparametric (homoscedastic) model is considered by \citet{CKC2015}. Their approach was modified and corrected by \citet{CvK2018}. The version of the nonparametric transformation estimator considered in the latter paper was then applied by \citet{CvK2018a} to suggest a new estimator of the transformation parameter if it is assumed that the transformation belongs to a parametric class. 

In general asymptotic theory for nonparametric transformation estimators is sophisticated and parametric transformation estimators show much better performance if the parametric model is true. A parametric transformation will thus lead to better estimates of the regression function. Moreover, parametric transformations are easier to interpret and allow for subsequent inference in the transformation model. For the latter purpose note that for transformation models with parametric transformation lack-of-fit tests for the regression function as well as tests for significance for covariate components have been suggested by \cite{CvK2016}, \cite{CvK2017}, \cite{AHM2018} and \cite{KloNeu2017}. Those tests cannot straightforwardly be generalized to nonparametric transformation models because known estimators in that model do not allow for uniform rates of convergence over the whole real line, see  \citet{CKC2015} and \citet{CvK2018}.  

However, before applying a transformation model with parametric transformation it would be appropriate to test the goodness-of-fit of the parametric transformation class. In the context of parametric quantile regression, \citet{MH2007} suggest such a goodness-of-fit test.  In the context of nonparametric mean regression \citet{NNvK2016} develop a goodness-of-fit test for the parametric transformation class based on an empirical independence process of pairs of residuals and covariates. The latter approach was modified by \citet{HMNP2018}, who applied empirical characteristic functions. In a linear regression model with transformation of the response \citet{Syd2017} suggests a goodness-of-fit test for the parametric transformation class that is based on a distance between the nonparametric transformation estimator considered by \citet{Chen2002} and the parametric class.  We will follow a similar approach but consider a nonparametric regression model. The aim of the transformations we consider is to induce independence between errors and covariates. The null hypothesis is that the unknown transformation belongs to a parametric class. Note that when applied to the special case of a class of transformations that contains as only element the identity, our test provides indication on whether a classical homoscedastic regression model (without transformation) is appropriate or whether first the response should be transformed. Our test statistic is based on a minimum distance between a nonparametric transformation and the parametric transformations. We present the asymptotic distribution of the test statistic under the null hypothesis of a parametric transformation  and under local alternatives of $n^{-1/2}$-rate. Under the null hypothesis the limit distribution is that of a degenerate U-statistic. With a flexible parametric class applying an appropriate transformation can reduce the dependence enormously, even if the `true' transformation does not belong to the class. Thus, for the first time in the context of transformation goodness-of-fit tests we consider testing for so-called precise or relevant hypotheses. Here the null hypothesis is that the distance between the true transformation and the parametric class is large. If this hypothesis is rejected, then the model with the parametric transformation fits well enough to be considered for further inference. Under the new null hypothesis the test statistic is asymptotically normally distributed. The term ``precise hypotheses" refers to \citet{BD1987}. \citet{DKV2018} considered precise hypotheses in the context of comparing mean functions in the context of functional time series. Note that the idea of precise hypotheses is related to that of equivalence tests, which originate from the field of pharmacokinetics (see \citet{Lak2017}). Throughout we assume that the nonparametric transformation estimator fulfills an asymptotic linear expansion. It is then shown that the estimator considered by \citet{CvK2018} fulfills this expansion and thus can be used for evaluating the test statistic. 

The remainder of the paper is organized as follows. In Section 2 we present the model and the test statistic. Asymptotic distributions under the null hypothesis of a parametric transformation class and under local alternatives are presented in Section 3, which also contains a consistency result and asymptotic results under relevant hypotheses. Section 4 presents a bootstrap algorithm and a simulation study.  Appendix \ref{assumptions} contains assumptions, while  Appendix \ref{estimation-of-h} treats a specific nonparametric transformation estimator and shows that it fulfills the required conditions. 
The proofs of the main results are given in Appendix \ref{proofs}. 
A supplement contains a rigorous treatment of bootstrap asymptotics.

\section{The model and test statistic} 

Assume we have observed $(X_i,Y_i)$, $i=1,\dots,n$, which are independent with the same distribution as $(X,Y)$ that fulfill the transformation regression model 
\begin{equation}\label{modeleq}
h(Y)=g(X)+\eps,
\end{equation}
where $E[\eps]=0$ holds and $\eps$ is independent of the covariate $X$, which is $\R^{d_X}$-valued, while $Y$ is univariate. The regression function $g$ will be modelled nonparametrically. The transformation $h:\R\to\R$ is strictly increasing. Throughout we assume that, given the joint distribution of $(X,Y)$ and some identification conditions, there exists a unique transformation $h$ such that this model is fulfilled. It then follows that the other model components are identified via $g(x)=E[h(Y)|X=x]$ and $\eps=h(Y)-g(X)$. See \citet{CKC2015} for conditions under which the identifiability of $h$ holds. In particular conditions are required to fix location and scale and we will assume throughout that
\begin{equation}\label{idcond}
h(0)=0\quad\textup{and}\quad h(1)=1.
\end{equation}
Now let $\{\Lt:\t\in\T\}$ be a class of strictly increasing parametric transformation functions $\Lt:\R\to\R$, where $\Theta\subseteq\mathbb{R}^{d_{\Theta}}$ is a finite dimensional parameter space. Our purpose is to test whether a semi-parametric transformation model holds, i.e.
\begin{equation*}
\Lambda_{\theta_0}(Y)=\tilde g(X)+\tilde\eps,
\end{equation*}
for some parameter $\theta_0\in\Theta$, where $\tilde\eps$ and $X$ are independent. Due to the assumed uniqueness of the transformation $h$ one obtains $h=h_0$ under validity of the semi-parametric model, where
\begin{equation*}
h_0(\cdot)=\frac{\Lambda_{\t_0}(\cdot)-\Lambda_{\t_0}(0)}{\Lambda_{\t_0}(1)-\Lambda_{\t_0}(0)}.
\end{equation*}
Thus we can write the null hypothesis as 
\begin{equation}\label{null}
H_0:\ h\in\bigg\{\frac{\Lambda_{\t}(\cdot)-\Lambda_{\t}(0)}{\Lambda_{\t}(1)-\Lambda_{\t}(0)}:\t\in\T\bigg\}
\end{equation}
which thanks to (\ref{idcond}) can be formulated equivalently as
\begin{equation}\label{nullhyp-2}
H_0:\ h\in\bigg\{\frac{\Lambda_{\t}(\cdot)-c_2}{c_1}:\t\in\T,c_1\in\R^+,c_2\in\R\bigg\}.
\end{equation}
Our test statistics will be based on the following $L^2$-distance
\begin{eqnarray}
d(\Lt,h)&=&\underset{c_1\in\R^+,c_2\in\R}{\min}\,E\big[w(Y)\{h(Y)c_1+c_2-\Lt(Y)\}^2\big],\label{defd}
\end{eqnarray}
where $w$ is a positive weight function with compact support $\cyw$. Its empirical counterpart is
$$d_n(\Lt,\hat{h}):=\underset{c_1\in C_1,c_2\in C_2}{\min}\,\frac{1}{n}\sjn w(Y_j)\{\hat{h}(Y_j)c_1+c_2-\Lt(Y_j)\}^2,$$
where $\hat h$ denotes a nonparametric estimator of the true transformation $h$ as discussed below, and $C_1\subset\R^+$, $C_2\subset\R$ are compact sets. Assumption \ref{A6} in Appendix \ref{assumptions} assures that the sets are large enough to contain the true values. 
The test statistic is defined as
\begin{eqnarray}
T_n=\min_{\t\in\Theta}d_n(\Lt,\hat{h})\label{teststatistic}
\end{eqnarray}
and the null hypothesis should be rejected for large values of the test statistic. We will derive the asymptotic distribution under the null hypothesis and local and fixed alternatives in Section \ref{asympt} and suggest a bootstrap version of the tests in Section \ref{simulations}.

\begin{rem} \citet{CvK2018} consider the estimator  
$$\hat{\t}:=\arg\underset{\t\in\T}{\min}\,d_n(\Lt,\hat{h})$$
for the parametric transformation (assuming $H_0$) and observe that  $\hat{\t}$ outperforms the version without minimization over $c_1,c_2$, i.e.\ $\tilde{\t}= \arg\min_{\t\in\T}n^{-1}\sjn w(Y_j)[\hat{h}(Y_j)\{\Lt(1)-\Lt(0)\}+\Lt(0)-\Lt(Y_j)]^2$ in simulations. 
%
\end{rem}

Nonparametric estimation of the transformation $h$ has been considered by \citet{CKC2015} and \citet{CvK2018}. For our main asymptotic results we need that $\hat h$ has a linear expansion, not only under the null hypothesis, but also under fixed alternatives and the local alternatives as defined in the next section. The linear expansion should have the form 
\begin{equation}\label{exprh}
	\hat{h}(y)-h(y)=\frac{1}{n}\sum_{i=1}^n\psi(Z_i,\mathcal{T}(y))+o_P(n^{-1/2})\mbox{ uniformly in }y\in\cyw. 
\end{equation}
 Here, $\psi$ needs to fulfil condition \ref{A8} in Appendix \ref{assumptions}
and we use the definitions  ($i=1,\dots,n$)
\begin{eqnarray}\label{U_i}\label{Z_i}
 Z_i&=&(U_i,X_i),\quad U_i\;=\;\mathcal{T}(Y_i),\quad \mathcal{T}(y)=\frac{F_Y(y)-F_Y(0)}{F_Y(1)-F_Y(0)},
 \end{eqnarray}
where $F_Y$ denotes the distribution of $Y$ and is assumed to be strictly increasing on the support of $Y$. To ensure that $\mathcal{T}$ is well defined the values $0$ and $1$ are w.l.o.g. assumed to belong to the support of $Y$, but can be replaced by arbitrary values $a<b\in\mathbb{R}$ (in the support of $Y$). The expansion (\ref{exprh}) could also be formulated with a linear term $n^{-1/2}\sum_{i=1}^n\tilde\psi(X_i,Y_i,y)$.
In Appendix \ref{estimation-of-h} we reproduce the definition of the estimator $\hat h$ that was suggested by  \citet{CvK2018} as modification of the estimator by \citet{CKC2015}. We give regularity assumptions under which the desired expansion holds, see Lemma \ref{expansion}. Other nonparametric estimators for the transformation that fulfill the expansion could be applied as well.

\section{Asymptotic results} \label{asympt}

In this section we will derive the asymptotic distribution under the null hypothesis and under local and fixed alternatives. For the formulation of the local alternatives consider the null hypothesis as given in (\ref{nullhyp-2}), i.e.\ $h(\cdot)c_1+c_2=\Lambda_{\t_0}(\cdot)$ for some $\t_0\in\Theta$, $c_1\in\R^+$, $c_2\in\R$, and instead assume 
$$H_{1,n}:
h(\cdot)c_1+c_2=\Lambda_{\t_0}(\cdot)+n^{-1/2}r(\cdot) \mbox{ for some }\t_0\in\Theta,c_1\in\mathbb{R}^+,c_2\in\mathbb{R}\mbox{ and some function } r.
$$
Due to the identifiability conditions (\ref{idcond}) one obtains $c_2=\Lambda_{\t_0}(0)+n^{-1/2}r(0)$ and $c_1=\Lambda_{\t_0}(1)-\Lambda_{\t_0}(0)+n^{-1/2}(r(1)-r(0))$. Assumption \ref{A5} yields boundedness of $r$, so that we rewrite the local alternative as
\begin{eqnarray}\nonumber
 h(\cdot)
 &=&\frac{\Lambda_{\t_0}(\cdot)-\Lambda_{\t_0}(0)+n^{-1/2}(r(\cdot)-r(0))}{\Lambda_{\t_0}(1)-\Lambda_{\t_0}(0) +n^{-1/2}(r(1)-r(0))}\\
 &=& h_0(\cdot)+n^{-1/2}r_0(\cdot)+o(n^{-1/2}),\label{h0r0}
\end{eqnarray}
where $h_0(\cdot)=(\Lambda_{\t_0}(\cdot)-\Lambda_{\t_0}(0))/(\Lambda_{\t_0}(1)-\Lambda_{\t_0}(0))$ and
\begin{align*}
r_0(\cdot)
&=\frac{r(\cdot)-r(0)-h_0(\cdot)(r(1)-r(0))}{\Lambda_{\theta_0}(1)-\Lambda_{\theta_0}(0)}.
\end{align*}
Note that the null hypothesis $H_0$ is included in the local alternative $H_{1,n}$ by considering $r\equiv 0$ which gives $h=h_0$. 
We assume the following data generating model under the local alternative $H_{1,n}$. Let the regression function $g$, the errors $\eps_i$ and the covariates $X_i$ be independent of $n$ and define $Y_i=h^{-1}(g(X_i)+\eps_i)$ ($i=1,\dots,n$), which under local alternatives depends on $n$ through the transformation $h$. 
Throughout we use the notation ($i=1,\dots,n$)
\begin{eqnarray}\label{S_i}
S_i=h(Y_i)=g(X_i)+\eps_i. 
\end{eqnarray}
Further, recall the definition of  $U_i$ in (\ref{U_i}). Note that the distribution of $U_i$ does not depend on $n$, even under local alternatives, because $F_Y(Y_i)$ is uniformly distributed on $[0,1]$, while 
$F_Y(0)=P(Y_i\leq 0)=P(h(Y_i)\leq h(0))=P(S_i\leq 0)$ due to (\ref{idcond}), and similarly $F_Y(1)=P(S_i\leq 1)$.

To formulate our main result we need some more notations. 	
 For notational convenience, define
$\gamma=(c_1,c_2,\theta)\in \Upsilon:=C_1\times C_2\times\T$, which is assumed to be compact (see \ref{A1} in Appendix \ref{assumptions}). Then, note that
$$T_n=\underset{\gamma=(c_1,c_2,\t)\in\Upsilon}{\min}\,\sjn w(Y_j)\{\hat{h}(Y_j)c_1+c_2-\Lt(Y_j)\}^2.$$
Further, with $Z_i$ from (\ref{Z_i}) and $S_i$ from (\ref{S_i}) define ($i=1,\dots,n$)
\begin{eqnarray}
\dot{\Lambda}_{\t}(y)&=&\bigg(\frac{\partial}{\partial \theta_k}\Lambda_{\theta}(y)\bigg)_{k=1,...,d_{\Theta}}\nonumber\\
R(s)&=& (s,1,-\dot{\Lambda}_{\t_0}(h_0^{-1}(s)))^t \label{R()}\\
\Gamma_0&=&E[w(h_0^{-1}(S_1))R(S_1)R(S_1)^t] \label{Gamma}\\
\varphi(z)&=& E[w(h_0^{-1}(S_2))\psi(Z_1,U_2)R(S_2)\mid Z_1=z] \label{phi()}\\
\zeta(z_1,z_2)&=& E\Big[w(h_0^{-1}(S_3))\big\{\psi(Z_1,U_3)-\varphi(Z_1)^t\Gamma_0^{-1}R(S_3)\big\} \nonumber\\
&&\qquad\times\big\{\psi(Z_2,U_3)-\varphi(Z_2)^t\Gamma_0^{-1}R(S_3)\big\}\mid Z_1=z_1,Z_2=z_2\Big] \label{zeta}\\
\bar{r}(s)&=&r_0(h_0^{-1}(s))-E[w(h_0^{-1}(S_1))r_0(h_0^{-1}(S_1))R(S_1)]^t\Gamma_0^{-1}R(s) \label{bar-r}\\
\tilde{\zeta}(z_1) &=& 2E[w(h_0^{-1}(S_2))\psi(Z_1,U_2)\bar{r}(S_2)\mid Z_1=z] \label{zeta-tilde}
\end{eqnarray}
and let $P_Z$ and $F^Z$ denote the law and distribution function, respectively, of $Z_i$.

\medskip

\begin{theo}\label{theolocalt}
Assume \ref{A1}--\ref{A8} given in Appendix \ref{assumptions}. Let $(\lambda_{k})_{k\in\{1,2,\dots\}}$ be  the eigenvalues of the operator 
$$K\rho(z_1):=\int\rho(z_2)\zeta(z_1,z_2)\,dF_{Z}(z_2)$$
 with corresponding  eigenfunctions $(\rho_{k})_{k\in\{1,2,\dots\}}$, which are orthonormal in the $L^2$-space corresponding to the distribution $F_Z$. Let $(W_k)_{k\in\{1,2,\dots\}}$ be independent and standard normally distributed random variables and let $W_0$ be centred normally distributed with variance $E[(\tilde\zeta(Z_1))^2]$ such that for all $K\in\mathbb{N}$ the random vector
	$(W_0,W_1,\dots,W_K)^t$
	follows a multivariate normal distribution with $\operatorname{Cov}(W_0,W_k)=E[\tilde{\zeta}(Z_1)\rho_k(Z_1)]$  for all $k=1,\dots,K$. Then, under the local alternative $H_{1,n}$, $T_n$ converges in distribution to  
	\begin{equation*}
	(\Lambda_{\t_0}(1)-\Lambda_{\t_0}(0))^2\Bigg(\sum_{k=1}^\infty\lambda_{k}W_k^2+W_0+ E\left[w(h_0^{-1}(S_1))\bar{r}(S_1)^2\right]\Bigg).
	\end{equation*}
	In particular, under $H_0$ (i.e. for $r\equiv 0$), $T_n$ converges in distribution to 
	\begin{equation*}
	T=(\Lambda_{\t_0}(1)-\Lambda_{\t_0}(0))^2\sum_{k=1}^\infty\lambda_{k}W_k^2.
	\end{equation*}
\end{theo}

\medskip

The proof is given in  Appendix \ref{proofs}. An asymptotic level-$\alpha$ test should reject $H_0$ if $T_n$ is larger than the $(1-\alpha)$-quantile of the distribution of $T$. As the distribution of $T$ depends in a complicated way on unknown quantities, we will propose a bootstrap procedure in Section \ref{simulations}.

\begin{rem}\label{rem-cov}
Note that $\zeta(z_1,z_2)=E[I(z_1)I(z_2)]$ with 
$$	I(z):=w(h_0^{-1}(S_1))^{1/2}\left(\psi(z,U_1)-\varphi(z)^t\Gamma_0^{-1}R(S_1)\right).$$
Thus, the operator $K$ defined in Theorem \ref{theolocalt} is positive semi-definite.
\end{rem}

Next we consider fixed alternatives of a transformation $h$ that do not belong to the parametric class, i.\,e.\
\begin{equation*}
H_1:\quad d(h,\Lambda_{\t})>0\quad\textup{for all }\t\in\T.
\end{equation*}

\medskip

\begin{theo}\label{theoconsistency}
	Assume \ref{A1}--\ref{A4}, \ref{A6'} given in the Appendix
	 and let $\hat{h}$ estimate $h$ uniformly consistently on compact sets. Then, under $H_1$, $\lim_{n\to\infty}P(T_n>q)=1$ for all $q\in\R$, that is, the proposed test is consistent.
\end{theo}

\medskip
\noindent
The proof is given in  Appendix \ref{proofs}.

The transformation model with a parametric transformation class might be useful in applications even if the model does not hold exactly. With a good choice of $\theta$ applying the transformation $\Lambda_\t$ can reduce the dependence between covariates and errors enormously. Estimating an appropriate $\t$ is much easier than estimating the transformation $h$ nonparametrically. Consequently, one might prefer the semiparametric transformation model over a completely nonparametric one. It is then of interest how far away we are from the true model. Therefore, in the following we consider testing precise hypotheses (relevant hypotheses)
\begin{equation*}
H'_0:\underset{\t\in\T}{\min}\,d(h,\Lambda_{\t})\geq\eta\quad\textup{and}\quad H'_1:\underset{\t\in\T}{\min}\,d(h,\Lambda_{\t})<\eta.
\end{equation*}
If a suitable test rejects $H_0'$ for some small $\eta$ (fixed beforehand by the experimenter) the model is considered ``good enough'' to work with, even if it does not hold exactly. To test those hypotheses we will use the same test statistic as before, but we have to standardize differently. Assume $H_0'$, then $h$ is a transformation which does not belong to the parametric class, i.e.\ the former fixed alternative $H_1$ holds.  Let
$$M(\gamma)=M(c_1,c_2,\t)=E\{w(Y)(h(Y)c_1+c_2-\Lambda_{\t}(Y))^2\},$$
and let
$$\gamma_0=(c_{1,0},c_{2,0},\t_0):=\arg\underset{(c_1,c_2,\t)\in\Upsilon}{\min}\,M(c_1,c_2,\t).$$
Note that $\underset{c_1\in C_1,c_2\in C_2}{\min}\,M(\gamma)=d(\Lt,h)$ for all $\theta\in\Theta$. Assume that
\begin{equation}\label{Gamma2}
\Gamma'=E\left[w(Y_1)
\left(\begin{array}{ccc}h(Y_1)^2&h(Y_1)&-h(Y_1)\dot{\Lambda}_{\t_0}(Y_1)
\\h(Y_1)&1&-\dot{\Lambda}_{\t_0}(Y_1)
\\-h(Y_1)\dot{\Lambda}_{\t_0}(Y_1)^t&-\dot{\Lambda}_{\t_0}(Y_1)^t&\Gamma'_{3,3}\end{array}\right)\right]
\end{equation}
is positive definite, where
$\Gamma'_{3,3}=\dot{\Lambda}_{\t_0}(Y_1)^t\dot{\Lambda}_{\t_0}(Y_1)-\ddot{\Lambda}_{\t_0}(Y_1)\tilde R_1$
with 
$$\ddot{\Lambda}_{\t}(y)=\bigg(\frac{\partial^2}{\partial \theta_k\partial \theta_\ell}\Lambda_{\theta}(y)\bigg)_{k,\ell=1,...,d_{\Theta}}
$$
and $\tilde R_i=h(Y_i)c_{1,0}+c_{2,0}-\Lambda_{\t_0}(Y_i)$ $ (i=1,\dots,n)$. 

\medskip

\begin{theo}\label{theointerchangedhyp}
	Assume \ref{A1}--\ref{A4}, \ref{A6'}, \ref{A8'}, let \ref{A7} hold with $\gamma_0$ from \ref{A6'}
	 and let $\Gamma'$ be positive definite. Then
	$$n^{1/2}(T_n/n-M(\gamma_0))\overset{\mathcal{D}}{\rightarrow}\mathcal{N}\big(0,\sigma^2\big)$$
	with $\sigma^2=\operatorname{Var}\big(w(Y_1)\tilde R_1^2+\delta(Z_1)\big)$,
	where $\delta(Z_1)=2E[w(Y_2)\psi(Z_1,U_2)\tilde R_2\mid Z_1]$.
\end{theo}
\medskip
The proof is given in  Appendix \ref{proofs}.
A consistent asymptotic level-$\alpha$-test rejects $H'_0$ if $(T_n-n\eta)/(n\hat{\sigma}^2)^{1/2}<u_{\alpha}$, where $u_{\alpha}$ is the $\alpha$-quantile of the standard normal distribution and $\hat{\sigma}^2$ is a consistent estimator for $\sigma^2$. Further research is required on suitable estimators for $\sigma^2$. For some intermediate sequence $m=m_n$ we considered
\begin{align*}
\hat{\sigma}^2&:=\frac{1}{q}\sum_{s=1}^{q}\bigg(\frac{2\sqrt{m_n}}{n}\sum_{k=1}^nw(Y_k)\big(\hat{h}^{(s)}(Y_k)-\hat{h}(Y_k)\big)(\hat{h}(Y_k)\hat{c}_1+\hat{c}_2-\Lambda_{\hat{\t}}(Y_k))
\\[0,2cm]&\quad\quad+\frac{1}{\sqrt{m_n}}\sum_{j=(s-1)m_n+1}^{sm_n}\bigg(w(Y_j)(\hat{h}(Y_j)\hat{c}_1+\hat{c}_2-\Lambda_{\hat{\t}}(Y_j))^2
\\[0,2cm]&\quad\quad-\frac{1}{n}\sum_{i=1}^nw(Y_i)(\hat{h}(Y_i)\hat{c}_1+\hat{c}_2-\Lambda_{\hat{\t}}(Y_i))^2\bigg)\bigg)^2
\end{align*}
as an estimator for $\sigma^2$, where $\hat{h}^{(s)}$ denotes the nonparametric estimator for $h$ depending on the subsample $(Y_{(s-1)m+1},X_{(s-1)m+1})...,(Y_{sm},X_{sm})$, but suitable choices for $m_n$ are still unclear.

%

\section{A bootstrap version and simulations}  \label{simulations}
Although Theorem \ref{theolocalt} shows how the test statistic behaves asymptotically under $H_0$, it is hard to extract any information about how to choose appropriate critical values of a test that rejects $H_0$ for large values of $T_n$.  The main reasons for this are that first for any function $\zeta$ the eigenvalues of the operator defined in Theorem \ref{theolocalt} are unknown, that second this function is unknown and has to be estimated as well, and that third even $\psi$ (which would be needed to estimate $\zeta$) mostly is unknown and rather complex (see e.g.\ Appendix \ref{estimation-of-h}). Therefore, approximating the $\alpha$-quantile, say $q_\alpha$, of the distribution of $T$ in Theorem \ref{theolocalt} in a direct way is difficult and instead we suggest a smooth bootstrap algorithm to approximate $q_\alpha$.

\begin{alg}\label{bootstrapalg}\rm
	Let $(Y_1,X_1),...,(Y_n,X_n)$ denote the observed data, define
	$$h_{\t}(y)=\frac{\Lambda_{\t}(y)-\Lambda_{\t}(0)}{\Lambda_{\t}(1)-\Lambda_{\t}(0)}\quad\textup{and}\quad g_{\t}(x)=E[h_{\t}(Y)|X=x]$$
	and let $\hat{g}$ be a consistent estimator of $g_{\t_0}$, where $\t_0$ is defined as in \ref{A6} under the null hypothesis and as in \ref{A6'} under the alternative (see Appendix \ref{assumptions} for the assumptions). Let $\kappa $ and $\ell $ be smooth Lebesgue densities on $\mathbb{R}^{d_X}$ and $\mathbb{R}$, respectively, where $\ell $ is strictly positive, $\kappa $ has bounded support and $\kappa (0)>0$. Let $(a_n)_n$ and $(b_n)_n$ be positive sequences with $a_n\to 0$, $b_n\to 0$, $na_n\to\infty$, $nb_n^{d_X}\to\infty$. Denote by $m\in\mathbb{N}$ the sample size of the bootstrap sample. 
	\begin{enumerate}[label=(\textbf{\arabic{*}})]
		\item Calculate $\hat{\gamma}=(\hat{c}_1,\hat{c}_2,\hat{\t})^t=\arg\underset{\gamma\in\Upsilon}{\min}\,\sum_{i=1}^nw(Y_i)(\hat{h}(Y_i)c_1+c_2-\Lambda_{\t}(Y_i))^2$. Estimate the parametric residuals $\varepsilon_i(\theta_0)=h_{\t_0}(Y_i)-g_{\t_0}(X_i)$ by
		$\hat{\varepsilon}_i=h_{\hat\t}(Y_i)-\hat{g}(X_i)$ 	and denote centered versions by $\tilde \varepsilon_i=\hat\varepsilon_i-n^{-1}\sum_{j=1}^n\hat\varepsilon_j$, $i=1,\dots,n$. 
		\item 
		Generate $X_j^*$, $j=1,\dots,m$, independently (given the original data) from the density
		$$f_{X^*}(x)=\frac{1}{nb_n^{d_X}}\sum_{i=1}^n\kappa \bigg(\frac{x-X_i}{b_n}\bigg)$$
		(which is a kernel density estimator for $f_X$ with kernel $\kappa $  and bandwidth $b_n$). 
		For $j=1,\dots,m$ define bootstrap observations as 
		\begin{equation}\label{defY*}
Y_j^*=(h^*)^{-1}\big(\hat{g}(X_j^*)+\varepsilon_j^*\big)\quad\textup{for}\quad h^*(\cdot)=\frac{\Lambda_{\hat{\t}}(\cdot)-\Lambda_{\hat{\t}}(0)}{\Lambda_{\hat{\t}}(1)-\Lambda_{\hat{\t}}(0)},
		\end{equation}
		where $\varepsilon_j^*$ is generated independently (given the original data) from the density 
		$$\frac 1n \sum_{i=1}^n \frac{1}{a_n}\ell \left(\frac{\tilde\varepsilon_i-\cdot}{a_n}\right)$$
		(which is a kernel density estimator for the density of $\varepsilon(\theta_0)$ with kernel $\ell $ and bandwidth $a_n$). 
		\item Calculate the bootstrap estimate $\hat{h}^*$ for $h^*$ from $(Y_j^*,X_j^*),j=1,...,m$.
		\item Calculate the bootstrap statistic $T_{n,m}^*=\underset{(c_1,c_2,\t)\in\Upsilon}{\min}\,\sum_{j=1}^mw(Y_j^*)(\hat{h}^*(Y_j^*)c_1+c_2-\Lambda_{\t}(Y_j^*))^2$.
		\item Let $B\in\mathbb{N}$. Repeat steps (\textbf{2})--(\textbf{4}) $B$ times to obtain the bootstrap statistics $T_{n,m,1}^*,\dots,T_{n,m,B}^*$. Let $q_{\alpha}^*$ denote the quantile of $T_{n,m}^*$ conditional on $(Y_i,X_i),i=1,...,n$. Estimate $q_{\alpha}^*$ by
		$$\hat{q}_{\alpha}^*=\min\,\bigg\{z\in\{T_{n,m,1}^*,...,T_{n,m,B}^*\}:\frac{1}{B}\sum_{k=1}^BI_{\{T_{n,m,k}^*\leq z\}}\geq\alpha\bigg\}.$$
	\end{enumerate}
\end{alg}

\begin{rem}\label{rembootsassumpCKC}
	\begin{enumerate}
		\item The properties $nb_n^{d_X}\rightarrow\infty$ and $\kappa (0)>0$ ensure that conditional on the original data $(Y_1,X_1),...,(Y_n,X_n)$ the support of $X^*$ contains that of $v$ (from assumption \ref{B7} in Appendix \ref{estimation-of-h}) with probability converging to one. Thus, $v$ can be used for calculating $\hat{h}^*$ as well.
		\item To proceed as in Algorithm \ref{bootstrapalg} it may be necessary to modify $h^*$ so that $S_j^*=\hat{g}(X_j^*)+\varepsilon_{j}^*$ belongs to the domain of $(h^*)^{-1}$ for all $j=1,...,m$. As long as these modifications do not have any influence on $h^*(y)$ for $y\in\mathcal{Y}_w$, the influence on the $\hat{h}^*$ and $T_{n,m}$ should be asymptotically negligible (which can be proven for the estimator by \citet{CvK2018}).
	\end{enumerate}
\end{rem}

The bootstrap algorithm should fulfil two properties: On the one hand, under the null hypothesis the algorithm has to provide, conditionally on the original data, consistent estimates of the quantiles of $T_n$, or to be precise its asymptotic distribution from Theorem \ref{theolocalt}. To formalize this, let $(\Omega,\mathcal{A},P)$ denote the underlying probability space. Assume that $(\Omega,\mathcal{A})$ can be written as $\Omega=\Omega_1\times\Omega_2$ and $\mathcal{A}=\mathcal{A}_1\otimes\mathcal{A}_2$ for some measurable spaces $(\Omega_1,\mathcal{A}_1)$ and $(\Omega_2,\mathcal{A}_2)$. Further, assume that $P$ is characterized as the product of a probability measure $P_1$ on $(\Omega_1,\mathcal{A}_1)$ and a Markov kernel
$$P_2^1:\Omega_1\times\mathcal{A}_2\rightarrow[0,1],$$
that is $P=P_1\otimes P_2^1$. 
While randomness with respect to the original data is modelled by $P_1$, randomness with respect to the bootstrap data and conditional on the original data is modelled by $P_2^1$. Moreover, assume
$$P_2^1(\omega,A)=P\big(\Omega_1\times A|(Y_1(\omega),X_1(\omega)),...,(Y_n(\omega),X_n(\omega))\big)\quad\textup{for all }\omega\in\Omega_1,A\in\mathcal{A}_2.$$
With these notations in mind 
for all $q\in(0,\infty)$ it would be desirable to obtain
\begin{equation}\label{bootstrapH0}
P_1\Big(\omega\in\Omega_1:\underset{m\rightarrow\infty}{\limsup}\,\big|P_2^1(\omega,\{T_{n,m}^*\leq q\})-P(T_n\leq q)\big|>\delta\Big)=o(1)
\end{equation}
for all $\delta>0$ and $n\rightarrow\infty$. Here, the convention
$$P_2^1(\omega,\{T_{n,m}^*\leq q\})=P_2^1\big(\omega,\big\{\tilde{\omega}\in\Omega_2:(\omega,\tilde{\omega})\in\{T_{n,m}^*\leq q\}\big\}\big)$$
is used. On the other hand, to be consistent under $H_1$ the bootstrap quantiles have to stabilize or at least converge to infinity with a rate less than that of $T_n$. To be precise, it is needed that
\begin{equation}\label{bootstrapH1}
P_1\Big(\omega\in\Omega_1:\underset{m\rightarrow\infty}{\limsup}\,P_2^1(\omega,\{T_n<T_{n,m}^*\})>\delta\Big)=o(1)
\end{equation}
for all $\delta>0$. 

In the supplement 
we give conditions under which the bootstrap Algorithm \ref{bootstrapalg} has the desired properties (\ref{bootstrapH0}) and (\ref{bootstrapH1}). In particular we need an expansion of $\hat h^*$ as bootstrap counterpart to (\ref{exprh}).
To formulate this, for any realisation $w\in\Omega_1$ define
$$F_{Y^*}(y)=P_2^1(\omega,\{Y^*_1\leq y\}),\;\mathcal{T}^*(y)=\frac{F_{Y^*}(y)-F_{Y^*}(0)}{F_{Y^*}(1)-F_{Y^*}(0)}\textup{ and } S^*=h^*(Y^*).$$
Then for any compact set $\mathcal{K}\subseteq\mathbb{R}$ and
$$A_{m,n,\delta}= \bigg\{\underset{y\in\mathcal{K}}{\sup}\,\bigg|\hat{h}^*(y)-h^*(y)-\frac{1}{m}\sum_{j=1}^m\psi^*(S_j^*,X_j^*,\mathcal{T}^*(y))\bigg|>\frac{\delta}{\sqrt{m}}\bigg\}$$
 we need 
\begin{equation}
P_1\left(\omega\in\Omega_1:\forall\delta>0:\underset{m\rightarrow\infty}{\limsup}\,P_2^1\left(\omega,A_{m,n,\delta}\right)=0\right)=1+o(1)\label{exprh*}
\end{equation}
for $n\rightarrow\infty$, where $\psi^*$ fulfils some assumptions given in the supplement 
(see assumption (\textbf{A8*}) for details). 
In the supplement 
we also give conditions under which for the transformation estimator of \citet{CvK2018} the expansion is valid (see Lemma \ref{lemmabootstrapalternative}). 

\medskip

\subsection*{Simulations}

Throughout this section,
$g(X)=4X-1$, $X\sim\mathcal{U}([0,1])$ and $\varepsilon\sim\mathcal{N}(0,1)$
are chosen. Moreover, the null hypothesis of $h$ belonging to the \cite{YJ2000} transformations
$$\Lambda_{\theta}(Y)=\left\{\begin{array}{rc}\frac{(Y+1)^{\theta}-1}{\theta},&\textup{if }Y\geq0,\theta\neq0\\
\log(Y+1),&\textup{if }Y\geq0,\theta=0\\
-\frac{(1-Y)^{2-\theta}-1}{2-\theta},&\textup{if }Y<0,\theta\neq2\\
-\log(1-Y),&\textup{if }Y<0,\theta= 2.\end{array}\right.$$
with parameter $\theta\in\Theta_0=[0,2]$ is tested. 
Under $H_0$ we generate data using the transformation $h=(\Lambda_{\theta_0} (\cdot)-\Lambda_{\theta_0}(0))/(\Lambda_{\theta_0} (1)-\Lambda_{\theta_0}(0))$ to match the identification constraints $h(0)=0,  h(1)=1$.
Under the alternative we choose transformations $h$ with an inverse given by the following convex combination,
\begin{equation}\label{simalt}
h^{-1}(Y)=\frac{(1-c)(\Lambda_{\t_0}^{-1}(Y)-\Lambda_{\t_0}^{-1}(0))+c(r(Y)-r(0))}{(1-c)(\Lambda_{\t_0}^{-1}(1)-\Lambda_{\t_0}^{-1}(0))+c(r(1)-r(0))}
\end{equation}
for some $\theta_0\in[0,2]$, some strictly increasing function $r$ and some $c\in[0,1]$. 
In general it is not clear if a growing factor $c$ leads to a growing distance (\ref{defd}). Indeed, the opposite might be the case, if $r$ is somehow close to the class of transformation functions considered in the null hypothesis.
Simulations were conducted for $r_1(Y)=5\Phi(Y)$, $r_2(Y)=\exp(Y)$ and $r_3(Y)=Y^3$, 
where $\Phi$ denotes the cumulative distribution function of a standard normal distribution, and $c=0,0.2,0.4,0.6,0.8,1$.
The prefactor in the definition of $r_1$ is introduced because the values of $\Phi$  are rather small compared to the values of $\Lambda_{\t}$, that is, even when using the presented convex combination in (\ref{simalt}), $\Lambda_{\t_0}$ (except for $c=1$) would dominate the ``alternative part'' $r$ of the transformation function without this factor. 
 Note that $r_2$ and $\Lambda_{0}$ only differ with respect to a different standardization. Therefore, if $h$ is defined via (\ref{simalt}) with $r=r_2$ the resulting function is for $c=1$ close to the null hypothesis case. 

For calculating the test statistic the weighting function $w$  was set equal to one. The nonparametric estimator of $h$ was calculated as in \cite{CvK2018} (see Appendix  \ref{estimation-of-h} for details) with  the Epanechnikov kernel
$K(y)=\frac{3}{4}(1-y^2)^2I_{[-1,1]}(y)$ and a normal reference rule bandwidth (see for example \cite{Sil1986})
$$h_u
=\bigg(\frac{40\sqrt{\pi}}{n}\bigg)^{\frac{1}{5}}\hat{\sigma}_u,\quad h_x=\bigg(\frac{40\sqrt{\pi}}{n}\bigg)^{\frac{1}{5}}\hat{\sigma}_x,$$
where $\hat{\sigma}_u^2$ and $\hat{\sigma}_x^2$ are estimators for the variance of $U=\mathcal{T}(Y)$ and $X$, respectively. The number of evaluation points $N_x$ for the nonparametric estimator of $h$ was set equal to $100$ (see Appendix \ref{estimation-of-h} for details). The integral in (\ref{defs}) was computed by applying the function \textit{integrate} implemented in \textit{R}.
 In each  simulation run $n=100$ independent and identically distributed random pairs $(Y_1,X_1),...,(Y_n,X_n)$ were generated as described before and $250$ bootstrap quantiles, which are based on $m=100$ bootstrap observations $(Y_1^*,X_1^*),...,(Y_m^*,X_m^*)$, were calculated as in Algorithm \ref{bootstrapalg} using $\kappa$ the $U([-1,1])$-density, $\ell$ the standard normal density and $ a_n=b_n=0.1$. To obtain more precise estimators of the rejection probabilities under the null hypothesis, $800$ simulation runs were performed for each choice of $\theta_0$ under the null hypothesis, whereas in the remaining alternative cases $200$ runs were conducted.
Among other things the nonparametric estimation of $h$, the integration in (\ref{defs}), the optimization with respect to $\theta$ and the number of bootstrap repetitions cause the simulations to be quite computationally demanding. Hence, an interface for C++ as well as parallelization were used to conduct the simulations.

\bigskip

\begin{table}[htbp]
\centering
\begin{small}
		\begin{tabular}{|l|l|r|r|r|r|r|r|r|r|}
		\hline
		\multicolumn{2}{|c|}{}&\multicolumn{2}{c|}{$\theta_0=0$}&\multicolumn{2}{c|}{$\theta_0=0.5$}&\multicolumn{2}{c|}{$\theta_0=1$}&\multicolumn{2}{c|}{$\theta_0=2$}\\
		\hline
		\multicolumn{2}{|c|}{level $\alpha$}&$0.05$&$0.10$&$0.05$&$0.10$&$0.05$&$0.10$&$0.05$&$0.10$\\
		\hline
		\multirow{6}{*}{$r_1$}
		&null hyp.&0.01000&0.04000&0.03125&0.08750&0.03125&0.07750&0.01625&0.05625\\
		\hline
		&c=0.2&0.000&0.010&0.075&0.105&0.010&0.015&0.000&0.020\\
		&c=0.4&0.000&0.000&0.020&0.045&0.000&0.015&0.120&0.200\\
		&c=0.6&0.100&0.155&0.035&0.050&0.085&0.150&0.415&0.545\\
		&c=0.8&0.685&0.765&0.110&0.210&0.505&0.645&0.785&0.890\\
		&c=1&0.965&0.990&0.925&0.975&0.975&0.985&0.985&0.990\\
		\hline
		\multirow{6}{*}{$r_2$}
		&c=0.2&0.010&0.035&0.030&0.045&0.515&0.640&0.885&0.965\\
		&c=0.4&0.015&0.040&0.000&0.005&0.060&0.135&0.870&0.980\\
		&c=0.6&0.035&0.085&0.000&0.005&0.005&0.005&0.625&0.815\\
		&c=0.8&0.020&0.040&0.010&0.040&0.000&0.005&0.185&0.325\\
		&c=1&0.020&0.065&0.030&0.090&0.025&0.095&0.050&0.105\\
		\hline
		\multirow{6}{*}{$r_3$}
		&c=0.2&0.330&0.505&0.730&0.855&0.810&0.905&0.930&0.995\\
		&c=0.4&0.730&0.865&0.815&0.945&0.875&0.970&0.915&0.990\\
		&c=0.6&0.880&0.940&0.895&0.960&0.950&0.995&0.940&0.990\\
		&c=0.8&0.895&0.965&0.925&0.975&0.935&0.990&0.915&0.980\\
		&c=1&0.980&0.990&0.960&0.990&0.939&0.990&0.940&0.985\\
		\hline
	\end{tabular}
\end{small}
	\caption[Rejection probabilities under the Models (\ref{altPhi})--(\ref{alty3})]{Rejection probabilities at $\theta_0\in\{0,0.5,1,2\}$ and $r\in \{r_1,r_2,r_3\}$.}
	\label{table1}
\end{table}\noindent

The main results of the simulation study are presented in Table \ref{table1}. There, the rejection probabilities of the settings with $h=(\Lambda_{\theta_0} (\cdot)-\Lambda_{\theta_0}(0))/(\Lambda_{\theta_0} (1)-\Lambda_{\theta_0}(0))$ under the null hypothesis, and $h$ as in (\ref{simalt}) under the alternative with $r\in\{r_1,r_2,r_3\}$, $c\in\{0,0.2,0.4,0.6,0.8,1\}$ and $\theta_0\in\{0,0.5,1,2\}$ are listed. The significance level was set equal to 0.05 and 0.10. 
Note that the test sticks to the level or is even a bit conservative. Under the alternatives the rejection probabilities not only differ between different choices of $r$, but also between different transformation parameters $\theta_0$ that are inserted in (\ref{simalt}). While the test shows high power for some alternatives, there are also cases, where  the rejection probabilities are  extremely small. There are certain reasons that explain these observations. First, the class of Yeo-Johnson transforms seems to be quite general and second the testing approach itself is rather flexible due to the minimization with respect to $\gamma$. Having a look at the definition of the test statistic in (\ref{teststatistic}), it attains small values if the true transformation function can be approximated by a linear transformation of $\Lambda_{\tilde{\t}}$ for some appropriate $\tilde{\t}\in[0,2]$. In the following, this issue will be explored further by analysing some graphics.
All of the figures that occur in the following have the same structure and consist of four panels. The upper left panel shows the true transformation function with inverse function (\ref{simalt}). Due to the choice of $g(X)=4X-1$ and $X\sim\mathcal{U}([0,1])$ the vertical axis reaches from $-1$ to 3, which would be the support of $h(Y)$ if the error is neglected. In the upper right panel the parametric estimator of this function is displayed. Both of these functions are then plotted against each other in the lower left panel. Finally, the function $Y\mapsto\Lambda_{\t_0}(Y(\Lambda_{\t_0}^{-1}(1)-\Lambda_{\t_0}^{-1}(0))+\Lambda_{\t_0}^{-1}(0))$, which represents the part of $h$ corresponding to the null hypothesis, is shown in the last panel.

\begin{figure}[htbp]
	\centering
	\includegraphics[width=14cm,height=9cm]{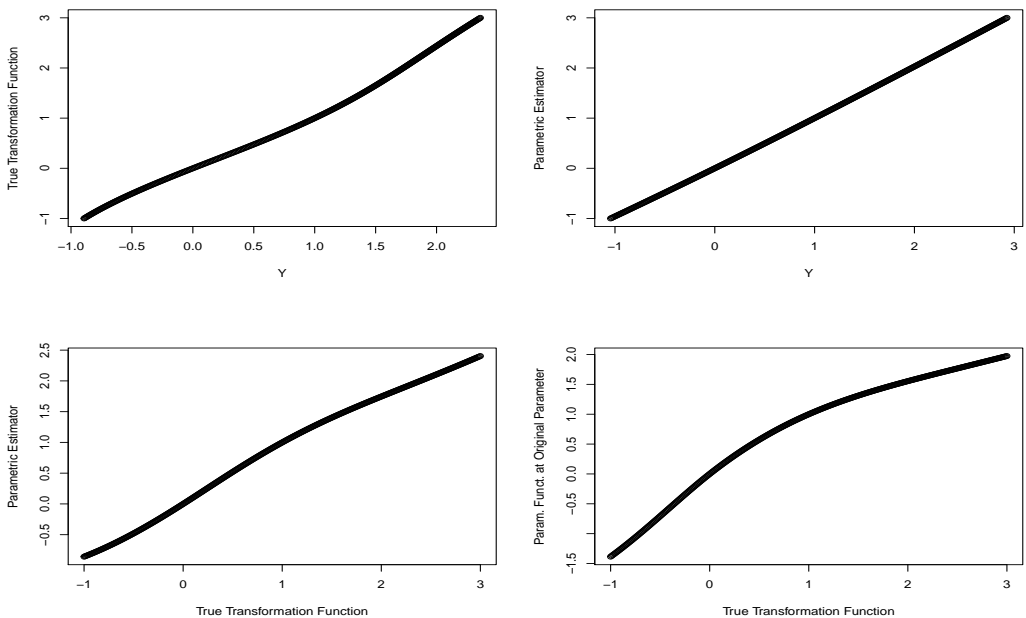}
	\caption[Transformation Functions under Model (\ref{altPhi}) for $\theta_0=0.5$]{Some transformation functions for $\theta_0=0.5,c=0.6$ and $r=r_1$.}
	\label{figureTrafotest_t05_5pnorm_c06}
\end{figure}

\begin{figure}[htbp]
	\centering
	\includegraphics[width=14cm,height=9cm]{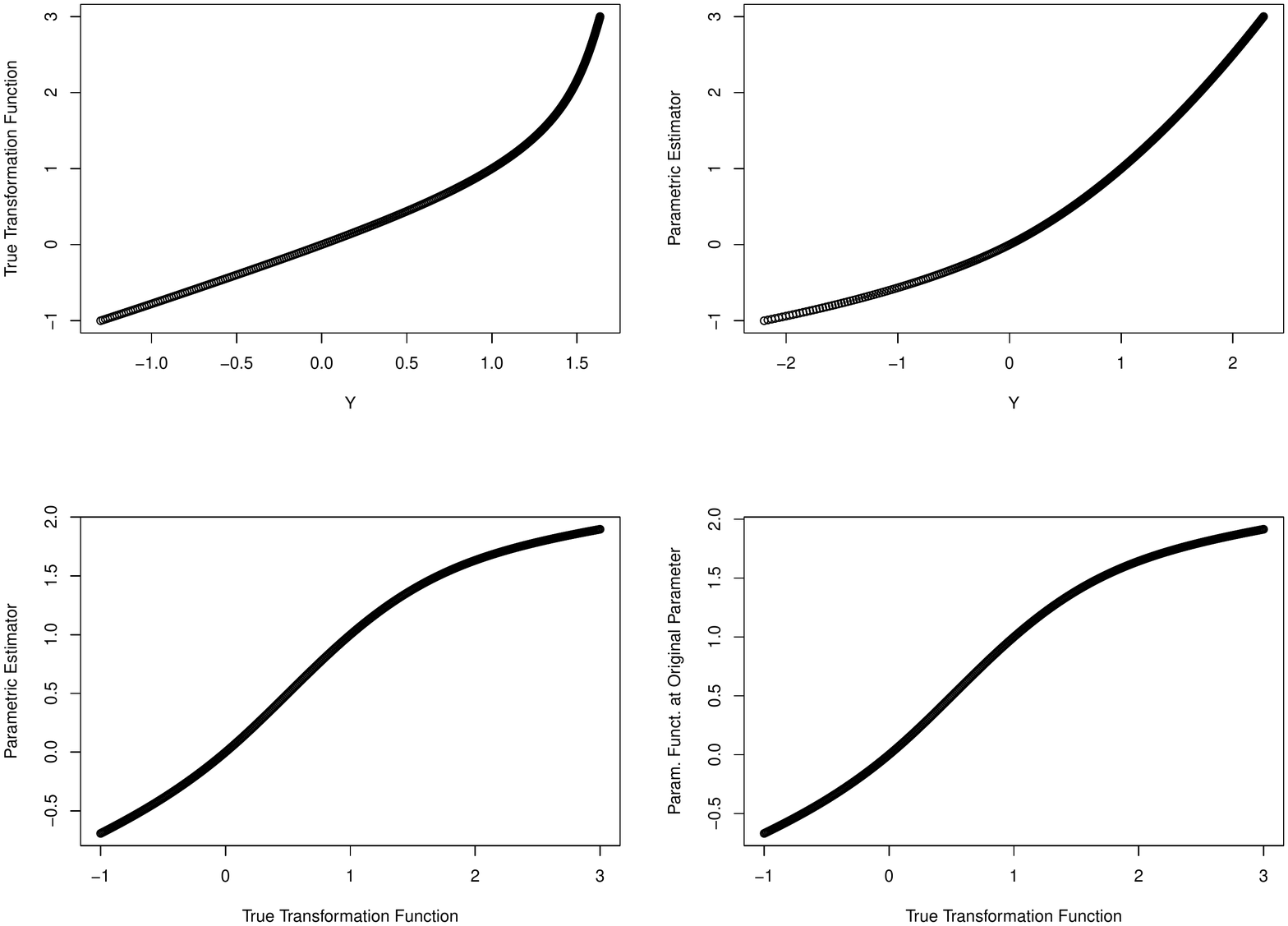}
	\caption[Transformation Functions under Model (\ref{altPhi}) for $\theta_0=2$]{Some transformation functions for $\theta_0=2,c=0.6$ and $r=r_1$.}
	\label{figureTrafotest_t2_5pnorm_c06}
\end{figure}

\begin{figure}[htbp]
	\centering
	\includegraphics[width=14cm,height=9cm]{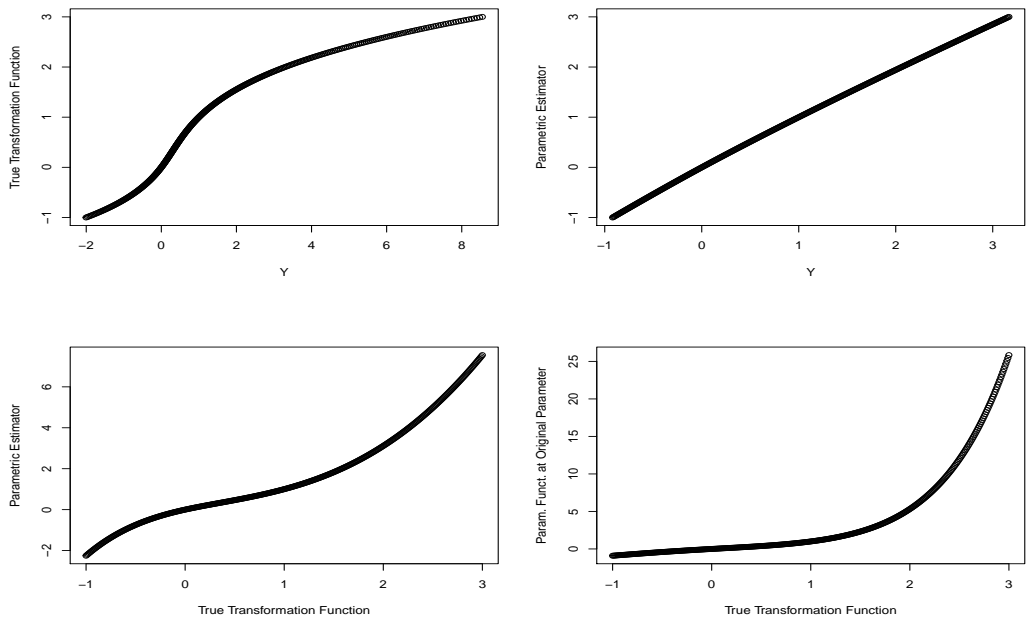}
	\caption[Transformation Functions under Model (\ref{alty3})]{Some transformation functions for $\theta_0=2,c=0.2$ and $r=r_3$.}
	\label{figureTrafotest_t2_y3_c02}
\end{figure}

In the lower left panel one can see if the true transformation function can be approximated by a linear transform of some $\Lambda_{\tilde{\t}},\tilde{\theta}\in[0,2]$, which is an indicator for rejecting or not rejecting the null hypothesis as was pointed out before.
As already mentioned, the rejection probabilities not only differ between different deviation functions $r$, but also within these settings. For example, when considering $r=r_1$ with $c=0.6$ the rejection probabilities for $\theta_0=0.5$ amount to $0.035$ for $\alpha=0.05$ and to $0.050$ for $\alpha=0.10$, while for $\theta_0=2$ they are $0.415$ and $0.545$. Figures \ref{figureTrafotest_t05_5pnorm_c06} and \ref{figureTrafotest_t2_5pnorm_c06} explain why the rejection probabilities differ that much.
While for $\theta_0=0.5$ the transformation function can be approximated quite well by transforming $\Lambda_{1.06}$ linearly, the best approximation for $\theta_0=2$ is given by $\Lambda_{1.94}$ and seems to be relatively bad. The best approximation for $c=1$ can be reached for $\theta$ around $1.4$.
In contrast to that, considering $\theta_0=2$ and $r=r_3$ results in a completely different picture. As can be seen in Figure \ref{figureTrafotest_t2_y3_c02} even for $c=0.2$ the resulting $h$ differs so much from the null hypothesis that it can not be linearly transformed into a Yeo-Johnson transform (see the lower left subgraphic). Consequently, the rejection probabilities are rather high. 

\noindent
\noindent

\noindent

\begin{figure}[htbp]
	\centering
	\includegraphics[width=14cm,height=9cm]{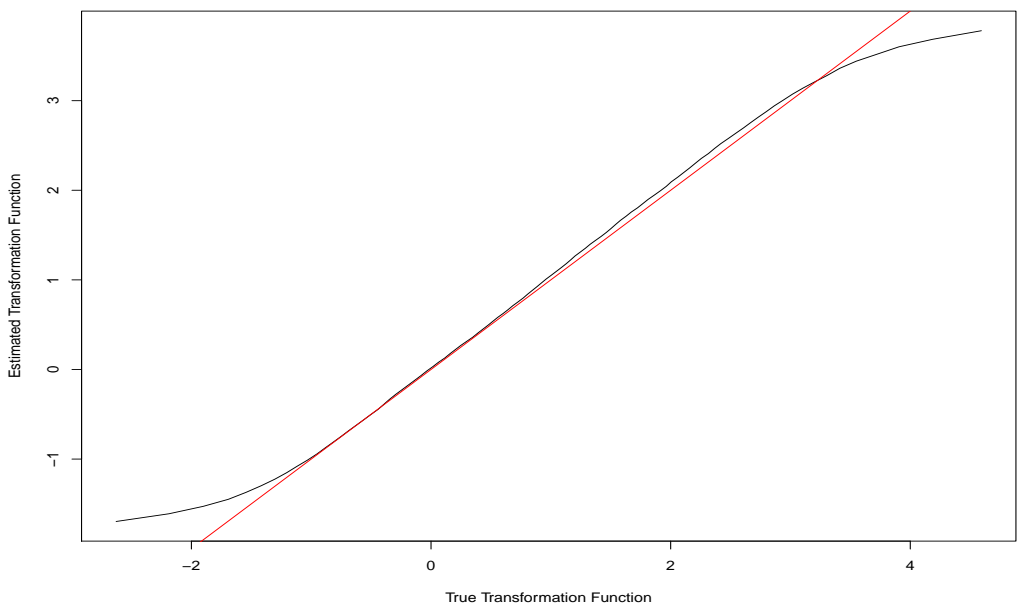}
	\caption[Bias of the Nonparametric Estimator]{The transformation function and its nonparametric estimator for $\theta_0=1,c=0.4$ and $r=r_1$.}
	\label{figureBias_nonparamest_t1c045pnorm}
\end{figure}

\noindent
Under some alternatives the rejection probabilities are even smaller than the level. This behaviour indicates that from the presented test's perspective these models seem to fulfil the null hypothesis more convincingly than the null hypothesis models themselves. The reason for this can be seen in Figure \ref{figureBias_nonparamest_t1c045pnorm} for the setting $\theta_0=1,c=0.4$ and $r=r_1$. There, the relationship between the nonparametric estimator of the transformation function and the true transformation function is shown. While the diagonal line represents the identity, the nonparametric estimator seems to flatten the edges of the transformation function. In contrast to this, using $r=r_1$ in (\ref{simalt}) steepens the edges so that both effects neutralize each other. Similar effects cause low rejection probabilities for $r=r_2$, although the reasoning is slightly more sophisticated and is also associated with the boundedness of the parameter space $\Theta_0=[0,2]$.\\
\noindent
One possible solution could consist in adjusting the weight function $w$ such that the boundary of the support of $Y$ does no longer belong to the support of $w$. In Table \ref{table4} the rejection probabilities for a modified weighting approach are presented. There, the weight function was chosen such that the smallest five percent and the largest five percent of observations were omitted to avoid the flattening effect of the nonparametric estimation. Indeed, the resulting rejection probabilities under the alternatives increase and lie above those under the null hypotheses.

\begin{table}[htbp]\footnotesize
	\centering
	\begin{footnotesize}
		\begin{tabular}{|l|l|c|c|c|c|}
			\hline
			&Alternative&\multicolumn{2}{|c|}{original framework}&\multicolumn{2}{|c|}{modified weighting}\\
			\hline
			Param.&Level&$\alpha=0.05$&$\alpha=0.10$&$\alpha=0.05$&$\alpha=0.10$\\
			\hline
			\multirow{7}{*}{$\theta_0=1$}&null hyp.&0.03125&0.07750&0.02875&0.07875\\
			&c=0.2&0.010&0.015&0.040&0.100\\
			&c=0.4&0.000&0.015&0.205&0.320\\
			&c=0.6&0.085&0.150&0.590&0.715\\
			&c=0.8&0.505&0.645&0.950&0.980\\
			&c=1&0.975&0.985&1.000&1.000\\
			\hline
			\multirow{7}{*}{$\theta_0=2$}&null hyp.&0.01625&0.05625&0.05500&0.10375\\
			&c=0.2&0.000&0.020&0.225&0.350\\
			&c=0.4&0.120&0.200&0.575&0.710\\
			&c=0.6&0.415&0.545&0.910&0.965\\
			&c=0.8&0.785&0.890&0.990&1.000\\
			&c=1&0.985&0.990&0.995&1.000\\
			\hline
		\end{tabular}
	\end{footnotesize}
	\caption{Rejection probabilities at $\theta_0=1$ and $\theta_0=2$ for $r=r_1$.}
	\label{table4}
\end{table}

\bigskip

{\bf Acknowledgements.} Natalie Neumeyer acknowledges financial support by the DFG (Research Unit FOR 1735 {\it Structural Inference in
Statistics: Adaptation and Efficiency}). Ingrid Van Keilegom acknowledges financial support by the European Research Council (2016- 2021, Horizon 2020 / ERC grant agreement No.\ 694409).

\bigskip

\bibliographystyle{plainnat}
\bibliography{References}

\appendix

\section{Assumptions for the main results}\label{assumptions}

In the following assumptions let $\cy$ denote the support of $Y$ (which depends on $n$ under local alternatives).
Further, $F_S$ denotes the distribution function of $S_1$ as in (\ref{S_i}) and $\mathcal{T}_S$ denotes the transformation $s\mapsto (F_S(s)-F_S(0))/(F_S(1)-F_S(0))$.

\begin{enumerate}[label=(\textbf{A\arabic{*}})]
	\item\label{A1} The sets $C_1,C_2$ and $\T$ are compact.
	\item\label{A2} The weight function $w$ is continuous with a compact support $\cyw\subset \cy$. 
	\item\label{A3} The map $(y,\t)\mapsto\Lambda_{\t}(y)$ is twice continuously differentiable on $\cyw$ with respect to $\t$ and the (partial) derivatives are continuous in $(y,\t)\in\cyw\times\T$.
	\item\label{A4} There exists a unique strictly increasing and continuous transformation $h$ such that model (\ref{modeleq}) holds with $X$ independent of $\eps$.
	\item\label{A5} The function $h_0$ defined in (\ref{h0r0}) is strictly increasing and continuously differentiable and $r$ is continuous on $\cyw$. $F_Y$ is strictly increasing on the support of $Y$.
	\item\label{A6} Minimizing the function $M:\Upsilon\rightarrow\R,\gamma=(c_1,c_2,\theta)\mapsto E\big[w(Y)(h_0(Y)c_1+c_2-\Lt(Y))^2\big]$ leads to a unique solution $\gamma_0=(c_{1,0},c_{2,0},\theta_0)$ in the interior of $\Upsilon$. For all $\t\neq\tilde{\t}$ it is $\underset{y\in\operatorname{supp}(w)}{\sup}\,\big|\frac{\Lambda_{\t}(y)-\Lambda_{\t}(0)}{\Lambda_{\t}(1)-\Lambda_{\t}(0)}-\frac{\Lambda_{\tilde{\t}}(y)-\Lambda_{\tilde{\t}}(0)}{\Lambda_{\tilde{\t}}(1)-\Lambda_{\tilde{\t}}(0)}\big|>0$.
	\item\label{A7} The Hessian matrix $\Gamma:=\operatorname{Hess}M(\gamma_0)$ is positive definite.
	\item\label{A8} The transformation estimator $\hat{h}$ fulfills (\ref{exprh}) for some function $\psi$. For some $\mathcal{U}_0$ (independent of $n$ under local alternatives) with $\mathcal{T}_S(h(\cyw))\subset \mathcal{U}_0$ the function class $\{z\mapsto\psi(z,t):t\in\mathcal{U}_0\}$ is Donsker with respect to $P^Z$ and $E[\psi(Z_1,t)]=0$ for all $t\in\mathcal{U}_0$. The fourth moment $E[w(h_0^{-1}(S_1))\psi(Z_1,U_1)^4]$ is finite and the conditional moments $E[w(h_0^{-1}(S_1))\psi(Z_1,U_2)^2|Z_1=z]$ are locally bounded.
	\end{enumerate}
	


\noindent When considering a fixed alternative $H_1$ or the relevant hypothesis $H'_0$, \ref{A6} and \ref{A8} are replaced by the following Assumptions \ref{A6'} and \ref{A8'} (assumption \ref{A8'} is only relevant for $H'_0$). Note that $h$ is a fixed function then, not depending on $n$. 

\begin{enumerate}[label=(\textbf{A6'})]
	\item\label{A6'} Minimizing the function $M:\Upsilon\rightarrow\R,\gamma=(c_1,c_2,\theta)\mapsto E\big[w(Y)(h(Y)c_1+c_2-\Lt(Y))^2\big]$
	leads to a unique solution $\gamma_0=(c_{1,0},c_{2,0},\theta_0)$ in the interior of $\Upsilon$. For all $\t\neq\tilde{\t}$ it is $\underset{y\in\operatorname{supp}(w)}{\sup}\,\big|\frac{\Lambda_{\t}(y)-\Lambda_{\t}(0)}{\Lambda_{\t}(1)-\Lambda_{\t}(0)}-\frac{\Lambda_{\tilde{\t}}(y)-\Lambda_{\tilde{\t}}(0)}{\Lambda_{\tilde{\t}}(1)-\Lambda_{\tilde{\t}}(0)}\big|>0$.
\end{enumerate}
\begin{enumerate}[label=(\textbf{A8'})]
	\item\label{A8'} The transformation estimator $\hat{h}$ fulfills (\ref{exprh}) for some function $\psi$. For some $\mathcal{U}_0\supset \mathcal{T}_S(h(\cyw))$ the function class $\{z\mapsto\psi(z,t):t\in\mathcal{U}_0\}$ is Donsker with respect to $P^Z$ and $E[\psi(Z_1,t)]=0$ for all $t\in\mathcal{U}_0$. Further, one has $E[\psi(Z_1,U_2)^2]<\infty$.
\end{enumerate}


\section{Nonparametric transformation estimation}\label{estimation-of-h} 

In this section we consider a transformation estimator which fulfills assumption \ref{A8} and in particular the expansion (\ref{exprh}). To this end we reproduce the definitions of \citet{CvK2018} and prove Lemma \ref{expansion} below. Denote the conditional distribution function of $U_1$ from (\ref{U_i}), given $X_1=x$,  by $F_{U|X}(\cdot|x)$ and estimate it by 
 $$\hat{F}_{U|X}(u|x):=\frac{\sum_{i=1}^nK_{h_x}(X_i-x)\mathcal{K}_{h_u}(u-\hat{U}_i)}{\sum_{i=1}^nK_{h_x}(X_i-x)}.$$
Here $h_x$ and $h_u$ are bandwidths and $K$ is an appropriate kernel function (as in assumptions \ref{B4} and \ref{B5} below),
$$K_h(u):=\frac{1}{h}K\bigg(\frac{u}{h}\bigg),\quad\textup{and}\quad\mathcal{K}_h(u):=\int_{-\infty}^uK_h(r)\,dr=\int_{-\infty}^{\frac{u}{h}}K(r)\,dr.$$
Further, consider some kernel $L$ and bandwidth $b\searrow0$ fulfilling assumption \ref{B6} and define
\begin{equation}\label{hatQ}
\hat{Q}(u)=\arg\underset{q\in\mathbb{R}}{\min}\,\int v(x)\left(\frac{\hat s_1(u,x)}{\hat s_1(1,x)}-q\right)\bigg(2L_b\bigg(\frac{\hat s_1(u,x)}{\hat s_1(1,x)}-q\bigg)-1\bigg)\,dx,
\end{equation}
where $L_b(\cdot)=\frac{1}{b}L\big(\frac{\cdot}{b}\big)$ and
\begin{equation}\label{defs}
\hat s_1(u,x):=\int_0^u\frac{\frac{\partial \hat F_{U|X}(r|x)}{\partial r}}{\frac{\partial \hat F_{U|X}(r|x)}{\partial x_1}}\,dr.
\end{equation}
\begin{rem}
In principle, the derivative with respect to any other component $x_i$ fulfilling assumption \ref{B3} below can be used in \ref{B2} as well (similar to \cite{CKC2015}). W.l.o.g.\ only the case $i=1$ is considered here.
\end{rem}
Let $\hat{F}_Y$ be the empirical distribution function of $Y_1,\dots,Y_n$ and define
\begin{equation}\label{hatT}
\hat{\mathcal{T}}(y):=\frac{\hat{F}_Y(y)-\hat{F}_Y(0)}{\hat{F}_Y(1)-\hat{F}_Y(0)}
\end{equation}
and estimate $U_i$ from (\ref{U_i}) by $\hat{U}_i:=\hat{\mathcal{T}}(Y_i)$. The estimator for $h$ can be defined as
\begin{equation}\label{hath}
\hat{h}(y)=\hat{Q}(\hat{\mathcal{T}}(y)).
\end{equation}
According to \citet{CvK2018} (see their proof of Propositions 6 and 7)
one has expansion (\ref{exprh}) for $\hat{h}$  under the null hypothesis of the parametric transformation class. For the definition of $\psi(Z_i,u)$ from (\ref{exprh}) in
(\ref{def-psi})  below we need the following notations. 
Let $f_{U,X}$ denote the joint density of $(U,X)$ and define
\begin{equation}\label{defpf}
p(u,x)=\int_{-\infty}^uf_{U,X}(r,x)\,dr\quad\textup{as well as}\quad f_X(x)=\int_{\mathbb{R}}f_{U,X}(r,x)\,dr.
\end{equation}
Then, the conditional distribution function of $U$ conditioned on $X=x$ can be written as
$$\Phi(u,x)=\frac{p(u,x)}{f_X(x)}.$$
Further, let
\begin{eqnarray*}
&& \Phi_u(u,x)=\frac{\partial}{\partial u}\Phi(u,x),\quad \Phi_1(u,x)=\frac{\partial}{\partial x_1}\Phi(u,x), \quad f_{X,1}(x)=\frac{\partial}{\partial x_1}f_X(x),
\\
&&D_{p,0}(u,x)=\frac{\Phi_u(x,u)f_{X,1}(x)}{\Phi_1^2(u,x)f_X^2(x)},\quad D_{p,u}(u,x)=\frac{1}{f_X(x)\Phi_1(u,x)},
\\
&&D_{p,1}(u,x)=-\frac{\Phi_u(u,x)}{f_X(x)\Phi_1^2(u,x)},\quad D_{f,0}(u,x)=-\frac{\Phi_u(u,x)\Phi(u,x)f_{X,1}(x)}{\Phi_1^2(u,x)f_X^2(x)},
\\
&&D_{f,1}(u,x)=\frac{\Phi_u(u,x)\Phi(u,x)}{\Phi_1^2(u,x)f_X(x)}.
\end{eqnarray*}
Define $\tilde{v}_1(u_0,x)=\frac{v(x)}{s_1(u_0,x)}$, $\tilde{v}_2(u_0,x)=\frac{v(x)s_1(u_0,x)}{s_1(1,x)^2}$
and (for $\tilde{v}\in\{\tilde{v}_1,\tilde{v}_2\}$)
\begin{align*}
\delta_j^{\tilde{v}}(u_0,u)&=\int_{\max(0,U_j)}^{\max(u,U_j)}\bigg(\tilde{v}(u_0,X_j)D_{p,0}(r,X_j)-\frac{\partial}{\partial x_1}\big(\tilde{v}(u_0,x)D_{p,1}(r,x)\big)\Big|_{x=X_j}\bigg)\,dr
\\[0,2cm]&\quad+\int_0^u\bigg(\tilde{v}(u_0,X_j)D_{f,0}(r,X_j)-\frac{\partial}{\partial x_1}\big(\tilde{v}(u_0,X_j)D_{f,1}(r,x)\big)\Big|_{x=X_j}\bigg)\,dr
\\[0,2cm]&\quad+(\mathds{1}_{\{U_j\leq u\}}-\mathds{1}_{\{U_j\leq0\}})\tilde{v}(u_0,X_j)D_{p,u}(U_j,X_j)
\\[0,2cm]&\quad+\int_0^u\bigg(\frac{\mathds{1}_{\{U_j\leq u\}}-\mathds{1}_{\{U_j\leq0\}}}{F_U(1)-F_U(0)}-r\bigg)
\\[0,2cm]&\quad\quad\int_{\mathcal{X}}\bigg(\bigg(\tilde{v}(u_0,x)D_{p,0}(r,x)+\frac{\partial}{\partial x_1}\big(\tilde{v}(u_0,x)D_{p,1}(r,x)\big)\bigg)
\\[0,2cm]&\quad\quad f_{U,X}(r,x)+\tilde{v}(u_0,x)D_{p,u}(r,x)\frac{\partial}{\partial r}f_{U,X}(r,x)\bigg)\,dx\,dr
\\[0,2cm]&\quad-\bigg(\frac{\mathds{1}_{\{U_j\leq 1\}}-\mathds{1}_{\{U_j\leq0\}}}{F_U(1)-F_U(0)}-1\bigg)\int_0^ur\int_{\mathcal{X}}\bigg(\tilde{v}(u_0,x)D_{p,0}(r,x)
\\[0,2cm]&\quad\quad-\tilde{v}(u_0,x)\frac{\partial}{\partial r}D_{p,u}(r,x)+\frac{\partial}{\partial x_1}\big(\tilde{v}(u_0,x)D_{p,1}(r,x)\big)\bigg)f_{U,X}(r,x)\,dx\,dr
\\[0,2cm]&\quad\quad-\bigg(\frac{\mathds{1}_{\{U_j\leq 1\}}-\mathds{1}_{\{U_j\leq0\}}}{F_U(1)-F_U(0)}-1\bigg)u\int_{\mathcal{X}}\tilde{v}(u_0,x)D_{p,u}(u,x)f_{U,X}(u,x)\,dx,
\end{align*}
see \citet{CvK2018} for details. Then, with $Q(\cdot)=h(\mathcal{T}^{-1}(\cdot))$ the function $\psi$ in the expansion (\ref{exprh}) can be written as
\begin{eqnarray}\nonumber
\psi(Z_j,u)&=&\delta_j^{\tilde{v}_1}(1,u)-\delta_j^{\tilde{v}_2}(u,1)+\frac{Q'(u)}{F_U(1)-F_U(0)}\big(\mathds{1}_{\{U_j\leq u\}}-\mathds{1}_{\{U_j\leq 0\}}-F_U(u)+F_U(0)\big)
\\&&{}-Q'(u)\frac{F_U(u)-F_U(0)}{(F_U(1)-F_U(0))^2}\big(\mathds{1}_{\{U_j\leq1\}}-\mathds{1}_{\{U_j\leq 0\}}-F_U(1)+F_U(0)\big).
\label{def-psi}
\end{eqnarray}
Note that $E[\psi(Z_j,u)]=0$ for all $u$.


In the following assumptions (adjusted from \cite{CvK2018}) are given which ensure \ref{A8} for the estimator $\hat h$ from (\ref{hath}). Let, as in \ref{A8}, $ \mathcal{T}_S(h(\cyw))\subset \mathcal{U}_0$, where $\mathcal{U}_0$ is independent from $n$ under local alternatives and lies in the interior of the support of $U$. Let $\cx\subset\R^{d_X}$ denote the support of $X$. 
\begin{enumerate}[label=(\textbf{B\arabic{*}})]
	\item\label{B1} The cumulative distribution function of $\eps$ is absolutely continuous and has a density that is continuous on its support. Furthermore, $X$ and $\eps$ are independent and $\mathcal{U}_0$ is a connected subset of $\R$.
	
	\item\label{B2} The transformation $Q(\cdot)=h(\mathcal{T}^{-1}(\cdot))$ is strictly increasing and continuously differentiable on  $\mathcal{U}_0$.
	\item\label{B3} The set
	$$\mathcal{X}_{\partial 1}:=\left\{x\in\cx: \frac{\partial F_{U|X}(u|x)}{\partial x_1}\neq0\quad\textup{for all }u\in \mathcal{U}_0\right\}$$
	is nonempty.
	\item\label{B4} The bandwidths $h_x$ and $h_u$ satisfy for an appropriate $q\in\N$
	$$\sqrt{n}h_x^q\rightarrow0,\sqrt{n}h_u^q\rightarrow0,\frac{\sqrt{n}h_x^{d_X+2}}{\log(n)}\rightarrow\infty,\quad\frac{\sqrt{n}h_x^{d_X}h_u^2}{\log(n)}\rightarrow\infty.$$
	\item\label{B5} The kernel $K$ is symmetric with a connected and compact support containing some neighbourhood around 0. Further, $K$ is $q$-times continuously differentiable with $K$ and $K'$ being of bounded variation. Moreover, $\int K(z)\,dz=1,\int z^lK(z)\,dz=0$ for all $l=1,...,q-1$.
	\item\label{B6} The kernel $L$ is twice continuously differentiable with uniformly bounded derivatives and with median 0, and $b=b_n>0$ is a bandwidth sequence that satisfies $nb^4\rightarrow\infty$ and $b\sqrt{n}h_x^{d_X}(\min(h_x,h_u))^2/\log(n)\rightarrow\infty$.
	\item\label{B7} $v$ is a weight function with compact support $\cx_0\subset\mathcal{X}_{\partial 1}$ with nonempty interior. Further, $\int_{\cx_0}v(x)\,dx=1$ and $v$ is $q$-times continuously differentiable and all these derivatives are uniformly bounded in the interior, i.e.,
	$$\underset{x\in\cx_0}{\sup}\,\bigg|\frac{\partial^{|\a|}}{\partial x_1^{\a_1}\cdots\partial x_{d_X}^{\a_{d_X}}}v(x)\bigg|<\infty,$$
	for all $\a_1,...,\a_{d_X}\in\{0,...,q-1\}$ with $|\a|=\sum_{i=1}^{d_X}\a_i\leq m$.
	\item\label{B8} The regression function $g$ is continuously differentiable with respect to $x_i$ on $\cx$ for $i=1,...,d_X$.
	\item\label{B9} The joint density function $f_{Y,X}(y,x)$ of $(Y,X)$ is uniformly bounded, $(q+2)$-times continuously differentiable and all these derivatives are uniformly bounded, i.e.,
	$$\underset{y:\,\mathcal{T}(y)\in \mathcal{U}_0,\,x\in\cx_0}{\sup}\,\bigg|\frac{\partial^{|\a|}}{\partial y^{\a_0}\partial x_1^{\a_1}\cdots\partial x_{d_X}^{\a_{d_X}}}f_{Y,X}(y,x)\bigg|<\infty,$$
	for all $\a_0,\a_1,...,\a_{d_X}\in\{0,...,q-1\}$ with $|\a|=\sum_{i=0}^{d_X}\a_i\leq q+2$. Further, we assume $\underset{y:\,\mathcal{T}(y)\in \mathcal{U}_0}{\inf}\,f_Y(y)>0$, where $f_Y$ is the density function of $Y$.
	\item\label{B10} Assume
	$$\underset{x\in\cx_0}{\inf}\,f_X(x)>0,\underset{(u,x)\in \mathcal{U}_0\times\cx_0}{\inf}\,\frac{\partial F_{U|X}(u|x)}{\partial x_1}>0\quad\textup{and}\quad\underset{x\in\cx_0}{\inf}\,\int_0^1\frac{\frac{\partial F_{U|X}(u|x)}{\partial u}(u,x)}{\frac{\partial F_{U|X}(u|x)}{\partial x_1}(u,x)}\,du>0.$$
\end{enumerate}


The following result holds under the null hypothesis $H_0$, under fixed alternatives $H_1$ and under local alternatives $H_{1,n}$. 

\begin{lemma}\label{expansion}
Assume \ref{A4} and \ref{B1}--\ref{B10}. Then, assumption \ref{A8} is fulfilled and especially the expansion (\ref{exprh}) holds for $\hat h$  
with $\psi$ from (\ref{def-psi}). 
\end{lemma}

\noindent
{\bf Proof.} In the case of a fixed transformation  $h$, the assertion is covered by Theorems 5.1 and 5.2 in \citet{CvK2018}. Therefore, only local alternatives need to be considered. To this end we consider the transformation class 
$$\H=\left\{h(\cdot)=\frac{\Lambda_{\t_0}(\cdot)-\Lambda_{\t_0}(0)+\a(r(\cdot)-r(0))}{\Lambda_{\t_0}(1)-\Lambda_{\t_0}(0)+\a(r(1)-r(0))}:\a\in A\right\},$$
for an appropriate set $A$. In case of local alternatives as in equation (\ref{h0r0}), consider for example $A=\{n^{-\frac{1}{2}}:n\in\mathbb{N},n\geq N\}$ for a sufficiently large $N\in\mathbb{N}$. The expansion is shown uniformly in $\a\in A$, that is uniformly in $h\in\H$, and uniformly in $y\in\cyw$.  Nevertheless, most arguments used for fixed $h$ as in \citet{CvK2018} are still valid. 
Note that in our framework $S=g(X)+\eps$ does not depend on $n$, and consider $Y=h^{-1}(S)$ for $h\in\H$.
First, note that neither the $U_i$ nor the $\hat{U}_i$ depend on $\a$, since $\mathds{1}_{\{Y_j\leq Y_i\}}=\mathds{1}_{\{h(Y_j)\leq h(Y_i)\}}$ is independent of $\a$ and $h(0)=0$ as well as $h(1)=1$ for all $h\in\H$. Hence, $Q$ from \ref{B2} and its estimator $\hat{Q}$ from (\ref{hatQ}) are independent of $\alpha\in A$. Moreover, $\hat{Q}'$ is uniformly consistent. By standard arguments it can be shown that for $\hat{\mathcal{T}}$ from (\ref{hatT}) and $\mathcal{T}$ from (\ref{U_i}) one has
$$\hat{\mathcal{T}}(y)-\mathcal{T}(y)=-\frac{F_Y(y)-F_Y(0)}{(F_Y(1)-F_Y(0))^2}(\hat{F}_Y(1)-\hat{F}_Y(0)-F_Y(1)+F_Y(0))+o_P(n^{-1/2})=O_P(n^{-1/2})$$
uniformly in $y\in\cyw$ and $\a\in A$, so that
\begin{align*}
\hat{h}(y)-h(y)&=\hat{Q}(\hat{\mathcal{T}}(y))-Q(\mathcal{T}(y)) 
\\[0,2cm]&=\hat{Q}(\mathcal{T}(y))-Q(\mathcal{T}(y))+\hat{Q}'(\mathcal{T}(y))(\hat{\mathcal{T}}(y)-\mathcal{T}(y))+o_P(n^{-1/2})
\\[0,2cm]&=\hat{Q}(\mathcal{T}(y))-Q(\mathcal{T}(y))+Q'(\mathcal{T}(y))(\hat{\mathcal{T}}(y)-\mathcal{T}(y))+o_P(n^{-1/2})
\\[0,2cm]&=\frac{1}{n}\sjn\psi(Z_j,\mathcal{T}(y))+o_P(n^{-1/2}).
\end{align*}
\hfill $\Box$

\section{Proofs of the main results} \label{proofs}

\bigskip

\noindent {\bf Proof of Theorem \ref{theolocalt}.}
For ease of presentation define 
\begin{equation}\label{Mn}
M_n(c_1,c_2,\t)=\sjn w(Y_j)(\hat{h}(Y_j)c_1+c_2-\Lt(Y_j))^2
\end{equation}
such that $T_n=\min_{\gamma\in\Upsilon}M_n(\gamma)$. Let $\hat{\gamma}=(\hat{c_1},\hat{c_2},\tilde\theta)$ denote the minimizer of $M_n$ and $\gamma_0$ be the vector such that
\begin{equation}\label{gamma_0}
h_0(y)c_{1,0}+c_{2,0}=\Lambda_{\t_0}(y).
\end{equation}
(see assumption \ref{A6}). Note that $c_{1,0}=\Lambda_{\t_0}(1)-\Lambda_{\t_0}(0)$, $c_{2,0}=\Lambda_{\t_0}(0)$. 

We have $\tilde \gamma-\gamma_0=o_P(1)$ because
for all $\gamma=(c_1,c_2,\t)$ in a compact set that does not contain $\gamma_0$ and an appropriate $\varepsilon>0$ one has
\begin{align*}
\frac{M_n(\gamma)}{n}&=E\big[w(Y)(h(Y)c_1+c_{2}-\Lambda_{\t}(Y))^2\big]+o_P(1)\geq \varepsilon +o_P(1)
\end{align*}
uniformly in $\gamma\in \Upsilon$ (see the proof of Theorem \ref{theoconsistency} for details).

Let in the following $\nabla f$ denote the gradient of a function $f$ and $\operatorname{Hess}f$ denote the Hessian matrix of $f$. Note that $\nabla M_n(\tilde \gamma)=0$ and thus by  Taylor expansion
\begin{eqnarray}\label{Taylor1}
M_n(\gamma_0) &=& M_n(\hat{\gamma})+(\gamma_0-\hat{\gamma})^t\operatorname{Hess}M_n(\gamma^*)(\gamma_0-\hat{\gamma}),
\end{eqnarray}
where $\gamma^*=(c_1^*,c_2^*,\t^*)$ is on the line between $\gamma_0$ and $\hat{\gamma}$. Further, with 
 $\hat\Gamma_{3,3}=\dot{\Lambda}_{{\t^*}}(Y_k)^t\dot{\Lambda}_{{\t^*}}(Y_k)-\ddot{\Lambda}_{{\t^*}}(Y_k)(\hat{h}(Y_k){c}_1^*+{c}_2^*-\Lambda_{{\t^*}}(Y_k))$
we have 
\begin{eqnarray}\nonumber
\frac{1}{n}\operatorname{Hess}M_n(\gamma^*)
&=&\frac{1}{n}\skn w(Y_k)
\left(\begin{array}{ccc}\hat{h}(Y_k)^2&\hat{h}(Y_k)&-\hat{h}(Y_k)\dot{\Lambda}_{{\t^*}}(Y_k)\\\hat{h}(Y_k)&1&-\dot{\Lambda}_{\t^*}(Y_k)\\-\hat{h}(Y_k)\dot{\Lambda}_{{\t^*}}(Y_k)^t&-\dot{\Lambda}_{{\t^*}}(Y_k)^t&\hat\Gamma_{3,3}\end{array}\right)
\nonumber\\
&=& \Gamma_0+o_P(1).\label{Gamma-limit} 
\end{eqnarray}
To obtain the last equality note that $\hat h$ converges to $h_0$ uniformly on $\Y_w$ thanks to (\ref{exprh}) under local alternatives. Further, $h^{-1}$ converges to $h_0^{-1}$ uniformly on compacta and, as $\tilde \gamma$ converges to $\gamma_0$, $\dot{\Lambda}_{\tilde{\t}}$ and $\ddot{\Lambda}_{\tilde{\t}}$ converge to $\dot{\Lambda}_{\t_0}$ and $\ddot{\Lambda}_{\t_0}$, respectively, uniformly on $\Y_w$. To obtain (\ref{Gamma-limit}) it remains to apply the law of large numbers and (\ref{gamma_0}).  

Since $\Gamma_0$ is positive definite by assumption and $M_n(\gamma_0)$ and $M_n(\hat{\gamma})$ are bounded in probability, one obtains from (\ref{Taylor1}) that $\hat{\gamma}-\gamma_0=O_P(n^{-1/2})$.

Now, again by Taylor expansion, for all values $\gamma$ with $\gamma-\gamma_0=O_P(n^{-1/2})$,
\begin{eqnarray*}
M_n(\gamma) &=& M_n(\gamma_0)+(\gamma_0-\gamma)^t\nabla M_n(\gamma_0)+(\gamma_0-\gamma)^t\operatorname{Hess}M_n(\gamma_0)(\gamma_0-\gamma)+o_P(1)\\
&=& q(\gamma_0-\gamma)+o_P(1),
\end{eqnarray*}
where the map $q$ is defined via
$$q(z)=M_n(\gamma_0)+nz^t\Gamma_0 z+2z^t\skn w(h_0^{-1}(S_k))(\hat{h}(Y_k)c_{1,0}+c_{2,0}-\Lambda_{\t_0}(Y_k))R(S_k)$$
with $R(\cdot)$ from (\ref{R()}). The minimizer $z_0$ of the quadratic function $q$ can easily be obtained as
\begin{eqnarray*}
z_0&=&-\Gamma_0^{-1}\frac{1}{n}\skn w(h_0^{-1}(S_k))(\hat{h}(Y_k)c_{1,0}+c_{2,0}-\Lambda_{\t_0}(Y_k))R(S_k)
\\&=& -\frac{c_{1,0}}{n}\Gamma_0^{-1}\skn\varphi(Z_k)-\frac{c_{1,0}}{\sqrt{n}}\Gamma_0^{-1}\beta+o_P(n^{-1/2})
\;=\;O_P(n^{-1/2}),
\end{eqnarray*}
where we have inserted the expansion from (\ref{exprh}) as well as (\ref{h0r0}) and use the definition $\beta=E[w(h_0^{-1}(S_1))r_0(h_0^{-1}(S_1))R(S_1)]$.
It is  also easy to see that (for some $\bar\gamma$)
\begin{align*}
T_n&=M_n(\tilde \gamma)=q(\hat{\gamma}-\gamma_0)+o_P(1)\geq \underset{z}{\min}\,q(z)+o_P(1) \\
&=q(z_0)+o_P(1)=M_n(\bar\gamma)+o_P(1)\geq T_n+o_P(1),
\end{align*}
so that it is sufficient to consider $q(z_0)=\min_z q(z)$ instead of $T_n$. 

Inserting the expansion for $z_0$ into $q(\cdot)$ as well as inserting the expansion for $\hat h$ from (\ref{exprh})  and (\ref{h0r0}) into $M_n(\gamma_0)$ gives
\begin{eqnarray*}
T_n&=& M_n(\hat{\gamma})=q(z_0)+o_P(1)\\
&=& c_{1,0}^2\Bigg(\frac{1}{n^2}\skn\sum_{i=1}^n\sjn w(h_0^{-1}(S_k))\psi(Z_i,U_k)\psi(Z_j,U_k)\\
&&\qquad{}+\frac{2}{n^{3/2}}\skn\sjn w(Y_k)r_0(Y_k)\psi(Z_j,U_k)
-\frac{1}{n}\sjn\skn\varphi(Z_j)^t\Gamma_0^{-1}\varphi(Z_k)\\
&&\qquad{}-\frac{2}{\sqrt{n}}\skn\varphi(Z_k)^t\Gamma_0^{-1}\beta
+E[w(h_0^{-1}(S))r_0(h_0^{-1}(S))^2]
-\beta^t\Gamma_0^{-1}\beta\Bigg) +o_P(1) .
\end{eqnarray*}
With some simple calculations of variances one shows that, after centering the multiple sums, those terms are negligible, where some of the indices coincide. Considering the (centred) multiple sums  with distinct indices only, for the resulting U-statistics Hoeff\-ding decompositions are applied (see, e.g.\ Section 1.6 in \citep{Lee1990}). Again with simple, but tedious calculations of variances one obtains the following dominating terms, 
\begin{eqnarray*}
T_n
&=& c_{1,0}^2\left(nU_n+b+n^{1/2}W_{0,n}+c\right),
\end{eqnarray*}
where $U_{n}$ is a U-statistic of order 2, i.e.\
$$U_n = \frac{1}{\binom{n}{2}}\sum_{i=1}^n\sum_{j=i+1}^n \zeta(Z_i,Z_j)$$
with degenerate kernel $$\zeta(z_1,z_2)=E[w(h_0^{-1}(S_3))\psi(Z_1,U_3)\psi(Z_2,U_3)\mid Z_1=z_1,Z_2=z_2]-\varphi(z_1)^t\Gamma_0^{-1}\varphi(z_2),$$
which coincides with $\zeta(z_1,z_2)$ from (\ref{zeta}). Further
$W_{0,n}=n^{-1}\sum_{i=1}^n \tilde{\zeta}(Z_i)$ with
\begin{eqnarray*}
\tilde{\zeta}(z) &=& 2E[w(h_0^{-1}(S_2))r_0(h_0^{-1}(S_2))\psi(Z_1,U_2)|Z_1=z]-2\varphi(z)^t\Gamma_0^{-1}\beta,
\end{eqnarray*}
which coincides with $\tilde\zeta(z)$ from (\ref{zeta-tilde}).
Furthermore
\begin{align}
b&=E[w(h_0^{-1}(S_1))\psi(Z_2,U_1)^2] -2E[\varphi(Z_1)^t\Gamma_0^{-1}\varphi(Z_1)] \;=\; E[\zeta(Z_1,Z_1)]\label{b}\\
c &=E[w(h_0^{-1}(S))r_0(h_0^{-1}(S))^2]\nonumber
-\beta^t\Gamma_0^{-1}\beta \;=\;E\left[w(h_0^{-1}(S_1))\bar{r}(S_1)^2\right]
\end{align}
with $\bar{r}$ from (\ref{bar-r}).\\
Note that $\zeta$ is symmetric. Hence, referring to \citet[p.~141]{WM1995} it can be written as
\begin{equation}\label{compactop}
\zeta(z_1,z_2)=\sum_{k=1}^\infty\lambda_k\rho_k(z_1)\rho_k(z_2)
\end{equation}
(in $L^2$ sense corresponding to the distribution $F_Z$) with notations from Theorem \ref{theolocalt}.
Referring to Remark \ref{rem-cov} $\zeta$ is positive semi-definite, which results in $\lambda_{k}\geq0,k\in\mathbb{N}$. From classical results on U-statistics, $nU_n$ converges to $\sum_{k=1}^\infty\lambda_{k}(W_k^2-1)$ in distribution (again with notations as in Theorem \ref{theolocalt}), see e.g.\ \cite{Lee1990}, Theorem 1 in Section 3.2.2. On the other hand, $n^{1/2}W_{0,n}$ converges to a normal distribution by the central limit theorem. As $U_n$ and $W_{0,n}$ are dependent, we have to go through some of the steps of the proof of Theorem 1 in \citet[p.~79]{Lee1990} to obtain the limiting distribution of $T_n$. 
\citet{Lee1990} uses the truncated sums $\sum_{k=1}^K\lambda_k\rho_k(Z_i)\rho_k(Z_j)\approx \zeta(Z_i,Z_j)$ (for large $K$)
to obtain the approximation
$$nU_n=\frac{1}{n-1}\sum_{i=1}^n \sum_{\overset{\scriptstyle j=1}{j\neq i}}^n\zeta(Z_i,Z_j) \approx \sum_{k=1}^K\lambda_k (W_{k,n}^2-v_{k,n})$$
with $W_{k,n}=n^{-1/2}\sum_{i=1}^n \rho_k(Z_i)$ and $v_{k,n}=n^{-1}\sum_{i=1}^n \rho_k^2(Z_i) =1+o_P(1)$ by the law of large numbers and the orthonormality of the eigenfunctions.  
Now to obtain convergence of $nU_n+n^{1/2}W_{0,n}$, note that applying the multivariate central limit theorem,
$(W_{0,n},W_{1,n},\dots,W_{K,n})^t$ converges in distribution to $(W_0,W_1,\dots,W_K)^t$ as defined in Theorem \ref{theolocalt}, for each $K$. 
Hence, by the continuous mapping theorem we obtain 
$$\sum_{k=1}^K\lambda_k(W_{k,n}^2-v_{k,n})+W_{0,n}\overset{\mathcal{D}}{\rightarrow}\sum_{k=1}^K\lambda_k(W_k^2-1)+W_0$$
for each $K$. Proceeding as in the proof of Theorem 1 in \citet[p.~79]{Lee1990} by letting $K\to\infty$, one obtains 
$\sum_{k=1}^\infty\lambda_k(W_k^2-1)+W_0$
as limit of $nU_n+n^{1/2}W_{0,n}$. 
Note further that (\ref{compactop}) especially leads to
$$\sum_{k=1}^\infty\lambda_k=\int\sum_{k=1}^\infty\lambda_k\rho_k(z_1)^2\,dP^Z(z_1)=E[\zeta(Z_1,Z_1)]=b,
$$
such that $nU_n+b+n^{1/2}W_{0,n}$ converges to $\sum_{k=1}^\infty\lambda_{k}W_k^2+W_0$, which completes the proof of Theorem \ref{theolocalt}. 
\hfill $\Box$

\bigskip

\noindent {\bf Proof of Theorem \ref{theoconsistency}.}
Note that the functions $f_{\gamma}(y)=w(y)(h(y)c_1+c_2-\Lambda_{\t}(y))^2$ are bounded, the parameter set $C_1\times C_2\times\T$ is compact and for every $y$ the map $\gamma=(c_1,c_2,\theta)\mapsto w(y)(h(y)c_1+c_2-\Lambda_{\t}(y))^2$ is continuous. Hence, following Lemma 6.1 in \citep{Wel2005} the class $\mathcal{F}=\{f_{\gamma}:\gamma\in C_1\times C_2\times\T\}$ is a Glivenko-Cantelli class. This leads to
\begin{align*}\frac{1}{n}T_n&=\underset{c_1,c_2,\t}{\min}\,\frac{1}{n}\sum_{k=1}^nw(Y_k)(\hat{h}(Y_k)c_1+c_2-\Lambda_{\t}(Y_k))^2
\\[0,2cm]&=\underset{c_1,c_2,\t}{\min}\,\frac{1}{n}\sum_{k=1}^nw(Y_k)(h(Y_k)c_1+c_2-\Lambda_{\t}(Y_k))^2+o_P(1)
\\[0,2cm]&=\underset{c_1,c_2,\t}{\min}\,E[w(Y)(h(Y)c_1+c_2-\Lambda_{\t}(Y))^2]+o_P(1)
\\[0,2cm]&=E[w(Y)(h(Y)c_1^*+c_2^*-\Lambda_{\t^*}(Y))^2]+o_P(1)
\end{align*}
for some $\gamma^*=(c_1^*,c_2^*,\t^*)^t\in\Upsilon$ (remind that $\Upsilon$ is compact). Under the fixed alternative $H_1$ one has $E[w(Y)(h(Y)c_1^*+c_2^*-\Lambda_{\t^*}(Y))^2]>0$, so that
$T_n\to\infty$. 
\hfill $\Box$

\bigskip

\noindent {\bf Proof of Theorem \ref{theointerchangedhyp}.}
The beginning of the proof is similar to the proof of Theorem \ref{theolocalt} and we will state the main differences. Again, we write $T_n=\min_{\gamma\in\Upsilon}M_n(\gamma)=M_n(\hat{\gamma})$ with $M_n$ as in (\ref{Mn}). Recall that by \ref{A6}, $\gamma_0=(c_{0,1},c_{0,2},\t_0)$ is the unique minimizer of
$M(c_1,c_2,\t)=E[w(Y_1)(h(Y_1)c_1+c_2-\Lambda_\t(Y_1))^2]$. 
We can again derive $\hat{\gamma}-\gamma_0=o_P(1)$ and further (\ref{Taylor1}) and (\ref{Gamma-limit}), but now with $\Gamma_0$ replaced by $\Gamma'$ as in (\ref{Gamma2}), which is the limit of  $n^{-1}\operatorname{Hess}M_n(\gamma_0)$. However, $M_n(\gamma_0)$ and $M_n(\hat{\gamma})$ are not bounded in probability under the model considered here. Instead we will show that
\begin{equation}\label{Mnrate}
M_n(\gamma_0)-M_n(\hat{\gamma})=O_P(n^{1/2})
\end{equation}
and thus we can derive from  (\ref{Taylor1}) that $\hat{\gamma}-\gamma_0=O_P(n^{-1/4})$. 
To obtain (\ref{Mnrate}), define
$$\tilde M_n(\gamma)=\sjn w(Y_j)(h(Y_j)c_1+c_2-\Lt(Y_j))^2$$
and let $\tilde M_n(\gamma^*)=\min_{\gamma\in\Upsilon}\tilde M_n(\gamma)$. Note that $E[M_n(\gamma)]=nM(\gamma)$. 
Now one can obtain
\begin{eqnarray}
\sup_{\gamma\in\Upsilon}|\tilde M_n(\gamma)-M_n(\gamma)| &=&O_P(n^{1/2}) \label{Mnexpansion}\\
 \sup_{\gamma\in\Upsilon} |\tilde M_n(\gamma)-nM(\gamma)| &=&O_P(n^{1/2}) \label{MnDonsker},
\end{eqnarray}
where for (\ref{Mnexpansion}) one applies (\ref{exprh}), whereas (\ref{MnDonsker}) holds because the empirical process $n^{-1/2}(\tilde M_n(\gamma)-nM(\gamma))$, $\gamma\in\Upsilon$, is Donsker. 
Because
\begin{eqnarray*}
|M_n(\hat{\gamma})-\tilde M_n(\gamma^*)| &=& |\inf_{\gamma\in\Upsilon} M_n(\gamma)-\inf_{\gamma\in\Upsilon} \tilde M_n(\gamma)| \;\leq\;\sup_{\gamma\in\Upsilon}|\tilde M_n(\gamma)-M_n(\gamma)| =O_P(n^{1/2}) 
\end{eqnarray*}
by (\ref{Mnexpansion}), to show (\ref{Mnrate}) it is sufficient to show 
\begin{equation}\label{Mnrate2}
M_n(\gamma_0)-\tilde M_n(\gamma^*)=O_P(n^{1/2}).
\end{equation}
To derive (\ref{Mnrate2}) note that
\begin{eqnarray*}
\tilde M_n(\gamma^*)-M_n(\gamma_0) &=& \inf_{\gamma\in\Upsilon}\tilde M_n(\gamma)-M_n(\gamma_0)\\
&\leq & \inf_{\gamma\in\Upsilon}nM(\gamma)+\sup_{\gamma\in\Upsilon} |\tilde M_n(\gamma)-nM(\gamma)| -M_n(\gamma_0)\\
&=&  nM(\gamma_0)+\sup_{\gamma\in\Upsilon} |\tilde M_n(\gamma)-nM(\gamma)|-M_n(\gamma_0)\\
&\leq & 2\sup_{\gamma\in\Upsilon} |\tilde M_n(\gamma)-nM(\gamma)|+ \sup_{\gamma\in\Upsilon}|\tilde M_n(\gamma)-M_n(\gamma)|\;=\;O_P(n^{1/2})
\end{eqnarray*}
by (\ref{Mnexpansion}) and (\ref{MnDonsker}). On the other hand 
\begin{eqnarray*}
M_n(\gamma_0)- \tilde M_n(\gamma^*)&=& M_n(\gamma_0) - \inf_\gamma \tilde M_n(\gamma)\\
&\leq& M_n(\gamma_0)-\inf_\gamma nM(\gamma) +\sup_\gamma|\tilde M_n(\gamma)-nM(\gamma)|\\
&=& M_n(\gamma_0)-nM(\gamma_0)+\sup_\gamma|\tilde M_n(\gamma)-nM(\gamma)|\\
&\leq& 2\sup_\gamma| M_n(\gamma)- \tilde M_n(\gamma)|+\sup_\gamma|\tilde M_n(\gamma)-nM(\gamma)| \;=\;O_P(n^{1/2})
\end{eqnarray*}
by (\ref{Mnexpansion}) and (\ref{MnDonsker}). 
Both inequalities together imply (\ref{Mnrate2}) and consequently (\ref{Mnrate}) holds. 

Again similar to the proof  of Theorem \ref{theolocalt} we obtain by Taylor expansion that, for all values $\gamma$ with $\gamma-\gamma_0=O_P(n^{-1/4})$,
\begin{eqnarray*}
M_n(\gamma) &=& M_n(\gamma_0)+(\gamma_0-\gamma)^t\nabla M_n(\gamma_0)+(\gamma_0-\gamma)^t\operatorname{Hess}M_n(\gamma_0)(\gamma_0-\gamma)+o_P(n^{1/2})\\
&=& q(\gamma_0-\gamma)+o_P(n^{1/2}),
\end{eqnarray*}
where the map $q$ is defined via
$$q(z)=M_n(\gamma_0)+nz^t\Gamma' z+2z^t\skn w(Y_k)(\hat{h}(Y_k)c_{1,0}+c_{2,0}-\Lambda_{\t_0}(Y_k))R_k$$
with $\Gamma'$ from (\ref{Gamma2}) and $\hat{R}_k= (\hat{h}(Y_k),1,-\dot{\Lambda}_{\t_0}(Y_k))^t$ ($k=1,\dots,n$).
The minimizer $z_0$ of the quadratic function $q$ can easily be obtained as
\begin{eqnarray*}
z_0&=&-\Gamma'^{-1}\frac{1}{n}\skn w(Y_k)(\hat{h}(Y_k)c_{1,0}+c_{2,0}-\Lambda_{\t_0}(Y_k))\hat{R}_k
\\&=& -\Gamma'^{-1}\frac{1}{n}\skn w(Y_k)(h(Y_k)c_{1,0}+c_{2,0}-\Lambda_{\t_0}(Y_k))R_k+O_P(n^{-1/2})
\;=\;O_P(n^{-1/2})
\end{eqnarray*}
with
$	R_k= (h(Y_k),1,-\dot{\Lambda}_{\t_0}(Y_k))^t$. 
To obtain the rate note that $E[w(Y_k)(h(Y_k)c_{1,0}+c_{2,0}-\Lambda_{\t_0}(Y_k))R_k]=0$ because $\gamma_0$ minimizes $M(\gamma)$ and thus one has $\nabla M(\gamma_0)=0$. 
As in the proof of Theorem \ref{theolocalt} it follows that instead of considering $T_n/n$ one can consider $q(z_0)/n$ to derive the limiting distribution. To this end note that
\begin{align*}
\frac{q(z_0)}{n}&=\frac{1}{n}\sum_{k=1}^nw(Y_k)(\hat{h}(Y_k)c_{1,0}+c_{2,0}-\Lambda_{\t_0}(Y_k))^2-z_0^t\Gamma' z_0
\\[0,2cm]&=\frac{c_{1,0}^2}{n}\sum_{k=1}^nw(Y_k)(\hat{h}(Y_k)-h(Y_k))^2+\frac{1}{n}\sum_{k=1}^nw(Y_k)(h(Y_k)c_{1,0}+c_{2,0}-\Lambda_{\t_0}(Y_k))^2
\\[0,2cm]&\quad+\frac{2c_{1,0}}{n}\sum_{k=1}^nw(Y_k)(\hat{h}(Y_k)-h(Y_k))(h(Y_k)c_{1,0}+c_{2,0}-\Lambda_{\t_0}(Y_k))+o_P(n^{-1/2}).
\end{align*}
The first term on the right hand side can be treated as in the proof of Theorem \ref{theolocalt} by inserting the expansion from (\ref{exprh}) to obtain the rate $o_P(n^{-1/2})$. The expectation of the second term on the right hand side is $M(\gamma_0)$.
Inserting the expansion from (\ref{exprh}) into the third term one obtains  
\begin{eqnarray}
n^{1/2}(T_n/n-M(\gamma_0))
&=&n^{1/2}(M_n(\gamma_0)/n-M(\gamma_0)) \nonumber\\
&&{}+\frac{2\sqrt{n}}{n^2}\sum_{j=1}^n\sum_{\overset{\scriptstyle k=1}{k\neq j}}^nw(Y_k)\psi(Z_j,U_k)\tilde R_k+o_P(1). \label{U-stat}
\end{eqnarray}
Applying a Hoeffding decomposition to the U-statistic term (\ref{U-stat}) one sees that the degenerate part is negligible and the dominating part is the (centred) linear term
$n^{-1/2}\sum_{i=1}^n \delta(Z_i)$ with $\delta$ from Theorem \ref{theointerchangedhyp}. The assertion of the theorem now follows from the classical central limit theorem.
\hfill $\Box$

\newpage

\noindent\Large{\bf Supplement to: \\``Specification testing in semi-parametric transformation models'' by Nick Kloodt, Natalie Neumeyer and Ingrid Van Keilegom}

\normalsize

\section{Bootstrap theory}\label{bootstrap-app}

In this section, we use the notations for the probability space as in section \ref{simulations}. The expectation with respect to $P_2^1(\omega,\cdot)$ is written as $E[\cdot|\omega]$. Note that the functions $h^*$ and $\psi^*$ depend on $\omega$ via the original sample. This is suppressed in the notation. 
We formulate the following additional assumptions. 

\begin{enumerate}[label=(\textbf{A8*})]
	\item\label{A8*} The following properties are meant conditional on the data $(Y_i,X_i),i=1,...,n,$ and thus define for fixed $n\in\mathbb{N}$ some subsets $A_n\in\mathcal{A}_1$ of $\Omega_1$, where these properties are valid. Thus, let $\omega\in A_n$, then we assume the following. 
	\begin{enumerate}[label=(\roman{*})]
		\item The transformation estimator $\hat{h}^*$ fulfils
		$$\underset{m\rightarrow\infty}{\limsup}\,P_2^1\bigg(\omega,\bigg\{\underset{y\in\mathcal{K}}{\sup}\,\bigg|\hat{h}^*(y)-h^*(y)-\frac{1}{m}\sum_{j=1}^m\psi^*(S_j^*,X_j^*,\mathcal{T}^*(y))\bigg|>\frac{\delta}{\sqrt{m}}\bigg\}\bigg)=0$$
		for all $\delta>0$, for some function $\psi^*$ and $h^*$ from Algorithm \ref{bootstrapalg}.
		\item For $\mathcal{U}_0$ from \ref{A8} and $\cyw$ from \ref{A2} we have $\mathcal{T}_{S^*}(h^*(\cyw))\subset \mathcal{U}_0$.
		\item Let $Z^*=(U^*,X^*)=(\mathcal{T}^*(Y^*),X^*)$. The function class $\{z\mapsto\psi^*(z,t):t\in\mathcal{U}_0\}$ is Donsker (for fixed $n$, but $m\rightarrow\infty$) with respect to $P_2^1(\omega,\cdot)^{Z^*}$ (distribution of $Z^*$ conditional on $\omega$) and
		$$E[\psi^*(Z_1^*,t)|\omega]=\int\psi^*(Z_1^*(\tilde{\omega}),t)P_2^1(\omega,d\tilde{\omega})=0\quad\textup{for all}\quad t\in\mathcal{U}_0.$$
		\item The fourth moment $E[w(h_0^{-1}(S_1^*))\psi^*(Z_1^*,U_1^*)^4|\omega]$ is finite and the conditional moments $E[w({h^*}^{-1}(S_1^*))\psi^*(Z_1^*,U_2^*)^2|Z_1^*=z,\omega]$ are locally bounded.
		\item For all compact sets $\mathcal{K}\subseteq\mathbb{R}$ we have
		\begin{equation*}\label{boundedpsi*}
		P_1\Big(\omega\in\Omega_1:\underset{y\in\mathcal{K},z\in\mathbb{R}^{d_X+1}}{\sup}\,|w(y)\psi^*(z,\mathcal{T}^*(y))|>\delta n\Big)=o(1)\quad\textup{for all }\delta>0.
		\end{equation*}
		\item One has $		\underset{y\rightarrow z}{\lim}\,E\left[|\psi^*(z,U_1^*)-\psi^*(y,U_1^*)|\,|\omega\right]=0$ for all $z$  in the support of $Z^*$.
	\end{enumerate} 
	For $A_n$ as defined above, we assume $P(A_n)\rightarrow1$ for $n\rightarrow\infty$.
\end{enumerate}
\begin{enumerate}[label=(\textbf{A9*})]
	\item\label{A9*} Define the distribution function of $Z^*$ for some $\omega\in\Omega_1$ by $F_{Z^*}(z)=P_2^1(\omega,\{Z^*\leq z\})$ and assume
	\begin{equation}\label{univconvdevFZ*}
	\underset{z\in\mathcal{R}^{d_X+1}}{\sup}\,|F_{Z^*}(z)-F_{{Z}}(z)|=o_P(1).
	\end{equation}
	Moreover, for all compact $\mathcal{K}\subseteq\mathbb{R}^{d_X+1}$ there exists an appropriate $C>0$, such that for $n\to\infty$
	\begin{equation}
	\underset{z\in\mathcal{K},s\in\mathbb{R}}{\sup}\,w((h^*)^{-1}(s))\psi^*(z,\mathcal{T}_{S^*}(s))\leq C+o_P(1)\label{boundedpsi^*}.
	\end{equation}
	Further,
	\begin{equation}
	w((h^*)^{-1}(s))(\psi^*(z,\mathcal{T}_{S^*}(s))-\psi(z,\mathcal{T}_{S^*}(s)))=o_P(1)\label{convpsi^*}
	\end{equation}
	for all $z\in\mathbb{R}^{d_X+1}$, $s\in\mathbb{R}$ and for $\psi$ from \ref{A8} for $n\to\infty$. 
\end{enumerate}


\subsection{Main bootstrap results and proofs}

\begin{theo}\label{theobootstrap}
	Assume $H_0$,\ref{A1}--\ref{A8},\ref{A8*},\ref{A9*}. Then, the bootstrap statistic $T_{n,m}^*$ computed by Algorithm \ref{bootstrapalg} fulfils (\ref{bootstrapH0}). If $q_{\alpha}^*$ denotes for all $\alpha\in(0,1)$ the corresponding bootstrap quantile described in Algorithm \ref{bootstrapalg}, it is
	$$P_1\Big(\omega\in\Omega_1:\underset{m\rightarrow\infty}{\limsup}\,|q_{\alpha}^*-q_{\alpha}|>\delta\Big)=o(1)$$
	for all $\delta>0$.
\end{theo}

\noindent {\bf Proof of Theorem \ref{theobootstrap}.}
As in the proof of Lemma \ref{lemmabootstrap}, the conditional distribution and expectation of $(Y_1^*,X_1^*),...,(Y_m^*,X_m^*)$ given $(Y_1,X_1),...,(Y_n,X_n)$ are denoted by $P^*$ and $E^*$, respectively. Consider $\omega\in A_n$ with $A_n$ from \ref{A8*}. The proof can be divided into two parts: First, the uniform convergence of some bootstrap components appearing in the asymptotic distribution of the bootstrap test statistic is proven and second, the assertion itself is shown by the convergence of the conditional distribution functions in probability. Referring to the definition of $h^*$, the following condition \ref{A6*} is valid.
\begin{enumerate}[label=(\textbf{A6*})]
	\item\label{A6*} With probability converging to one, minimizing the function $M^*:\Upsilon\rightarrow\R,$
	\begin{align*}
	\gamma&\ \:=(c_1,c_2,\theta)\mapsto E^*\big[w(Y^*)(h^*(Y^*)c_1+c_2-\Lt(Y^*))^2\big]
	\end{align*}
	leads to a unique solution
	$$\gamma^*=(c^*_{1,0},c^*_{2,0},\theta^*_0)=\big(\Lambda_{\hat{\theta}}(1)-\Lambda_{\hat{\theta}}(0),\Lambda_{\hat{\theta}}(0),\hat{\theta}\big)$$
	in the interior of $\Upsilon$.
\end{enumerate}
Here, uniqueness follows due to $E^*[w(Y^*)(h^*(Y^*)(\Lambda_{\hat{\theta}}(1)-\Lambda_{\hat{\theta}}(0))+\Lambda_{\hat{\theta}}(0)-\Lambda_{\hat{\t}}(Y^*))^2]=0$ from \ref{A6}. With the notations
\begin{align}
R^*(s)&:=(s,1,-\dot{\Lambda}_{\hat{\t}}((h^*)^{-1}(s)))^t,\nonumber\\[0,2cm]
\Gamma^*&:=E^*[w((h^*)^{-1}(S_1^*))R^*(S_1^*)R^*(S_1^*)^t],\nonumber\\[0,2cm]
\varphi^*(z)&:=E^*[w((h^*)^{-1}(S_2^*))\psi^*(Z_1^*,U_2^*)R^*(S_2^*)\mid Z_1^*=z],\nonumber
\end{align}
a function $\zeta^*$ can be defined as
\begin{align}
\zeta^*(z_1,z_2)&:=E^*\Big[w((h^*)^{-1}(S_3^*))\big(\psi^*(Z_1^*,U_3^*)-\varphi^*(Z_1^*)^t(\Gamma^*)^{-1}R^*(S_3^*)\big) \nonumber\\[0,2cm]
&\qquad\times\big(\psi^*(Z_2^*,U_3^*)-\varphi^*(Z_2^*)^t(\Gamma^*)^{-1}R^*(S_3^*)\big)\mid Z_1^*=z_1,Z_2^*=z_2\Big]\nonumber.
\end{align}
Moreover, define
\begin{equation}\label{asympdistH0bootstrap}
T_n^*={c_{1,0}^*}^2\sum_{k=1}^\infty\lambda_{k}^*W_k^2\quad\textup{with}\quad c_{1,0}^*=\Lambda_{\hat{\t}}(1)-\Lambda_{\hat{\t}}(0)
\end{equation}
and $b^*=E[\zeta^*(Z_1^*,Z_1^*)|\omega]$, where $W_1,W_2,...$ are independent and standard normally distributed and $\lambda_1^*,\lambda_2^*,...$ are the eigenvalues of the operator
$$K^*\rho(z_1):=\int\rho(z_2)\zeta^*(z_1,z_2)\,dF_{Z^*}(z_2).$$
Therefore for fixed $n$ and conditional on $(X_1,Y_1),\dots,(X_n,Y_n)$, one can proceed exactly  as in the proof of Theorem \ref{theolocalt} to obtain that $T_{n,m}^*$ converges in distribution to $T_n^*$ for $m\rightarrow\infty$. We  have for $n\to\infty$
\begin{equation}\label{approvTnm*zeta*}
P_1\bigg(\omega\in\Omega_1:P_2^1\bigg(\omega,\bigg\{\underset{m\rightarrow\infty}{\limsup}\,\bigg|T_{n,m}^*-\frac{1}{m-1}\sum_{i=1}^m\sum_{\overset{\scriptstyle j=1}{j\neq i}}^m\zeta^*(Z_i^*,Z_j^*)-b^*\bigg|>0\bigg\}\bigg)>0\bigg)=o(1),
\end{equation}
as well as
$$P_1\Big(\omega\in\Omega_1:\underset{m\rightarrow\infty}{\limsup}\,\underset{t\in\mathbb{R}}{\sup}\,\big|P_2^1(\omega,\{T_{n,m}^*\leq t\})-P_2^1(\omega,\{T_n^*\leq t\})\big|=0\Big)=1+o(1).$$

\textbf{Convergence of the bootstrap components:} In the following, the convergence in probability of $h^*,R^*,\Gamma^*,\varphi^*,\zeta^*,b^*$ to $h_0$ (the true transformation under $H_0$), $R(s)$ from (\ref{R()}), $\Gamma_0$ from (\ref{Gamma}), $\varphi$ from (\ref{phi()}), $\zeta$ from (\ref{zeta}) and $b$ from (\ref{b}) is shown.\\
One has $\hat{\t}=\t_0+o_P(1)$ for $n\to\infty$. From (\ref{defY*}) follows $h^*=h+o_P(1)$ for $n\to\infty$ uniformly on compact sets and thus $w((h^*)^{-1}(s))R^*(s)=w(h^{-1}(s))R(s)+o_P(1)$ uniformly in $s\in\mathbb{R}$. Further, there exists some $\tilde{C}>0$, such that $h^*$ is bijective on $\mathcal{Y}_w$ and $|h^*|<\tilde{C}$ as well as $|(h^*)^{-1}|<\tilde{C}$ on $\cyw$ or $h^*(\cyw)$ with probability converging to one, that is
$$P_1\Big(\omega\in\Omega_1:\ h^*\textup{ bijective,}\quad|h^*(y)|\leq\tilde{C}\ \forall y\in\mathcal{Y}_w,\quad |(h^*)^{-1}(s)|\leq\tilde{C}\ \forall y\in h^*(\mathcal{Y}_w)\Big)=1+o(1)$$
for $n\to\infty$. 
This in turn means that (see Remark \ref{rembootsassumpCKC} for a possible adjustment of $h^*$)
$$P_1\Big(\omega\in\Omega_1:w((h^*)^{-1}(s))||R^*(s)||^2\leq C\ \forall s\in\mathbb{R}\Big)=1+o(1)$$
for $n\to\infty$
for some sufficiently large $C>0$. Let $f_{S},f_{S^*}$ denote the densities of $S$ and $S^*$, respectively, conditioned on $(Y_i,X_i),i=1,...,n$. The dominated convergence theorem leads to (the inequality is meant componentwise)
\begin{align*}
&P_1\big(\omega\in\Omega_1:|\Gamma^*-\Gamma|>\delta\big)
\\[0,2cm]&=P_1\big(\omega\in\Omega_1:\big|E[w((h^*)^{-1}(S_1^*))R^*(S_1^*)R^*(S_1^*)^t|\omega]-E[w(h^{-1}(S_1))R(S_1)R(S_1)^t]\big|>\delta\big)
\\[0,2cm]&=P_1\bigg(\omega\in\Omega_1:\bigg|\int w((h^*)^{-1}(s))R^*(s)R^*(s)^tf_{S^*}(s)\,ds-\int w(h^{-1}(s))R(s)R(s)^tf_{S}(s)\,ds\bigg|>\delta\bigg)
\\[0,2cm]&\leq P_1\bigg(\omega\in\Omega_1:\bigg|\int w((h^*)^{-1}(s))R^*(s)R^*(s)^tf_{S^*}(s)\,ds-\int w(h^{-1}(s))R(s)R(s)^tf_{S}(s)\,ds\bigg|>\delta,
\\[0,2cm]&\qquad w((h^*)^{-1}(s))||R^*(s)||^2\leq C\ \forall s\in\mathbb{R}\bigg)+P_1\big(w((h^*)^{-1}(s))||R^*(s)||^2>C\textup{ for some }s\in\mathbb{R}\big)
\\[0,2cm]&=o(1)
\end{align*}
for all $\delta>0$. Consequently,
$$\Gamma^*=\Gamma+o_{P_1}(1)$$
for $n\rightarrow\infty$. Due to part (\ref{boundedpsi^*}) of \ref{A9*}, the map $(Z_1^*,U_2^*)\mapsto w((h^*)^{-1}(S_2^*))\psi^*(Z_1^*,U_2^*)$ is bounded by some constant $C$ uniformly over compact sets with probability converging to one. Together with the dominated convergence theorem, this leads to boundedness of $\varphi^*$ as well as $\varphi^*(z)=\varphi(z)+o_P(1)$ and finally boundedness of $\zeta^*$ and
\begin{align*}
\zeta^*(z_1,z_2)&=\int w((h^*)^{-1}(s))\big(\psi^*(z_1,\mathcal{T}_{S^*}(s))-\varphi^*(z_1)^t(\Gamma^*)^{-1}R^*(s)\big)
\\[0,2cm]&\quad\quad\big(\psi^*(z_2,\mathcal{T}_{S^*}(s))-\varphi^*(z_2)^t(\Gamma^*)^{-1}R^*(s)\big)f_{S^*}(s)\,ds
\\[0,2cm]&=\int w(h^{-1}(s))\big(\psi(z_1,\mathcal{T}_{S}(s))-\varphi(z_1)^t\Gamma^{-1}R(s)\big)\big(\psi(z_2,\mathcal{T}_{S}(s))-\varphi(z_2)^t\Gamma^{-1}R(s)\big)f_{S}(s)\,ds
\\[0,2cm]&\quad+o_{P_1}(1)
\\[0,2cm]&=\zeta(z_1,z_2)+o_{P_1}(1)
\end{align*}
for all $z_1,z_2\in\mathbb{R}^{d_X+1}$. Additionally, one has
$$P_1\bigg(\omega\in\Omega_1:\underset{z_1,z_2\in\mathbb{R}^{d_X+1}}{\sup}\,|\zeta^*(z_1,z_2)|>C\bigg)\rightarrow0$$
for some $C>0$, so that for all $\delta>0$
\begin{align*}
&P_1\big(\omega\in\Omega_1:|b^*-b|>\delta\big)
\\[0,2cm]&=P_1\big(\omega\in\Omega_1\big|E[\zeta^*(Z_1^*,Z_1^*)|\omega]-E[\zeta(Z_1,Z_1)]\big|>\delta\big)
\\[0,2cm]&\leq P_1\bigg(\omega\in\Omega_1:\big|E[\zeta^*(Z_1^*,Z_1^*)|\omega]-E[\zeta(Z_1,Z_1)]\big|>\delta,\underset{z_1,z_2\in\mathbb{R}^{d_X+1}}{\sup}\,|\zeta^*(z_1,z_2)|\leq C\bigg)
\\[0,2cm]&\quad+P_1\bigg(\omega\in\Omega_1:\underset{z_1,z_2\in\mathbb{R}^{d_X+1}}{\sup}\,|\zeta^*(z_1,z_2)|>C\bigg)
\\[0,2cm]&=o(1),
\end{align*}
that is $b^*=b+o_{P_1}(1)$ for $n\to\infty$. Now, all ingredients to prove the convergence of the distribution functions $F_{T_n^*}(t)=P^*(T_n^*\leq t)$ to $P(T\leq t)$ in probability have been presented.
\\[0,2cm]
\textbf{Convergence of the distribution functions:} Let $t\in\mathbb{R},\varepsilon>0$ be arbitrary and $\tilde{\varepsilon},\delta_{\tilde{\varepsilon}}>0,M_{\tilde{\varepsilon}}\in\mathbb{N}$ such that
$$P(T>t-2\tilde{\varepsilon} c_{1,0}^2)-P(T>t+2\tilde{\varepsilon} c_{1,0}^2)+2\tilde{\varepsilon}\leq\frac{\varepsilon}{2},$$
\begin{equation}\label{convdistfunct}
\underset{s\in\mathbb{R}}{\sup}\,\bigg|P\bigg(c_{1,0}^2\bigg(\frac{m}{\binom{m}{2}}\sum_{i=1}^m\sum_{j=i+1}^m\zeta(Z_i,Z_j)+b\bigg)>s\bigg)-P(T>s)\bigg|\leq\tilde{\varepsilon}\quad\textup{for all }m\geq M_{\tilde{\varepsilon}}
\end{equation}
(due to the proof of Theorem \ref{theolocalt} for sample size $n$ replaced by $m$)
and
$$|c_{1,0}^*-c_{1,0}|+|b^*-b|\leq\delta_{\tilde{\varepsilon}}\quad\Rightarrow\quad\bigg|\frac{t}{{c_{1,0}^*}^2}-b^*-\frac{t}{c_{1,0}^2}+b\bigg|\leq\tilde{\varepsilon}.$$
For a moment consider $m$ as fixed and define
$$\mathcal{K}_{\tilde{\varepsilon},m}:=\Big\{(z_1,z_2)\in\mathbb{R}^{2(d+1)}:|\zeta^*(z_1,z_2)-\zeta(z_1,z_2)|\leq\frac{\tilde{\varepsilon}}{m}\Big\}.$$
Then, $P((Z_1^*,Z_2^*)\in\mathcal{K}_{\tilde{\varepsilon},m})\overset{n\rightarrow\infty}{\longrightarrow}1$. For all $m\geq M_{\varepsilon}$ one has
\begin{align*}
&P^*\bigg({c_{1,0}^*}^2\bigg(\frac{m}{\binom{m}{2}}\sum_{i=1}^m\sum_{j=i+1}^m\zeta^*(Z_i^*,Z_j^*)+b^*\bigg)>t\bigg)
\\[0,2cm]&\leq P^*\bigg({c_{1,0}^*}^2\bigg(\frac{m}{\binom{m}{2}}\sum_{i=1}^m\sum_{j=i+1}^m\zeta^*(Z_i^*,Z_j^*)+b^*\bigg)>t,(Z_i^*,Z_j^*)\in\mathcal{K}_{\tilde{\varepsilon},m}\ \forall i\neq j\in\{1,...,m\}\bigg)
\\[0,2cm]&\quad+P^*\big(\exists i\neq j\in\{1,...,m\}:(Z_i^*,Z_j^*)\notin\mathcal{K}_{\tilde{\varepsilon},m}\big)
\\[0,2cm]&\leq P^*\bigg(\frac{m}{\binom{m}{2}}\sum_{i=1}^m\sum_{j=i+1}^m\zeta(Z_i^*,Z_j^*)>\frac{t}{c_{1,0}^2}-b-2\tilde{\varepsilon}\bigg)+P^*(|c_{1,0}^*-c_{1,0}|+|b^*-b|>\delta_{\tilde{\varepsilon}})
\\[0,2cm]&\quad+P^*\big(\exists i\neq j\in\{1,...,m\}:(Z_i^*,Z_j^*)\notin\mathcal{K}_{\tilde{\varepsilon},m}\big)
\\[0,2cm]&\leq\int\dots\int I_{\big\{\frac{m}{\binom{m}{2}}\sum_{i=1}^m\sum_{j=i+1}^m\zeta(z_i,z_j)>\frac{t}{c_{1,0}^2}-b-2\tilde{\varepsilon}\big\}}\,dF_{Z^*}(z_1)\dots dF_{Z^*}(dz_m)
\\[0,2cm]&\quad+P^*(|c_{1,0}^*-c_{1,0}|+|b^*-b|>\delta_{\tilde{\varepsilon}})+P^*\big(\exists i\neq j\in\{1,...,m\}:(Z_i^*,Z_j^*)\notin\mathcal{K}_{\tilde{\varepsilon},m}\big)
\\[0,2cm]&=\int\dots\int I_{\big\{\frac{m}{\binom{m}{2}}\sum_{i=1}^m\sum_{j=i+1}^m\zeta(z_i,z_j)>\frac{t}{c_{1,0}^2}-b-2\tilde{\varepsilon}\big\}}\,dF_{Z}(z_1)\dots dF_{Z}(dz_m)+o_P(1)
\\[0,2cm]&\quad+P^*(|c_{1,0}^*-c_{1,0}|+|b^*-b|>\delta_{\tilde{\varepsilon}})+P^*\big(\exists i\neq j\in\{1,...,m\}:(Z_i^*,Z_j^*)\notin\mathcal{K}_{\tilde{\varepsilon},m}\big)
\\[0,2cm]&\leq P(T>t-2\tilde{\varepsilon}c_{1,0}^2)+\tilde{\varepsilon}+o_P(1).
\end{align*}
Here, the equality follows from  the Portmanteau theorem due to (\ref{univconvdevFZ*}).  The same reasoning leads to
$$P^*\bigg({c_{1,0}^*}^2\bigg(\frac{m}{\binom{m}{2}}\sum_{i=1}^m\sum_{j=i+1}^m\zeta^*(Z_i^*,Z_j^*)+b^*\bigg)>t\bigg)\geq P(T>t+2\tilde{\varepsilon}c_{1,0}^2)-\tilde{\varepsilon}+o_P(1)$$
for all $\varepsilon>0$ and thus
\begin{equation}\label{bootscondconvfixedm}
\bigg|P^*\bigg({c_{1,0}^*}^2\bigg(\frac{m}{\binom{m}{2}}\sum_{i=1}^m\sum_{j=i+1}^m\zeta^*(Z_i^*,Z_j^*)+b^*\bigg)\leq t\bigg)-F_{T}(t)\bigg|\leq\frac{\varepsilon}{2}+o_P(1).
\end{equation}
Let $M\in\mathbb{N}$ and define
$$A_{M,\varepsilon}:=\bigg\{\omega\in\Omega_1:\forall m\geq M: \bigg|P^*\bigg({c_{1,0}^*}^2\bigg(\frac{m}{\binom{m}{2}}\sum_{i=1}^m\sum_{j=i+1}^m\zeta^*(Z_i^*,Z_j^*)+b^*\bigg)>t\bigg)-P^*(T_n^*>t)\bigg|<\varepsilon\bigg\}$$
as well as
$$B_{n}:=\bigg\{\omega\in\Omega_1:\underset{m\rightarrow\infty}{\limsup}\,\bigg|P^*\bigg({c_{1,0}^*}^2\bigg(\frac{m}{\binom{m}{2}}\sum_{i=1}^m\sum_{j=i+1}^m\zeta^*(Z_i^*,Z_j^*)+b^*\bigg)>t\bigg)-P^*(T_n^*>t)\bigg|=0\bigg\}.$$
Note that $P_1(B_{n})\overset{n\rightarrow\infty}{\longrightarrow}1$ because of (\ref{approvTnm*zeta*}) and $P_1(B_{n}\cap A_{M,\varepsilon})\overset{M\rightarrow\infty}{\longrightarrow}P_1(B_{n})$. Let $N\in\mathbb{N}$ such that $P_1(B_{N})\geq 1-\varepsilon$ and let $M_{\varepsilon}$ fulfil (\ref{convdistfunct}) and $P_1(B_{N})\leq P_1(B_{N}\cap A_{M_{\varepsilon},\varepsilon})+\varepsilon$. Then, it is
\begin{small}
	\begin{align*}
	&\underset{n\rightarrow\infty}{\limsup}\,P_1(|F_{T_n^*}(t)-F_T(t)|>3\varepsilon)
	\\[0,2cm]&\ \ =\underset{n\rightarrow\infty}{\limsup}\,P_1\bigg(\underset{m\rightarrow\infty}{\limsup}\,\bigg|P^*\bigg({c_{1,0}^*}^2\bigg(\frac{m}{\binom{m}{2}}\sum_{i=1}^m\sum_{j=i+1}^m\zeta^*(Z_i^*,Z_j^*)+b^*\bigg)>t\bigg)-P(T>t)\bigg|>3\varepsilon,B_{n}\bigg)
	\\[0,2cm]&\ \ =\underset{n\rightarrow\infty}{\limsup}\,\underset{M\rightarrow\infty}{\lim}\,P_1\bigg(\underset{m\geq M}{\sup}\,\bigg|P^*\bigg({c_{1,0}^*}^2\bigg(\frac{m}{\binom{m}{2}}\sum_{i=1}^m\sum_{j=i+1}^m\zeta^*(Z_i^*,Z_j^*)+b^*\bigg)>t\bigg)-P(T>t)\bigg|>3\varepsilon,B_{n}\cap A_{M,\varepsilon}\bigg)
	\\[0,2cm]&\ \ \leq\underset{n\rightarrow\infty}{\limsup}\,\underset{M\rightarrow\infty}{\lim}\,P_1\bigg(\underset{m\geq M}{\sup}\,\bigg|P^*\bigg({c_{1,0}^*}^2\bigg(\frac{m}{\binom{m}{2}}\sum_{i=1}^m\sum_{j=i+1}^m\zeta^*(Z_i^*,Z_j^*)+b^*\bigg)>t\bigg)-P(T>t)\bigg|>3\varepsilon
	\\[0,2cm]&\ \ \quad\quad\quad\quad\quad\quad\quad\quad\quad,B_{N}\cap A_{M_{\varepsilon},\varepsilon}\bigg)+2\varepsilon
	\\[0,2cm]&\ \ \leq\underset{n\rightarrow\infty}{\limsup}\,P_1\bigg(\bigg|P^*\bigg({c_{1,0}^*}^2\bigg(\frac{M_{\varepsilon}}{\binom{M_{\varepsilon}}{2}}\sum_{i=1}^{M_{\varepsilon}}\sum_{j=i+1}^{M_{\varepsilon}}\zeta^*(Z_i^*,Z_j^*)+b^*\bigg)>t\bigg)-P(T>t)\bigg|>\varepsilon,B_{N}\cap A_{M_{\varepsilon},\varepsilon}\bigg)+2\varepsilon
	\\[0,2cm]&\ \ \leq\underset{n\rightarrow\infty}{\limsup}\,P_1\bigg(\bigg|P^*\bigg({c_{1,0}^*}^2\bigg(\frac{M_{\varepsilon}}{\binom{M_{\varepsilon}}{2}}\sum_{i=1}^{M_{\varepsilon}}\sum_{j=i+1}^{M_{\varepsilon}}\zeta^*(Z_i^*,Z_j^*)+b^*\bigg)>t\bigg)-P(T>t)\bigg|>\varepsilon\bigg)+2\varepsilon
	\\[0,2cm]&\overset{(\ref{bootscondconvfixedm})}{=}2\varepsilon.
	\end{align*}
\end{small}
Since $\varepsilon>0$ can be chosen arbitrarily small, one has
$$P_1\big(|F_{T_n^*}(t)-F_{T}(t)|>\varepsilon\big)=o(1)\quad\textup{for all }t\in\mathbb{R},\varepsilon>0.$$
In total, (\ref{bootstrapH0}) was proven, that is,
$$P_1\Big(\omega\in\Omega_1:\underset{m\rightarrow\infty}{\limsup}\,\big|P_2^1(\omega,\{T_{n,m}^*\leq q\})-P(T_n\leq q)\big|>\delta\Big)=o(1)\quad\textup{for all }\delta>0.$$
It remains to deduce
$$P_1\Big(\omega\in\Omega_1:\underset{m\rightarrow\infty}{\limsup}\,|q_{\alpha}^*-q_{\alpha}|>\delta\Big)=o(1)\quad\textup{for all }\delta>0$$
from this. Let $\delta>0$ be arbitrarily small and let $q_{\alpha}$ be the $\alpha$-quantile of $T$ and define
$$\varepsilon=\min\bigg(\frac{\alpha-P(T_n\leq q_{\alpha}-\delta)}{2},\frac{P(T_n\leq q_{\alpha}+\delta)-\alpha}{2}\bigg)>0$$
which depends on $n$, but converges to a positive value. 
Then, if
$$\big|P_2^1(\omega,\{T_{n,m}^*\leq q_{\alpha}-\delta\})-P(T_n\leq q_{\alpha}-\delta)\big|\leq\varepsilon,$$
one has
$$P_2^1(\omega,\{T_{n,m}^*\leq q_{\alpha}-\delta\})<\alpha$$
for $n$ sufficiently large, 
and thus $q_{\alpha}^*>q_{\alpha}-\delta$. Analogously, $\big|P_2^1(\omega,\{T_{n,m}^*\leq q_{\alpha}+\delta\})-P(T_n\leq q_{\alpha}+\delta)\big|\leq\varepsilon$ implies $P_2^1(\omega,\{T_{n,m}^*\leq q_{\alpha}+\delta\})>\alpha$ and $q_{\alpha}^*<q_{\alpha}+\delta$, so that in total
\begin{align*}
&P_1\Big(\omega\in\Omega_1:\underset{m\rightarrow\infty}{\limsup}\,|q_{\alpha}^*-q_{\alpha}|\leq\delta\Big)
\\[0,2cm]&\geq P_1\Big(\omega\in\Omega_1:\underset{m\rightarrow\infty}{\limsup}\,\big|P_2^1(\omega,\{T_{n,m}^*\leq q_{\alpha}-\delta\})-P(T_n\leq q_{\alpha}-\delta)\big|\leq\varepsilon,
\\[0,2cm]&\quad\quad\underset{m\rightarrow\infty}{\limsup}\,\big|P_2^1(\omega,\{T_{n,m}^*\leq q_{\alpha}+\delta\})-P(T_n\leq q_{\alpha}+\delta)\big|\leq\varepsilon\Big)
\\[0,2cm]&=1+o(1).
\end{align*}
\hfill $\Box$

\begin{theo}\label{theobootstrapalternative}
	Assume $H_1$,\ref{A1}--\ref{A4},\ref{A6'},\ref{A8*}. Then, the bootstrap statistic $T_{n,m}^*$ computed by Algorithm \ref{bootstrapalg} fulfils (\ref{bootstrapH1}). If $q_{\alpha}^*$ denotes for all $\alpha\in(0,1)$ the corresponding bootstrap quantile described in Algorithm \ref{bootstrapalg}, it is
	$$P_1\Big(\omega\in\Omega_1:T_n>\underset{m\rightarrow\infty}{\limsup}\,q_{\alpha}^*\Big)=1+o(1).$$
\end{theo}

\noindent {\bf Proof of Theorem \ref{theobootstrapalternative}.}
Borrowing the notations for $R^*,\Gamma^*$ and $\varphi^*$ from the proof of Theorem \ref{theobootstrap} remind that $\zeta^*$ can be written as
\begin{align}
\zeta^*(z_1,z_2)&:=E^*\Big[w((h^*)^{-1}(S_3^*))\big(\psi^*(Z_1^*,U_3^*)-\varphi^*(Z_1^*)^t(\Gamma^*)^{-1}R^*(S_3^*)\big) \nonumber\\[0,2cm]
&\qquad\times\big(\psi^*(Z_2^*,U_3^*)-\varphi^*(Z_2^*)^t(\Gamma^*)^{-1}R^*(S_3^*)\big)\mid Z_1^*=z_1,Z_2^*=z_2\Big]\nonumber.
\end{align}
Similar to the proof of Theorem \ref{theobootstrap} Assumption \ref{A6'} leads to \ref{A6*} from the proof of Theorem \ref{theobootstrap}. As before, Equations (\ref{asympdistH0bootstrap}) and (\ref{approvTnm*zeta*}) remain valid. Especially, it is
$$P_1\bigg(\omega\in\Omega_1:P_2^1\bigg(\omega,\bigg\{\underset{m\rightarrow\infty}{\limsup}\,\bigg|T_{n,m}^*-\frac{1}{m-1}\sum_{i=1}^m\sum_{\overset{\scriptstyle j=1}{j\neq i}}^m\zeta^*(Z_i^*,Z_j^*)-b^*\bigg|>0\bigg\}\bigg)>0\bigg)=o(1)$$
for $b^*=E[\zeta^*(Z_1^*,Z_1^*)|\omega]$.
Since \ref{A8*} ensures that
$$P_1\Big(\omega\in\Omega_1:\underset{y\in\mathcal{K},z\in\mathbb{R}^{d_X+1}}{\sup}\,|w(y)\psi^*(z,\mathcal{T}^*(y))|>\delta n\Big)=o(1)\quad\textup{for all }\delta>0$$
and  all compact sets $\mathcal{K}\subseteq\mathbb{R}$, it is
$$P_1\Big(\omega\in\Omega_1:\underset{z_1,z_2\in\mathbb{R}^{d_X+1}}{\sup}|\zeta^*(z_1,z_2)|>\delta n\Big)=o(1)\quad\textup{for all }\delta>0$$
(boundedness of $R^*$ and $\Gamma^*$ follows as in the proof of Theorem \ref{theobootstrap}). The same reasoning as in the proof of Theorem \ref{theolocalt} leads to
$$E\bigg[\bigg(\frac{1}{m-1}\sum_{i=1}^m\sum_{\overset{\scriptstyle j=1}{j\neq i}}^m\zeta^*(Z_i^*,Z_j^*)\bigg)^2\,\Big|\omega\bigg]=\frac{2m^2}{(m-1)^2}E\big[\zeta^*(Z_1^*,Z^*_2)^2\,|\omega\big],$$
so that
$$\underset{m\rightarrow\infty}{\limsup}\bigg|\frac{1}{m-1}\sum_{i=1}^m\sum_{\overset{\scriptstyle j=1}{j\neq i}}^m\zeta^*(Z_i^*,Z_j^*)+b^*\bigg|=o_P(n)$$
and thus $\underset{m\rightarrow\infty}{\limsup}|T_n^*|=o_P(n)$. Referring to Theorem \ref{theointerchangedhyp}, one has
$$\frac{1}{n}T_n=M(\gamma_0)+o_P(1)$$
with $M(\gamma_0)>0$. Consequently,
\begin{align*}
&P_1\Big(\omega\in\Omega_1:\underset{m\rightarrow\infty}{\limsup}\,P_2^1\big(\omega,\{T_n\leq T_{n,m}^*\}\big)>0\Big)
\\[0,2cm]&=P_1\Bigg(\omega\in\Omega_1:\underset{m\rightarrow\infty}{\limsup}\,P_2^1\bigg(\omega,\bigg\{\frac{1}{n}T_n\leq \frac{1}{n}T_{n,m}^*\bigg\}\bigg)>0\Bigg)
\\[0,2cm]&\leq P_1\Bigg(\omega\in\Omega_1:\underset{m\rightarrow\infty}{\limsup}\,P_2^1\bigg(\omega,\bigg\{\frac{M(\gamma_0)}{2}\leq \frac{1}{n}T_{n,m}^*\bigg\}\bigg)>0\Bigg)+P_1\bigg(\omega\in\Omega_1:\frac{T_n}{n}<\frac{M(\gamma_0)}{2}\bigg)
\\[0,2cm]&=o(1).
\end{align*}
In total, (\ref{bootstrapH1}) was proven, that is,
$$P_1\Big(\omega\in\Omega_1:\underset{m\rightarrow\infty}{\limsup}\,P_2^1(\omega,\{T_n<T_{n,m}^*\})>\delta\Big)=o(1)\quad\textup{for all }\delta>0.$$
It remains to deduce
$$P_1\Big(\omega\in\Omega_1:T_n>\underset{m\rightarrow\infty}{\limsup}\,q_{\alpha}^*\Big)=1+o(1)$$
for all $\alpha\in(0,1)$ from this. Let $\alpha\in(0,1)$ be fixed. If $T_n>q_{\alpha}^*$, one has
$$P_2^1(\omega,\{T_n<T_{n,m}^*\})<P_2^1(\omega,\{q_{\alpha}^*<T_{n,m}^*\})\leq1-\alpha,$$
so that
\begin{align*}
P_1\Big(\omega\in\Omega_1:T_n>\underset{m\rightarrow\infty}{\limsup}\,q_{\alpha}^*\Big)&\geq P_1\bigg(\omega\in\Omega_1:\underset{m\rightarrow\infty}{\limsup}\,P_2^1(\omega,\{T_n<T_{n,m}^*\})\leq1-\alpha\bigg)
\\[0,2cm]&=1+o(1).
\end{align*}
\hfill $\Box$

\subsection{Nonparametric transformation estimation in the bootstrap case}

The aim of this subsection is to show that the estimating approach developed by \citet{CvK2018} can be applied in this context. To this end, for the estimator $\hat h$ from Appendix \ref{estimation-of-h} respectively its bootstrap analog $\hat h^*$ validity of assumptions \ref{A8*} and \ref{A9*} needs to be shown, such that Theorems \ref{theobootstrap} and \ref{theobootstrapalternative} apply and Algorithm \ref{bootstrapalg} gives valid approximation of the critical value. Denote the conditional density of $\varepsilon({\theta_0})$ (defined in Algorithm \ref{bootstrapalg}) given $X$ by $f_{\varepsilon({\theta_0})|X}$. Then we need the following assumptions. 
\begin{enumerate}[label=(\textbf{A10})]
	\item\label{A10} Denote the conditional density of $\varepsilon({\theta_0})$ as defined in Algorithm \ref{bootstrapalg} given $X$ by $f_{\varepsilon({\theta_0})|X}$. 	Let $\mathcal{K}\subseteq\mathbb{R}$ be compact and $\ell $ be bounded and $r$-times continuously differentiable with bounded derivatives and denote the $k$-th derivative of $\ell $ by $\ell ^{(k)}$. Further, assume
	\begin{equation}\label{expectlambda}
	\underset{u\in\mathcal{K}}{\sup}\,E\bigg[\int ||g_{\t_0}(X)+u-a_ne||^l|\ell ^{(j)}(e)|f_{\varepsilon({\theta_0})|X}(u-a_ne|X)\,de\bigg]<C
	\end{equation}
	for $l\in\{0,j\}$ and
	\begin{equation}\label{expectdevlambda}
	\underset{u\in\mathcal{K}}{\sup}\,E\bigg[\int \big|\big|\dot{\Lambda}_{\t_0}(\Lambda_{\t_0}^{-1}(g_{\t_0}(X)+u-a_ne))\big|\big|^j|\ell ^{(j)}(e)|f_{\varepsilon({\theta_0})|X}(u-a_ne|X)\,de\bigg]<C
	\end{equation}
	as well as
	\begin{equation}\label{orderHess}
	\frac{||\hat{\t}-\t_0||^j}{na_n}\,\underset{u\in\mathcal{K}}{\sup}\,\sum_{i=1}^n\bigg|\ell ^{(j)}\bigg(\frac{u-\varepsilon_i(\theta_0)}{a_n}\bigg)\bigg|\underset{||\t-\t_0||<\delta}{\sup}\,||Hess\Lambda_{\t}(Y_i)||^j=\mathcal{O}_p(1)
	\end{equation}
	for sufficiently large $C>0,n\in\mathbb{N}$ and all $j=1,...,r-1$. Moreover, let
	\begin{equation}\label{expectLambdadevLambda2}
	E\big[|\Lambda_{\t_0}(Y)|^r+||\dot{\Lambda}_{\t_0}(Y)||^r\big]<\infty
	\end{equation}
	and
	\begin{equation}\label{orderHess2}
	\frac{||\tilde{\t}-\t_0||^r}{n}\sum_{i=1}^n\underset{||\t-\t_0||<\delta}{\sup}\,||Hess\Lambda_{\t}(Y_i)||^r=\mathcal{O}_p(1).
	\end{equation}
\end{enumerate}
\begin{enumerate}[label=(\textbf{B11*})]
	\item\label{B11*} In the following, the notations from Algorithm \ref{bootstrapalg} are employed.
	Let $1/(na_n)\to 0$, $(\log n)/(nb_n^{d_X+4})\to 0$.
	\begin{enumerate}[label=(\textbf{\arabic{*}})]
		\item Let $\hat{g}$ be $(q+2)$-times continuously differentiable (same $q$ as in \ref{B4}).
		\item Assume 
		\begin{enumerate}[label=$\bullet$]
			\item $\kappa $ and $\ell $ are $(q+2)$-times continuously differentiable,
			\item $\kappa $ has bounded support,
			\item $\kappa (0)>0$ and $\ell >0$,
			\item it is either $\big|\frac{\partial}{\partial u}\ell (u)\big|\leq K$,
			\begin{equation}\label{boundedfxi}
			|\tilde{f}(u)|<|u|^{-\nu}\quad\textup{and}\quad\Big|\frac{\partial}{\partial u}\tilde{f}(u)\Big|\leq K|u|^{-\nu}
			\end{equation}
			for all $|u|>L$ and $\tilde{f}\in\Big\{\ell ,\frac{\partial}{\partial u}\ell \Big\}$ or
			$$\Big|\Big|\frac{\partial}{\partial x}\tilde{f}(x)\Big|\Big|\leq K\quad\textup{and}\quad\Big|\Big|\frac{\partial}{\partial x}\tilde{f}(x)\Big|\Big|I_{\{||x||>L\}}\leq K||x||^{-\nu}I_{\{||x||>L\}}$$
			for some $\nu>1,K,L\in(0,\infty)$ and all $\tilde{f}\in\Big\{\kappa ,\frac{\partial}{\partial x_i}\kappa ,\frac{\partial^2}{\partial x_i^2}\kappa \Big\}$, where the same $i$ was used as in \ref{B3}. From now on, the case $i=1$ is considered w.l.o.g.
		\end{enumerate}
	\end{enumerate}
\end{enumerate}

\begin{lemma}\label{lemmaan}
	Let $\mathcal{K}\subset\mathbb{R}$ be compact. Assume \ref{A10} and
	\begin{equation}\label{convghat}
	\frac{\underset{i=1,...,n}{\max}\,|\hat{g}(X_i)-g_{\t_0}(X_i)|^r}{a_n^{r+1}}=o_P(1)\quad\textup{and}\quad\frac{||\hat{\t}-\t_0||^r}{a_n^{r+1}}=o_P(1).
	\end{equation}
	Then,
	$$\underset{u\in\mathcal{K}}{\sup}\,\bigg|\frac{1}{n}\sum_{i=1}^nF_{\xi}\bigg(\frac{u-\hat{\varepsilon}_i}{a_n}\bigg)-\frac{1}{n}\sum_{i=1}^nF_{\xi}\bigg(\frac{u-\varepsilon_i}{a_n}\bigg)\bigg|=o_P(1)$$
	and
	$$\underset{u\in\mathcal{K}}{\sup}\,\bigg|\frac{1}{na_n}\sum_{i=1}^n\ell \bigg(\frac{u-\hat{\varepsilon}_i}{a_n}\bigg)-\frac{1}{na_n}\sum_{i=1}^n\ell \bigg(\frac{u-\varepsilon_i}{a_n}\bigg)\bigg|=o_P(1).$$
\end{lemma}

\begin{rem}
	\begin{enumerate}[label=(\roman*)]
		\item Later, it will be shown that $\frac{1}{n}\sum_{l=1}^n\hat{\varepsilon}_l=o_P(1)$ in the proof of Lemma \ref{lemmahansenbootstrap}. Hence, the Assertion of Lemma \ref{lemmaan} can be extended to the centred residuals $\tilde{\varepsilon}_i$ by considering a slightly larger set $\mathcal{K}$.
		\item Many of the assumptions in \ref{A10} can be replaced by less complex, but more restrictive versions. For example, due to (\ref{convghat}) assumption (\ref{orderHess}) is implied by
		$$E\bigg[\underset{||\t-\t_0||<\delta}{\sup}\,||Hess\Lambda_{\t}(Y_i)||^j\bigg]<\infty$$
		or
		$$\underset{k=1,...,n}{\max}\underset{||\t-\t_0||<\delta}{\sup}\,||Hess\Lambda_{\t}(Y_i)||^j=O_P\big(||\hat{\t}-\t_0||^{-j}\big)\quad\Big(=O_P\big(n^{\frac{j}{4}}\big)\quad\textup{under the alternative}\Big)$$
		for all $j=1,...,r-1.$
		\item It is $||\hat{\t}-\t_0||=O_P\big(n^{-\frac{1}{4}}\big)$ under the alternative. Unfortunately, the experimenter in advance does not know, if the null hypothesis or the alternative holds, so that in general (\ref{convghat}) limits $a_n$ to $a_n^{-1}=o\big(n^{\frac{r}{4(r+1)}}\big)=o\big(n^{\frac{1}{4}}\big)$.
		\item $g_{\t_0}$ can be estimated using the Nadaraya-Watson approach as in \cite{HSvK2015}. Under some additional assumptions, their Proposition 6.1 or to be precise its extension in the supplementary material of \citet[p.~7]{CvK2016} yields $\hat{g}(x)-g_{\t_0}(x)=O_P(((\log n)/(nh_x^{d_X}))^{1/2}+||\hat{\t}-\t_0||)$ uniformly on compact sets. When assuming the existence of some compact set $\mathcal{K}\subseteq\operatorname{supp}(f_X)$ such that $\operatorname{supp}(v)$ is contained in the interior of $\mathcal{K}$, Equation (\ref{convghat}) (and a counterpart of (\ref{unifconvghat}) below for compact sets) can be obtained when discarding those $(Y_i,X_i),(Y_j^*,X_j^*)$ such that $X_i,X_j^*\notin\mathcal{K}$ (note that an equation similar to (\ref{exprh}) can still be derived (although in general with another $\psi$)). Then, (\ref{convghat}) requires 
		$a_n^{-1}=o(((\log n)/(nh_x^{d_X}))^{-\frac{r}{2(r+1)}}+n^{\frac{r}{4(r+1)}})$ for $||\hat{\t}-\t_0||=O_P\big(n^{-\frac{1}{4}}\big)$.
	\end{enumerate}
\end{rem}

\noindent {\bf Proof of Lemma \ref{lemmaan}.}
Only the second assertion is shown since the first one can be concluded similarly. The proof uses similar techniques as \citet{Han2008}. First, for the deviation terms $R_{i}=\hat{\varepsilon}_i-\varepsilon_i(\t_0)$ and appropriate $u_i^*,i=1,...,n,$ a Taylor expansion leads to
\begin{equation}\label{Taylor}
\frac{1}{na_n}\sum_{i=1}^n\ell \bigg(\frac{u-\hat\varepsilon_i}{a_n}\bigg)=\frac{1}{na_n}\sum_{i=1}^n\Bigg(\sum_{j=0}^{r-1}\frac{(-R_i)^j}{a_n^jj!}\ell ^{(j)}\bigg(\frac{u-\varepsilon_i(\theta_0)}{a_n}\bigg)+\frac{(-R_i)^r}{a_n^{r}r!}\ell ^{(r)}(u_i^*)\Bigg)
\end{equation}
For appropriate $\tilde{\theta}_y$ between $\hat{\theta}$ and $\theta_0$ the $R_i$ can be split into
\begin{align}
R_i&=\hat{\varepsilon}_i-\varepsilon_i(\theta_0)\nonumber
\\[0,2cm]&=\frac{\Lambda_{\hat{\t}}(Y_i)-\Lambda_{\hat{\t}}(0)}{\Lambda_{\hat{\t}}(1)-\Lambda_{\hat{\t}}(0)}-\hat{g}(X_i)-\varepsilon_i(\theta_0)\nonumber
\\[0,2cm]&=\frac{\Lambda_{\hat{\t}}(Y_i)-\Lambda_{\hat{\t}}(0)}{\Lambda_{\hat{\t}}(1)-\Lambda_{\hat{\t}}(0)}-\frac{\Lambda_{\t_0}(Y_i)-\Lambda_{\t_0}(0)}{\Lambda_{\t_0}(1)-\Lambda_{\t_0}(0)}+g_{\t_0}(X_i)-\hat{g}(X_i)\nonumber
\\[0,2cm]&=\frac{1}{(\Lambda_{\hat{\t}}(1)-\Lambda_{\hat{\t}}(0))(\Lambda_{\t_0}(1)-\Lambda_{\t_0}(0))}(\Lambda_{\hat{\t}}(Y_i)-\Lambda_{\hat{\t}}(0))(\Lambda_{\t_0}(1)-\Lambda_{\hat{\t}}(1)+\Lambda_{\hat{\t}}(0)-\Lambda_{\t_0}(0))\nonumber
\\[0,2cm]&\quad+\frac{\Lambda_{\hat{\t}}(Y_i)-\Lambda_{\t_0}(Y_i)+\Lambda_{\t_0}(0)-\Lambda_{\hat{\t}}(0)}{\Lambda_{\t_0}(1)-\Lambda_{\t_0}(0)}+g_{\t_0}(X_i)-\hat{g}(X_i)\nonumber
\\[0,2cm]&=\frac{1}{(\Lambda_{\hat{\t}}(1)-\Lambda_{\hat{\t}}(0))(\Lambda_{\t_0}(1)-\Lambda_{\t_0}(0))}\Big(\Lambda_{\t_0}(Y_i)-\Lambda_{\t_0}(0)+\big(\dot{\Lambda}_{\t_0}(Y_i)-\dot{\Lambda}_{\t_0}(0)\big)(\hat{\t}-\t_0)\nonumber
\\[0,2cm]&\quad\quad+\frac{1}{2}(\hat{\t}-\t_0)^t\big(Hess\Lambda_{\tilde{\t}_{Y_i,0}}(Y_i)-Hess\Lambda_{\tilde{\t}_{Y_i,0}}(0)\big)(\hat{\t}-\t_0)\Big)\nonumber
\\[0,2cm]&\quad\quad\Big(-\dot{\Lambda}_{\t_0}(1)(\hat{\t}-\t_0)-\frac{1}{2}(\hat{\t}-\t_0)^tHess\Lambda_{\tilde{\t}_1}(1)(\hat{\t}-\t_0)\nonumber
\\[0,2cm]&\quad\quad+\dot{\Lambda}_{\t_0}(0)(\hat{\t}-\t_0)+\frac{1}{2}(\hat{\t}-\t_0)^tHess\Lambda_{\tilde{\t}_0}(0)(\hat{\t}-\t_0)\Big)\nonumber
\\[0,2cm]&\quad+\frac{1}{\Lambda_{\t_0}(1)-\Lambda_{\t_0}(0)}\Big(\dot{\Lambda}_{\t_0}(Y_i)(\hat{\t}-\t_0)+\frac{1}{2}(\hat{\t}-\t_0)^tHess\Lambda_{\tilde{\t}_{Y_i}}(Y_i)(\hat{\t}-\t_0)\nonumber
\\[0,2cm]&\quad\quad-\dot{\Lambda}_{\t_0}(0)(\hat{\t}-\t_0)-\frac{1}{2}(\hat{\t}-\t_0)^tHess\Lambda_{\tilde{\t}_0}(0)(\hat{\t}-\t_0)\Big)+g_{\t_0}(X_i)-\hat{g}(X_i)\nonumber
\\[0,2cm]&=\tilde{R}_i+g_{\t_0}(X_i)-\hat{g}(X_i).\label{bootsepsilon}
\end{align}
Therefore,
\begin{align*}
\frac{1}{na_n}\sum_{i=1}^n\bigg|\frac{R_i^j}{a_n^{j}j!}\ell ^{(j)}\bigg(\frac{u-\varepsilon_i(\theta_0)}{a_n}\bigg)\bigg|&\leq\frac{C}{na_n}\sum_{i=1}^n\bigg|\frac{\tilde{R}_i^j}{a_n^{j}j!}\ell ^{(j)}\bigg(\frac{u-\varepsilon_i(\theta_0)}{a_n}\bigg)\bigg|
\\[0,2cm]&\quad+\frac{C}{na_n}\sum_{i=1}^n\bigg|\frac{(\hat{g}(X_i)-g_{\t_0}(X_i))^j}{a_n^{j}j!}\ell ^{(j)}\bigg(\frac{u-\varepsilon_i(\theta_0)}{a_n}\bigg)\bigg|
\end{align*}
for all $j=1,...,r-1$ and
$$\frac{1}{na_n}\sum_{i=1}^n\bigg|\frac{R_i^r}{a_n^{r}r!}\ell ^{(r)}(u_i^*)\bigg|\leq\frac{C}{na_n}\sum_{i=1}^n\bigg|\frac{\tilde{R}_i^r}{a_n^{r}r!}\ell ^{(r)}(u_i^*)\bigg|+\frac{C}{na_n}\sum_{i=1}^n\bigg|\frac{(\hat{g}(X_i)-g_{\t_0}(X_i))^r}{a_n^{r}r!}\ell ^{(r)}(u_i^*)\bigg|$$
for some sufficiently large constant $C>0$, so that it suffices to treat the cases $R_i^{(1)}=\hat{g}(X_i)-g_{\t_0}(X_i)$ and $R_i^{(2)}=\tilde{R}_i$ separately.\\
When inserting $R_i^{(1)}$ in equation (\ref{Taylor}) negligibility of the last summand directly follows from (\ref{convghat}) and the boundedness of $\ell ^{(r)}$. Thanks to \citet{Han2008}, to prove
$$\underset{u\in\mathcal{K}}{\sup}\,\frac{1}{na_n}\sum_{i=1}^n\bigg|\ell ^{(j)}\bigg(\frac{u-\varepsilon_i(\theta_0)}{a_n}\bigg)\bigg|\leq C+o_P(1)$$
for all $j=0,...,r-1$ and some constant $C>0$, it suffices to show uniform (with respect to $u$) boundedness of the expectation. Hence one has
\begin{align*}
E\bigg[\frac{1}{na_n}\sum_{i=1}^n\bigg|\ell ^{(j)}\bigg(\frac{u-\varepsilon_i(\theta_0)}{a_n}\bigg)\bigg|\bigg]&=E\bigg[\frac{1}{a_n}\bigg|\ell ^{(j)}\bigg(\frac{u-\varepsilon(\theta_0)}{a_n}\bigg)\bigg|\bigg]
\\[0,2cm]&=\int\frac{1}{a_n}\bigg|\ell ^{(j)}\bigg(\frac{u-e}{a_n}\bigg)\bigg|f_{\varepsilon(\t_0)}(e)\,de
\\[0,2cm]&=\int|\ell ^{(j)}(e)|f_{\varepsilon(\t_0)}(u-a_ne)\,de
\\[0,2cm]&\leq C
\end{align*}
for some constant $C>0$ (see (\ref{expectlambda})) and thus
\begin{align*}
&\frac{1}{na_n}\sum_{i=1}^n\bigg|\frac{(\hat{g}(X_i)-g_{\t_0}(X_i))^j}{a_n^jj!}\ell ^{(j)}\bigg(\frac{u-\varepsilon_i(\theta_0)}{a_n}\bigg)\bigg|
\\[0,2cm]&\leq\underbrace{\frac{\underset{k=1,...,n}{\max}\,|\hat{g}(X_k)-g_{\t_0}(X_k)|^j}{a_n^jj!}}_{=o_P(1)}\underbrace{\frac{1}{na_n}\sum_{i=1}^n\bigg|\ell ^{(j)}\bigg(\frac{u-\varepsilon_i(\theta_0)}{a_n}\bigg)\bigg|}_{=O_P(1)}
\\[0,2cm]&=o_P(1)
\end{align*}
for all $j=1,...,r-1$. Further $\tilde{R}_i$ can be written as
\begin{align*}
\tilde{R}_i&=O_P(||\hat{\t}-\t_0||)+O_P(||\hat{\t}-\t_0||)\Big(\Lambda_{\t_0}(Y_i)+\dot{\Lambda}_{\t_0}(Y_i)(\hat{\t}-\t_0)
\\[0,2cm]&\quad\quad+\frac{1}{2}(\hat{\t}-\t_0)^t\big(Hess\Lambda_{\tilde{\t}_{Y_i,0}}(Y_i)-Hess\Lambda_{\tilde{\t}_{Y_i,0}}(0)\big)(\hat{\t}-\t_0)\Big)
\\[0,2cm]&\quad+\frac{1}{\Lambda_{\t_0}(1)-\Lambda_{\t_0}(0)}\Big(\dot{\Lambda}_{\t_0}(Y_i)(\hat{\t}-\t_0)+\frac{1}{2}(\hat{\t}-\t_0)^tHess\Lambda_{\tilde{\t}_{Y_i}}(Y_i)(\hat{\t}-\t_0)\Big),
\end{align*}
where the $O_P$-terms are independent of $i$. When inserting $\tilde{R}_i$ in equation (\ref{Taylor}), one has for any $\delta>0$
\begin{align*}
&\frac{1}{na_n}\sum_{i=1}^n\bigg|\frac{\tilde{R}_i^j}{a_n^jj!}\ell ^{(j)}\bigg(\frac{u-\varepsilon_i(\theta_0)}{a_n}\bigg)\bigg|
\\[0,2cm]&\leq O_P\bigg(\frac{||\hat{\t}-\t_0||^j}{a_n^j}\frac{1}{na_n}\sum_{i=1}^n\bigg|\ell ^{(j)}\bigg(\frac{u-\varepsilon_i(\theta_0)}{a_n}\bigg)\bigg|\Big(1+|\Lambda_{\t_0}(Y_i)|^j+||\dot{\Lambda}_{\t_0}(Y_i)||^j
\\[0,2cm]&\quad+||\hat{\t}-\t_0||^j\underset{||\t-\t_0||<\delta}{\sup}\,||Hess\Lambda_{\t}(Y_i)||^j\Big)\bigg)
\end{align*}
for all $j=1,...,r-1$. By assumptions (\ref{expectlambda}) and (\ref{expectdevlambda}) the expected value of the sum
$$\frac{1}{na_n}\sum_{i=1}^n\bigg|\ell ^{(j)}\bigg(\frac{u-\varepsilon_i(\theta_0)}{a_n}\bigg)\bigg|\Big(1+|\Lambda_{\t_0}(Y_i)|^j+||\dot{\Lambda}_{\t_0}(Y_i)||^j+||\tilde{\t}-\t_0||^j\underset{||\t-\t_0||<\delta}{\sup}\,||Hess\Lambda_{\t}(Y_i)||^j\Big)$$
can be bounded by some constant $C>0$, so that
\begin{align*}
&\frac{1}{na_n}\sum_{i=1}^n\bigg|\frac{\tilde{R}_i^j}{a_n^jj!}\ell ^{(j)}\bigg(\frac{u-\varepsilon_i(\theta_0)}{a_n}\bigg)\bigg|
\\[0,2cm]&\leq O_P\bigg(\frac{||\hat{\t}-\t_0||^{2j}}{a_n^j}\frac{1}{na_n}\sum_{i=1}^n\bigg|\ell ^{(j)}\bigg(\frac{u-\varepsilon_i(\theta_0)}{a_n}\bigg)\bigg|\underset{||\t-\t_0||<\delta}{\sup}\,||Hess\Lambda_{\t}(Y_i)||^j\bigg)+o_P(1)
\\[0,2cm]&=o_P(1)
\end{align*}
by (\ref{orderHess}) and (\ref{convghat}). The remaining term can be treated similarly by applying (\ref{expectLambdadevLambda2}) and (\ref{orderHess2}) to obtain
\begin{align*}
&\frac{1}{na_n}\sum_{i=1}^n\bigg|\frac{\tilde{R}_i^r}{a_n^rr!}\ell ^{(r)}(u_i^*)\bigg|
\\[0,2cm]&\leq O_P\bigg(\frac{||\hat{\t}-\t_0||^r}{a_n^{r+1}}\frac{1}{n}\sum_{i=1}^n\Big(|\Lambda_{\t_0}(Y_i)|^r+||\dot{\Lambda}_{\t_0}(Y_i)||^r+||\hat{\t}-\t_0||^r\underset{||\t-\t_0||<\delta}{\sup}\,||Hess\Lambda_{\t}(Y_i)||^r\Big)\bigg)
\\[0,2cm]&=o_P(1).
\end{align*}
Altogether one obtains 
$$\frac{1}{na_n}\sum_{i=1}^n\ell \bigg(\frac{u-\hat\varepsilon_i}{a_n}\bigg)=\frac{1}{na_n}\sum_{i=1}^n\ell \bigg(\frac{u-\varepsilon_i(\theta_0)}{a_n}\bigg)+o_P(1)$$
uniformly on compact sets.
\hfill $\Box$

\bigskip

\begin{lemma}\label{lemmabootstrap}
	Let Assumptions \ref{B11*} and (\ref{convghat}) be fulfilled. Further, assume \ref{A1}--\ref{A8},\ref{A10},\ref{B1}--\ref{B10},
	\begin{equation}\label{unifconvghat}
	\underset{x\in\mathbb{R}^{d_X}}{\sup}\,|\hat{g}(x)-g_{\t_0}(x)|=o_P(1)\quad\textup{and}\quad\underset{x\in\mathbb{R}^{d_X}}{\sup}\,\bigg|\frac{\partial}{\partial x_1}\hat{g}(x)-\frac{\partial}{\partial x_1}g_{\t_0}(x)\bigg|=o_P(1).
	\end{equation}
	Further, assume the existence of a neighbourhood $\tilde{\Theta}$ of $\theta_0$ such that the map $y\mapsto\Lambda_{\theta}(y)$ is $(q+2)$-times continuously differentiable for all $\theta\in\tilde{\Theta}$. Let $\hat{h}^*$ be the estimator from (\ref{hath}) based on the bootstrap data $(Y_j^*,X_j^*),j=1,...,m$. Assume that the density of $\varepsilon(\theta_0)$ is continuous and
	\begin{equation}\label{boundedepsfeps}
	\underset{e\in\mathbb{R}}{\sup}\,|ef_{\varepsilon(\t_0)}(e)|<\infty.
	\end{equation}
	Then, assumptions \ref{A8*} and \ref{A9*} are fulfilled. 
\end{lemma}

\noindent {\bf Proof of Lemma \ref{lemmabootstrap}.}
Note that conditional on $(Y_1,X_1),...,(Y_n,X_n)$ the random variables $(Y_1^*,X_1^*),...,(Y_m^*,X_m^*)$ are independent as well as identically distributed. Moreover, after conditioning on the original data, the Assumptions \ref{B1}--\ref{B10} are valid for the bootstrap sample with probability converging to one, so that due to Remark \ref{rembootsassumpCKC} the same reasoning as in \citep{CvK2018} can be applied to obtain (\ref{exprh*}).\\
For notational convenience the conditional distribution of $(Y_1^*,X_1^*),...,(Y_m^*,X_m^*)$ conditional on $(Y_1,X_1),...,(Y_n,X_n)$ is written as $P^*$ and the expectation with respect to $P^*$ is written as $E^*$. 
Let $F_{Y^*|X^*}$ denote the conditional distribution function of $Y_1^*$ conditioned on $X_1^*$ (and $(Y_1,X_1),...,(Y_n,Y_n)$). To verify \ref{A8*} $\psi^*$ has to be examined further and to define $\psi^*$ some further notations are needed. Let $v$ be the weighting function from assumption \ref{B7} and define
$$s_1^*(u,x)=\int_0^u\frac{\frac{\partial F_{Y^*|X^*}(y|x)}{\partial y}}{\frac{\partial F_{Y^*|X^*}(y|x)}{\partial x_1}}\,dy,\quad\tilde{v}_1^*(u_0,x)=\frac{v(x)}{s_1^*(u_0,x)},\quad\tilde{v}^*_2(u_0,x)=\frac{v(x)s_1^*(u_0,x)}{s_1^*(1,x)^2}$$
and (for $\tilde{v}^*=\tilde{v}_1^*,\tilde{v}_2^*$)
\begin{align*}
{\delta_j^*}^{\tilde{v}^*}(u_0,u)&=\int_{\max(0,U_j^*)}^{\max(u,U_j^*)}\bigg(\tilde{v}^*(u_0,X_j^*)D_{p,0}^*(r,X_j^*)-\frac{\partial}{\partial x_1}\big(\tilde{v}^*(u_0,x)D_{p,1}^*(r,x)\big)\Big|_{x=X_j^*}\bigg)\,dr
\\[0,2cm]&\quad+\int_0^u\bigg(\tilde{v}^*(u_0,X_j^*)D_{f,0}^*(r,X_j^*)-\frac{\partial}{\partial x_1}\big(\tilde{v}^*(u_0,X_j^*)D_{f,1}^*(r,x)\big)\Big|_{x=X_j^*}\bigg)\,dr
\\[0,2cm]&\quad+(\mathds{1}_{\{U_j^*\leq u\}}-\mathds{1}_{\{U_j^*\leq0\}})\tilde{v}^*(u_0,X_j^*)D_{p,u}^*(U_j^*,X_j^*)
\\[0,2cm]&\quad+\int_0^u\bigg(\frac{\mathds{1}_{\{U_j^*\leq u\}}-\mathds{1}_{\{U_j^*\leq0\}}}{F_{U^*}(1)-F_{U^*}(0)}-r\bigg)
\\[0,2cm]&\quad\quad\int_{\mathcal{X}}\bigg(\bigg(\tilde{v}^*(u_0,x)D_{p,0}^*(r,x)+\frac{\partial}{\partial x_1}\big(\tilde{v}^*(u_0,x)D_{p,1}^*(r,x)\big)\bigg)
\\[0,2cm]&\quad\quad f_{U^*,X^*}(r,x)+\tilde{v}^*(u_0,x)D_{p,u}^*(r,x)\frac{\partial}{\partial r}f_{U^*,X^*}(r,x)\bigg)\,dx\,dr
\\[0,2cm]&\quad-\bigg(\frac{\mathds{1}_{\{U_j^*\leq 1\}}-\mathds{1}_{\{U_j^*\leq0\}}}{F_{U^*}(1)-F_{U^*}(0)}-1\bigg)\int_0^ur\int_{\mathcal{X}}\bigg(\tilde{v}^*(u_0,x)D_{p,0}^*(r,x)
\\[0,2cm]&\quad\quad-\tilde{v}^*(u_0,x)\frac{\partial}{\partial r}D_{p,u}^*(r,x)+\frac{\partial}{\partial x_1}\big(\tilde{v}^*(u_0,x)D_{p,1}^*(r,x)\big)\bigg)f_{U^*,X^*}(r,x)\,dx\,dr
\\[0,2cm]&\quad\quad-\bigg(\frac{\hat{F}_{U^*}(1)-\hat{F}_{U^*}(0)}{F_{U^*}(1)-F_{U^*}(0)}-1\bigg)u\int_{\mathcal{X}}\tilde{v}^*(u_0,x)D^*_{p,u}(u,x)f_{U^*,X^*}(u,x)\,dx,
\end{align*}
where $D_{p,0}^*(u,x),...,D_{f,1}^*(u,x)$ are defined as
\begin{align}
D_{p,0}^*(u,x)&=\frac{\Phi^*(u,x)\frac{\partial}{\partial x_1}f_{X^*}(x)}{\Phi_i^*(u,x)^2f_{X^*}(x)^2},\label{Dp0*}
\\[0,2cm]D_{p,u}^*(u,x)&=\frac{1}{f_{X^*}(x)\Phi_i^*(u,x)},\nonumber
\\[0,2cm]D_{p,1}^*(u,x)&=\frac{-\Phi_u^*(u,x)}{f_{X^*}(x)\Phi_i^*(u,x)^2}\nonumber
\\[0,2cm]D_{f,0}^*(u,x)&=\frac{-\Phi_u^*(u,x)\Phi^*(u,x)\frac{\partial}{\partial x_1}f_{X^*}(x)}{\Phi_i^*(u,x)^2f_{X^*}(x)^2}\nonumber
\\[0,2cm]D_{f,1}^*(u,x)&=\frac{\Phi_u^*(u,x)\Phi^*(u,x)}{\Phi_i^*(u,x)^2f_{X^*}(x)}\label{Dfi*}
\end{align}
with $\Phi^*(u|x)=F_{U^*|X^*}(u,x)=\frac{p^*(u,x)}{f_{X^*}(x)},\Phi^*_u(u|x)=\frac{\partial}{\partial _u}\Phi^*(u,x)$, $\Phi^*_i(u|x)=\frac{\partial}{\partial x_1}\Phi^*(u,x)$ and $f_{X^*}$ from Algorithm \ref{bootstrapalg}. Then, $\psi^*$ is defined as
\begin{align}
\psi^*(Z_j^*,u)&={\delta_j^*}^{\tilde{v}_1^*}(1,u)-{\delta_j^*}^{\tilde{v}_2^*}(u,1)+\frac{{Q^*}'(u)}{F_{U^*}(1)-F_{U^*}(0)}\big(\mathds{1}_{\{U_j^*\leq u\}}-\mathds{1}_{\{U_j^*\leq 0\}}-F_{U^*}(u)+F_{U^*}(0)\big)\nonumber
\\[0,2cm]&\quad-{Q^*}'(u)\frac{F_{U^*}(u)-F_{U^*}(0)}{(F_{U^*}(1)-F_{U^*}(0))^2}\big(\mathds{1}_{\{U_j^*\leq1\}}-\mathds{1}_{\{U_j^*\leq 0\}}-F_{U^*}(1)+F_{U^*}(0)\big).\label{defpsi^*}
\end{align}
Condition \ref{A8*} for $\psi^*$ is implied by the same reasoning as in \cite{CvK2018}. Note that the first part of Remark \ref{rembootsassumpCKC} ensures that $v$ can be used as the weighting function for the bootstrap data as well.\\
To prove \ref{A9*}, an auxiliary lemma is shown in the following. Thanks to the expressions above for $\psi^*,D_{p,0}^*(u,x),...,D_{f,1}^*(u,x)$, equations (\ref{univconvdevFZ*}) and (\ref{convpsi^*}) will be a direct consequence of Lemma \ref{lemmahansenbootstrap}, while (\ref{boundedpsi^*}) follows from expression (\ref{defpsi^*}), so that $\psi^*$ fulfils \ref{A9*} then.
\begin{lemma}\label{lemmahansenbootstrap}
	Let $\mathcal{C}\subseteq\big(-\frac{F_Y(0)}{F_Y(1)-F_Y(0)},\frac{1-F_Y(0)}{F_Y(1)-F_Y(0)}\big)$ be a compact set and define
	$$F_{U^*,X^*}(u,x)=P^*(U^*\leq u,X^*\leq x),\quad\textup{and}\quad p^*(u,x)=\frac{\partial^{d_X}}{\partial x_1...\partial x_{d_X}}F_{U^*,X^*}(u,x).$$
	Under the assumptions of Lemma \ref{lemmabootstrap}, one has
	\begin{align}
	\underset{x\in\operatorname{supp}(v)}{\sup}\,|f_{X^*}(x)-{f}(x)|&=o_P(1),\nonumber
	\\[0,2cm]\underset{x\in\operatorname{supp}(v)}{\sup}\,\bigg|\frac{\partial}{\partial x_1}f_{X^*}(x)-\frac{\partial}{\partial x_1}{f}(x)\bigg|&=o_P(1),\nonumber
	\\[0,2cm]\underset{x\in\operatorname{supp}(v)}{\sup}\,\bigg|\frac{\partial^2}{\partial x_1^2}f_{X^*}(x)-\frac{\partial^2}{\partial x_1^2}{f}(x)\bigg|&=o_P(1),\nonumber
	\\[0,2cm]\underset{u\in\mathcal{C},x\in\operatorname{supp}(v)}{\sup}\,|p^*(u,x)-{p}(u,x)|&=o_P(1),\nonumber
	\\[0,2cm]\underset{u\in\mathcal{C},x\in\operatorname{supp}(v)}{\sup}\,\bigg|\frac{\partial}{\partial u}p^*(u,x)-\frac{\partial}{\partial u}{p}(u,x)\bigg|&=o_P(1),\nonumber
	\\[0,2cm]\underset{u\in\mathcal{C},x\in\operatorname{supp}(v)}{\sup}\,\bigg|\frac{\partial}{\partial x_1}p^*(u,x)-\frac{\partial}{\partial x_1}{p}(u,x)\bigg|&=o_P(1),\nonumber
	\\[0,2cm]\underset{u\in\mathcal{C},x\in\operatorname{supp}(v)}{\sup}\,\bigg|\frac{\partial^2}{\partial u\partial x_1}p^*(u,x)-\frac{\partial^2}{\partial u\partial x_1}{p}(u,x)\bigg|&=o_P(1),\nonumber
	\\[0,2cm]\underset{u\in\mathcal{C}}{\sup}\,|F_{U^*}(u)-F_{{U}}(u)|&=o_P(1),\nonumber
	\\[0,2cm]\underset{z\in\mathcal{R}^{d_X+1}}{\sup}\,|F_{Z^*}(z)-F_{{Z}}(z)|&=o_P(1)\nonumber.
	\end{align}
	Here, $p$ is defined as in \ref{defpf} and $f_{X^*}$ is defined as in Algorithm \ref{bootstrapalg}.
\end{lemma}
\textbf{Proof:} Most of the proof contains in applying the results in \cite{Han2008} for kernel estimates. While doing so, note that due to (\ref{boundedfxi}) and (\ref{boundedepsfeps}) kernel estimates like
$$\hat{f}_{\varepsilon(\t_0)}(e)=\frac{1}{na_n}\sum_{i=1}^n\ell \bigg(\frac{e-\varepsilon(\t_0)}{a_n}\bigg)$$
converge uniformly in $e\in\mathbb{R}$ to their expectation (see Theorem 4 in \cite{Han2008}).\\
The results for $f_{X^*}$ directly follow from \citet{Han2008} (note that $\kappa $ is a kernel of order $2$). The assertion for $\frac{\partial}{\partial x_1}f_{X^*}$ and $\frac{\partial^2}{\partial x_1^2}f_{X^*}$ follows similarly by applying for example $E\big[\frac{\partial^2}{\partial v_i^2}f_{X^*}(v)\big]=\frac{\partial^2}{\partial x_1^2}f_X(x)+\mathcal{O}(b_n^2).$ Moreover, $p^*$ can be expressed for any $j\in\{1,...,m\}$ as
\begin{align*}
&p^*(\mathcal{T}^*((h^*)^{-1}(u)),x)
\\[0,2cm]&=\frac{\partial^{d_X}}{\partial x_1...\partial x_{d_X}}F_{U^*,X^*}(\mathcal{T}^*((h^*)^{-1}(u)),x)
\\[0,2cm]&=\frac{\partial^{d_X}}{\partial x_1...\partial x_{d_X}} P^*(S^*_j\leq u,X_j^*\leq x)
\\[0,2cm]&=\frac{\partial^{d_X}}{\partial x_1...\partial x_{d_X}}\int_{(-\infty,x]}P^*(S^*_j\leq u|X^*_j=z)f_{X^*}(z)\,dz
\\[0,2cm]&=\frac{\partial^{d_X}}{\partial x_1...\partial x_{d_X}}\int_{(-\infty,x]}\frac{1}{n}\sum_{k=1}^nF_{\xi}\bigg(\frac{u-\hat{g}(z)-\hat{\varepsilon}_k+\frac{1}{n}\sum_{l=1}^n\hat{\varepsilon}_l}{a_n}\bigg)f_{X^*}(z)\,dz
\\[0,2cm]&=\frac{1}{n}\sum_{k=1}^nF_{\xi}\bigg(\frac{u-\hat{g}(x)-\hat{\varepsilon}_k+\frac{1}{n}\sum_{l=1}^n\hat{\varepsilon}_l}{a_n}\bigg)f_{X^*}(x),
\end{align*}
where $(-\infty,x]=\times_{i=1}^{d_X}(-\infty,x_i]$ and $F_{\xi}$ denotes the cumulative distribution function corresponding to $\ell $. From now on, only $\frac{\partial}{\partial x_1}p^*$ is considered, since the other terms can be treated analogously. Due to (\ref{bootsepsilon}), one obtains (if $||\hat{\t}-\t_0||\leq\delta$)
\begin{align*}
\bigg|\frac{1}{n}\sum_{l=1}^n\hat{\varepsilon}_l\bigg|&\leq\bigg|\frac{1}{n}\sum_{l=1}^n\varepsilon_l(\theta_0)\bigg|+O_P(||\hat{\t}-\t_0||)\frac{1}{n}\sum_{l=1}^n\Big(1+|\Lambda_{\t_0}(Y_l)|+||\dot{\Lambda}_{\t_0}(Y_l)||
\\[0,2cm]&\quad\quad+||\hat{\t}-\t_0||\underset{||\t-\t_0||<\delta}{\sup}\,||Hess\Lambda_{\t}(Y_l)||\Big)+\underset{k=1,...,n}{\max}\,|\hat{g}(X_k)-g(X_k)|
\\[0,2cm]&=O_P\bigg(\frac{1}{\sqrt{n}}\bigg)+O_P\big(||\hat{\t}-\t_0||\big)+O_P\bigg(\underset{k=1,...,n}{\max}\,|\hat{g}(X_k)-g(X_k)|\bigg)
\\[0,2cm]&=o_P(1).
\end{align*}
Since $P(||\hat{\t}-\t_0||\leq\delta)\rightarrow1$, this can be used together with Lemma \ref{lemmaan} to obtain
\begin{align}
&\frac{\partial}{\partial x_1}p^*(\mathcal{T}^*((h^*)^{-1}(u)),x)\nonumber
\\[0,2cm]&=\frac{1}{n}\sum_{k=1}^n\frac{\partial}{\partial x_1}F_{\xi}\bigg(\frac{u-\hat{g}(x)-\hat{\varepsilon}_k+\frac{1}{n}\sum_{l=1}^n\hat{\varepsilon}_l}{a_n}\bigg)f_{X^*}(x)\nonumber
\\[0,2cm]&=\frac{1}{n}\sum_{k=1}^n\bigg[-\frac{f_{X^*}(x)}{a_n}\ell \bigg(\frac{u-\hat{g}(x)-\hat{\varepsilon}_k+\frac{1}{n}\sum_{l=1}^n\hat{\varepsilon}_l}{a_n}\bigg)\frac{\partial}{\partial x_1}\hat{g}(x)\nonumber
\\[0,2cm]&\quad+{F_{\xi}}\bigg(\frac{u-\hat{g}(x)-\hat{\varepsilon}_k+\frac{1}{n}\sum_{l=1}^n\hat{\varepsilon}_l}{a_n}\bigg)\frac{\partial}{\partial x_1}f_{X^*}(x)\bigg]\nonumber
\\[0,2cm]&=-\frac{f_{X^*}(x)\frac{\partial}{\partial x_1}\hat{g}(x)}{na_n}\sum_{k=1}^n\ell \bigg(\frac{u-\hat{g}(x)+\frac{1}{n}\sum_{l=1}^n\hat{\varepsilon}_l-\varepsilon_k(\t_0)}{a_n}\bigg)\nonumber
\\[0,2cm]&\quad+\frac{\frac{\partial}{\partial x_1}f_{X^*}(x)}{n}\sum_{k=1}^nF_{\xi}\bigg(\frac{u-\hat{g}(x)+\frac{1}{n}\sum_{l=1}^n\hat{\varepsilon}_l-\varepsilon_k(\t_0)}{a_n}\bigg)+o_P(1)\nonumber
\\[0,2cm]&=-f_{\varepsilon(\theta_0)}\bigg(u-\hat{g}(x)+\frac{1}{n}\sum_{l=1}^n\hat{\varepsilon}_l\bigg)f_{X^*}(x)\frac{\partial}{\partial x_1}\hat{g}(x)+F_{\varepsilon(\t_0)}\bigg(u-\hat{g}(x)+\frac{1}{n}\sum_{l=1}^n\hat{\varepsilon}_l\bigg)\frac{\partial}{\partial x_1}f_{X^*}(x)\nonumber
\\[0,2cm]&\quad+o_P(1)\nonumber
\\[0,2cm]&=-f_{\varepsilon(\theta_0)}(u-g(x))f_X(x)\frac{\partial}{\partial x_1}g(x)+F_{\varepsilon(\t_0)}(u-g(x))\frac{\partial}{\partial x_1}f_X(x)+o_P(1)\nonumber
\\[0,2cm]&=\frac{\partial}{\partial x_1}p(\mathcal{T}((h)^{-1}(u)),x)+o_P(1)\label{unifconvdevxp*}
\end{align}
uniformly with respect to $x\in\operatorname{supp}(v)$ and with respect to $u$ belonging to some compact set $\mathcal{K}$, where the third last equality again follows from Theorem 4 in \cite{Han2008}. The same reasoning for $\frac{\partial}{\partial u}p^*$ results in
$$\frac{\partial}{\partial u}p^*(\mathcal{T}^*((h^*)^{-1}(u)),x)=\frac{\partial}{\partial u}p(\mathcal{T}((h)^{-1}(u)),x)$$
uniformly in $(u,x)\in\mathcal{K}\times\operatorname{supp}(v)$. Similarly, one can show
$$\mathcal{T}^*((h^*)^{-1}(u))-\tilde{\mathcal{T}}(\tilde{h}^{-1}(u))=o_P(1)$$
as well as
$$\frac{\partial}{\partial u}\mathcal{T}^*((h^*)^{-1}(u))-\frac{\partial}{\partial u}\mathcal{T}(h^{-1}(u))=o_P(1)$$
uniformly on compact sets. Hence, after possibly adjusting the set of admissible values for $u$, (\ref{unifconvdevxp*}) leads to
\begin{align*}
\frac{\partial^2}{\partial u\partial x_1}p^*(u,x)&=\frac{\partial}{\partial x_1}\frac{\frac{\partial}{\partial u}p^*(\mathcal{T}^*((h^*)^{-1}(t)),x)\big|_{u=h^*((\mathcal{T}^*)^{-1}(u))}}{\frac{\partial}{\partial u}\mathcal{T}^*((h^*)^{-1}(t))\big|_{u=h^*((\mathcal{T}^*)^{-1}(u))}}
\\[0,2cm]&=\frac{\partial}{\partial x_1}\frac{f_{\varepsilon}(h^*((\mathcal{T}^*)^{-1}(u))-g(x))\cdot f_X(x)}{\frac{\partial}{\partial u}\mathcal{T}^*((h^*)^{-1}(t))\big|_{u=h^*((\mathcal{T}^*)^{-1}(u))}}+o_P(1)
\\[0,2cm]&=\frac{\partial}{\partial x_1}\frac{f_{\varepsilon}(h(\mathcal{T}^{-1}(u))-g(x))\cdot f_X(x)}{\frac{\partial}{\partial u}\mathcal{T}((h)^{-1}(t))\big|_{u=h((T)^{-1}(u))}}+o_P(1)
\\[0,2cm]&=\frac{\partial^2}{\partial u\partial x_1}p(u,x)+o_P(1)
\end{align*}
uniformly on $(u,x)\in\mathcal{C}\times\operatorname{supp}(v)$.
\hfill $\Box$\\[0,2cm]

\begin{rem}
	Roughly speaking, the proof of Lemma \ref{lemmabootstrap} was based on the convergence of $\psi^*$ to $\psi$. If the alternative holds, it is not even clear if $\psi^*$ stabilizes in some sense (see Assumption \ref{boundedpsi*}). Hence, additional assumptions are needed. For that purpose define
	\begin{align*}
	F_{\varepsilon(\theta)}(e)&=P(\varepsilon(\theta)\leq e),
	\\[0,2cm]F_S^B(u)&=\int F_{\varepsilon(\theta_0)}(u-g_{\t_0}(x))f_X(x)\,dx,
	\\[0,2cm]\mathcal{T}_S^B(u)&=\frac{F_S^B(u)-F_S^B(0)}{F_S^B(1)-F_S^B(0)},
	\\[0,2cm]\tilde{\Phi}(u|x)&=F_{\varepsilon(\theta_0)}\big((\mathcal{T}_S^B)^{-1}(u)-g_{\theta_0}(x)\big).
	\end{align*}
	While doing so, assume $F_S^B(0)<F_S^B(1)$ to ensure that $\mathcal{T}_S^B$ is well defined, and define
	$$(\mathcal{T}_S^B)^{-1}(u)=\left\{\begin{array}{r}-\infty\\\infty\end{array}\right\},\quad\textup{if}\quad T_S^B(y)\left\{\begin{array}{r}>\\<\end{array}\right\}u\quad\textup{for all }y\in\mathbb{R}.$$
	$\tilde{\Phi}$ plays a similar role under the alternative as $\Phi(u|x)=F_{U|X}(u|x)$ under the null hypothesis and thus needs to be continuously differentiable on $\mathcal{U}_0\times\operatorname{supp}(v)$ with (again, the same $i$ as in \ref{B3} is used)
	\begin{equation}\label{devsPhitilde}
	\underset{(u,x)\in\mathcal{U}_0\times\operatorname{supp}(v)}{\inf}\,\frac{\partial}{\partial u}\tilde{\Phi}(u|x)>0\quad\textup{and}\quad\underset{(u,x)\in\mathcal{U}_0\times\operatorname{supp}(v)}{\inf}\,\frac{\partial}{\partial x_1}\tilde{\Phi}(u|x)>0.
	\end{equation}
\end{rem}

\begin{lemma}\label{lemmabootstrapalternative}
	Let Assumptions \ref{B11*} and (\ref{convghat}) be fulfilled. Further, assume $H_1$,\ref{A1}--\ref{A4},\ref{A6'},\ref{A7},\ref{A8'},\ref{A10}, \ref{B1}--\ref{B10},(\ref{unifconvghat}), (\ref{boundedepsfeps}) and (\ref{devsPhitilde}). Further, assume the existence of a neighbourhood $\tilde{\Theta}$ of $\theta_0$ such that the map $y\mapsto\Lambda_{\theta}(y)$ is $(q+2)$-times continuously differentiable for all $\theta\in\tilde{\Theta}$. Let $\hat{h}^*$ be the estimator from (\ref{hath}) based on the bootstrap data $(Y_j^*,X_j^*),j=1,...,m$. Then, assumption \ref{A8*} is fulfilled.
\end{lemma}

\bigskip

\noindent {\bf Proof of Lemma \ref{lemmabootstrapalternative}.}
Only equation \ref{boundedpsi*} needs to be proven, since the remaining conditions follow as in the proof of Lemma \ref{lemmabootstrap} by the results of \citet{CvK2018}. In contrast to the proof of Lemma \ref{lemmabootstrap}, it is $\hat{\t}-\t_0=O_P\big(n^{-\frac{1}{4}}\big)$ (here and in the following, $o_P$- and $O_P$-terms are with respect to $P_1$ and for $n\rightarrow\infty$) and the asymptotic behaviour of $\psi^*$ can not be reduced to the convergence to $\psi$. Nevertheless, $\psi^*$ can be expressed as in (\ref{defpsi^*}) with $D_{p,0}^*,...,D_{f,1}^*$ as in (\ref{Dp0*})--(\ref{Dfi*}). The main idea is to prove uniform convergence of $\frac{\partial
}{\partial u}\Phi^*$ and $\frac{\partial
}{\partial u}\Phi^*$ on $\mathcal{U}_0\times\operatorname{supp}(v)$ to $\frac{\partial
}{\partial u}\tilde{\Phi}$ and $\frac{\partial
}{\partial x_1}\tilde{\Phi}$, respectively, while the remaining parts of ${\delta_j^*}^{\tilde{v}_1^*},{\delta_j^*}^{\tilde{v}_2^*}$ and ${Q^*}'$ are bounded in probability.\\
Due to (\ref{devsPhitilde}) it is $\mathcal{U}_0\subseteq\big(-\frac{F_S^B(0)}{F_S^B(1)-F_S^B(0)},\frac{1-F_S^B(0)}{F_S^B(1)-F_S^B(0)}\big)$. In the following it is proven that under the assumptions of Lemma \ref{lemmabootstrapalternative} one has
\begin{align}
\underset{t\in\mathbb{R}}{\sup}\,|\mathcal{T}_{S^*}(t)-\mathcal{T}_S^B(t)|&=o_P(1),\label{unifconvTS*}		\\[0,2cm]\underset{u\in\mathcal{U}_0}{\sup}\,|(\mathcal{T}_{S^*})^{-1}(u)-(\mathcal{T}_S^B)^{-1}(u)|&=o_P(1),\label{unifconvTS*^-1}
\\[0,2cm]\underset{(u,x)\in\mathcal{U}_0,\times\operatorname{supp}(v)}{\sup}\,|\Phi^*(u,x)-\tilde{\Phi}(u,x)|&=o_P(1),
\\[0,2cm]\underset{(u,x)\in\mathcal{U}_0,\times\operatorname{supp}(v)}{\sup}\,\bigg|\frac{\partial}{\partial u}\Phi^*(u,x)-\frac{\partial}{\partial u}\tilde{\Phi}(u,x)\bigg|&=o_P(1),
\\[0,2cm]\underset{(u,x)\in\mathcal{U}_0,\times\operatorname{supp}(v)}{\sup}\,\bigg|\frac{\partial}{\partial x_1}\Phi^*(u,x)-\frac{\partial}{\partial x_1}\tilde{\Phi}(u,x)\bigg|&=o_P(1).
\end{align}
Due to Lemma \ref{lemmaan} $F_{S^*}$ can be written for appropriate $u_{i,w}^*\in\mathbb{R},i=1,...,n,\ w\in\operatorname{supp}(\kappa ),$ as
\begin{align*}
F_{S^*}(u)&=\frac{1}{n^2}\sum_{i=1}^n\sum_{k=1}^n\int F_{\xi}\Bigg(\frac{u-\hat{g}(X_i+b_nw)-\hat{\varepsilon}_k+\frac{1}{n}\sum_{l=1}^n\hat{\varepsilon}_l}{a_n}\Bigg)\kappa (w)\,dw
\\[0,2cm]&=\frac{1}{n^2}\sum_{i=1}^n\sum_{k=1}^n\int F_{\xi}\Bigg(\frac{u-g_{\t_0}(X_i+b_nw)-\hat{\varepsilon}_k}{a_n}\Bigg)\kappa (w)\,dw
\\[0,2cm]&\quad+\frac{1}{n^2}\sum_{i=1}^n\sum_{k=1}^n\int \ell (u_{i,w}^*)\kappa (w)\bigg(\frac{\hat{g}(X_i+b_nw)-g_{\t_0}(X_i+b_nw)+\frac{1}{n}\sum_{l=1}^n\hat{\varepsilon}_l}{a_n}\bigg)\,dw
\\[0,2cm]&=\frac{1}{n^2}\sum_{i=1}^n\sum_{k=1}^n\int F_{\xi}\bigg(\frac{u-g_{\theta_0}(X_i+b_nw)-\varepsilon_k(\t_0)}{a_n}\bigg)\kappa (w)\,dw+o_P(1).
\end{align*}
As a distribution function $F_{\xi}$ is bounded so that
$$\operatorname{Var}\Bigg(\frac{1}{n^2}\sum_{i=1}^n\sum_{k=1}^n\int F_{\xi}\bigg(\frac{u-g_{\theta_0}(X_i+b_nw)-\varepsilon_k(\t_0)}{a_n}\bigg)\kappa (w)\,dw\Bigg)\rightarrow0,$$
that is
\begin{align*}
F_{S^*}(u)&=E\Bigg[\int F_{\xi}\bigg(\frac{u-g_{\theta_0}(X_1+b_nw)-\varepsilon_2(\t_0)}{a_n}\bigg)\kappa (w)\,dw\Bigg]+o_P(1)
\\[0,2cm]&=\int\int E\Bigg[F_{\xi}\bigg(\frac{u-g_{\theta_0}(x+b_nw)-\varepsilon_2(\t_0)}{a_n}\bigg)\Bigg]\kappa (w)f_X(x)\,dx\,dw+o_P(1)
\\[0,2cm]&=\int\int E\big[I_{\{u-g_{\theta_0}(x+b_nw)-\varepsilon_2(\t_0)\geq0\}}\big]\kappa (w)f_X(x)\,dx\,dw+o_P(1)
\\[0,2cm]&=\int\int F_{\varepsilon(\t_0)}(u-g_{\theta_0}(x+b_nw)\big)\kappa (w)f_X(x)\,dw\,dx+o_P(1)
\\[0,2cm]&=\int\int F_{\varepsilon(\t_0)}(u-g_{\t_0}(x)\big)\kappa (w)f_X(x)\,dw\,dx+o_P(1).
\end{align*}
Since $F_{S^*}$ and $F_S^B$ are distribution functions, this leads to the uniform convergence
$$\underset{t\in\mathbb{R}}{\sup}\,|F_{S^*}(t)-F_S^B(t)|=o_P(1)$$
and thus to (\ref{unifconvTS*}). To prove (\ref{unifconvTS*^-1}) write $F_{\varepsilon(\theta_0)}$ as
\begin{align}
F_{\varepsilon(\theta_0)}(e)&=P(\varepsilon(\theta_0)\leq e)\nonumber
\\[0,2cm]&=P\big(g(X)+\varepsilon\leq h(h_{\theta_0}^{-1}(e+g_{\theta_0}(X)))\big)\nonumber
\\[0,2cm]&=\int F_{\varepsilon}\big(h(h_{\theta_0}^{-1}(e+g_{\theta_0}(x)))-g(x)\big)f_X(x)\,dx,\nonumber
\end{align}
which implies
$$F_S^B(u)=\int\int F_{\varepsilon}\big(h(h_{\theta_0}^{-1}(u-g_{\theta_0}(z)+g_{\theta_0}(x)))-g(x)\big)f_X(x)f_X(z)\,dx\,dz.$$
for
$$h_{\theta_0}^{-1}(u)=\left\{\begin{array}{r}-\infty\\\infty\end{array}\right\},\quad\textup{if}\quad h_{\theta_0}(y)\left\{\begin{array}{r}>\\<\end{array}\right\}u\quad\textup{for all }y\in\mathbb{R}.$$
Since $h$ and $h_{\theta_0}$ are strictly increasing and $F_S^B$ is continuous, it is $F_S^B(u_1)<F_S^B(u_2)$ for all $u_1<u_2\in (\mathcal{T}_S^B)^{-1}(\mathcal{U}_0)\subseteq (\mathcal{T}_S^B)^{-1}\Big(\big(-\frac{F_S^B(0)}{F_S^B(1)-F_S^B(0)},\frac{1-F_S^B(0)}{F_S^B(1)-F_S^B(0)}\big)\Big)$. Especially, $(\mathcal{T}_S^B)^{-1}$ is strictly increasing on $\mathcal{U}_0$, that is, (\ref{unifconvTS*^-1}) follows from (\ref{unifconvTS*}).\\
Finally, this can be used to obtain
\begin{align*}
\Phi^*_i(u|x)&=\frac{\partial}{\partial x_1}P^*(U^*_j\leq u|X^*_j=x)
\\[0,2cm]&=\frac{1}{a_nn}\sum_{k=1}^n\ell \Bigg(\frac{\mathcal{T}_{S^*}^{-1}(u)-\hat{g}(x)+\frac{1}{n}\sum_{l=1}^n\hat{\varepsilon}_l-\hat{\varepsilon}_{k}}{a_n}\Bigg)\frac{\partial}{\partial x_1}\hat{g}(x)
\\[0,2cm]&=-\frac{\partial}{\partial x_1}\hat{g}(x)\frac{1}{na_n}\sum_{k=1}^n\ell \Bigg(\frac{\mathcal{T}_{S^*}^{-1}(u)-\hat{g}(x)+\frac{1}{n}\sum_{l=1}^n\hat{\varepsilon}_l-\varepsilon_k(\theta_0)}{a_n}\Bigg)+o_P(1)
\\[0,2cm]&=-\frac{\frac{\partial}{\partial x_1}\hat{g}(x)}{a_n}\int \ell \Bigg(\frac{\mathcal{T}_{S^*}^{-1}(u)-\hat{g}(x)+\frac{1}{n}\sum_{l=1}^n\hat{\varepsilon}_l-e}{a_n}\Bigg)f_{\varepsilon(\theta_0)}(e)\,de+o_P(1)
\\[0,2cm]&=-\frac{\partial}{\partial x_1}\hat{g}(x)\int \ell (e)f_{\varepsilon(\theta_0)}\bigg(\mathcal{T}_{S^*}^{-1}(u)-\hat{g}(x)+\frac{1}{n}\sum_{l=1}^n\hat{\varepsilon}_l\bigg)\,de+o_P(1)
\\[0,2cm]&=-\frac{\partial}{\partial x_1}g_{\theta_0}(x)\int \ell (e)f_{\varepsilon(\theta_0)}\big((\mathcal{T}_{S}^B)^{-1}(u)-g_{\theta_0}(x)\big)\,de+o_P(1)
\\[0,2cm]&=-\frac{\partial}{\partial x_1}g_{\theta_0}(x)f_{\varepsilon(\theta_0)}\big((\mathcal{T}_{S}^B)^{-1}(u)-g_{\theta_0}(x)\big)+o_P(1)
\end{align*}
uniformly in $(u,x)\in\mathcal{U}_0\times\operatorname{supp}(v)$, where the second last equality follows from the continuity of $f_{\varepsilon(\theta_0)}$. The bootstrap functions
\begin{align*}
\frac{\partial}{\partial u}F_{S^*}(u)&=\frac{1}{n^2a_n}\sum_{i=1}^n\sum_{k=1}^n\int \ell \Bigg(\frac{u-\hat{g}(X_i+b_nw)-\hat{\varepsilon}_k+\frac{1}{n}\sum_{l=1}^n\hat{\varepsilon}_l}{a_n}\Bigg)\kappa (w)\,dw,
\end{align*}
$\Phi^*$ and $\frac{\partial}{\partial u}\Phi^*$ can be treated by similar arguments to obtain
$$\underset{(u,x)\in\mathcal{U}_0\times\operatorname{supp}(v)}{\sup}\,|\Phi^*(u|x)-\tilde{\Phi}(u|x)|+\bigg|\frac{\partial}{\partial x_1}\Phi^*(u|x)-\frac{\partial}{\partial x_1}\tilde{\Phi}(u|x)\bigg|=o_P(1).$$
Since $Q^*=\mathcal{T}_{S^*}^{-1}$ equations (\ref{unifconvTS*^-1}) and (\ref{Dp0*})--(\ref{Dfi*}) lead to
$$\underset{y\in\mathbb{R},z\in\mathbb{R}^{d_X+1}}{\sup}\,|w(y)\psi^*(z,\mathcal{T}^*(y))|=O_P(1).$$
\hfill $\Box$

\end{document}